\documentclass[pre,floatfix,superscriptaddress,showpacs]{revtex4}%
\usepackage{amsmath,amsthm,amssymb,amsbsy,graphicx,epsfig,citesort}

\pdfoutput=1

\begin{document}

\newtheorem{thm}{Theorem}[section]
\newtheorem{cor}[thm]{Corollary}
\newtheorem{lem}[thm]{Lemma}
\newtheorem{definition}[thm]{Definition}
\newtheorem{note}[thm]{Note}
\newtheorem{conjecture}[thm]{Conjecture}

\def\Bu{\mathbf{u}}
\def\Bv{\mathbf{v}}
\def\Bw{\mathbf{w}}
\def\Bomega{\mathbf{\omega}}
\def\Bn{\mathbf{n}}
\def\Btau{\mathbf{\tau}}
\def\Br{\mathbf{r}}
\def\Bf{\mathbf{f}}
\def\BR{\mathbf{R}}
\def\Bx{\mathbf{x}}
\def\By{\mathbf{y}}
\def\Beta{\mathbf{\eta}}
\def\Bxi{{\mathbf{\xi}}}
\def\T{\mathcal{T}}
\def\BF{\mathbf{F}}
\def\BT{\mathbf{T}}
\def\BI{\mathbf{I}}
\def\Re{\mathcal{R}}
\def\BU{\mathbf{U}}
\def\BW{\mathbf{W}}

\def\lap{\triangle}

\def\tphi{\tilde{\phi}}
\def\ve{\varepsilon}
\def\a{\alpha}
\def\b{\beta}
\def\g{\gamma}
\def\d{\partial}
\def\l{\lambda}
\def\s{\sigma}
\def\D{\nabla}
\def\R{\mathbb{R}}
\def\W {\Omega}
\def\w{\omega}
\def\lap{\triangle}
\def\trace{\operatorname{trace}}
\def\curl {\operatorname{curl}}
\def\Div{\operatorname{div}}
\def\degr{\operatorname{deg}}
\def\diam{\operatorname{diam}}
\def\diiv{\operatorname{div}}
\def\aint{\,-\hspace{-4.364mm}\int}
\def\epi{\operatorname{epi}}
\def\esssup{\operatorname{ess}\sup}
\def\argmin{\operatorname{arg}\min}
\def\sign{\operatorname{sign}}

\title {A model of hydrodynamic interaction between swimming bacteria}

\author{Vitaliy Gyrya}
\affiliation{Department of Mathematics, Pennsylvania State
University, 418 McAllister Building, University Park, PA 16802 }

\author{Igor S.~Aranson}
\affiliation{Materials Science Division, Argonne National
Laboratory, 9700 South Cass Avenue, Argonne, IL 60439}

\author{Leonid V.~Berlyand}
\affiliation{Department of Mathematics, Pennsylvania State
University, 337 McAllister Building, University Park, PA 16802}

\author{Dmitry Karpeev}
\affiliation{Mathematics and Computer Science Division, Argonne National
Laboratory, 9700 South Cass Avenue, Argonne, IL 60439}

\date{\today}

\begin{abstract}
We study the dynamics and interaction of two swimming bacteria, modeled by
self-propelled dumbbell-type  structures. We focus on alignment dynamics of a coplanar pair
of elongated swimmers, which propel themselves either by ``pushing'' or ``pulling''   both in
three- and quasi-two-dimensional
geometries of space.
We derive asymptotic expressions for the dynamics
of the pair, which, complemented by numerical experiments, indicate that
the tendency of bacteria to swim in or swim off depends strongly on the position of the
propulsion force.
In particular, we observe that positioning of the effective propulsion force inside the dumbbell
results in qualitative agreement with the dynamics observed in experiments, such as
mutual alignment of converging bacteria.
\end{abstract}

\pacs{87.16.-b, 05.65.+b, 87.17.Jj}

\maketitle

\section{Introduction}

Modeling of bacterial suspensions and, more generally, of suspensions of active microparticles
has recently become an increasingly active area of research.
One of the motivating factors
behind this trend is the study of the dynamics of large populations of aquatic single-cellular
\cite{WuLib00,MenBouWilAndWat99,DomCisChaGolKes04,Kes00,SokAraKesGol07} and  multicellular organisms \cite{Short06,Kitsunezaki07}.
In particular, there is significant interest in understanding the mechanism of formation
of coherent structures on a scale much larger than individual microorganisms in suspension
(see, e.g., \cite{PedKes92,SokAraKesGol07,AraSokGolKes07}).
Recently, bacterial suspensions have also emerged as a prototypical system for the study and engineering of novel
biomaterials with unusual rheological properties \cite{KimBre04,Khatavkar07}.
Here the idea is to exploit the active nature of the particles in the suspension in order to generate specific effects,
such as enhancement of transport and diffusion of tracers  relative to that of
the solvent \cite{KimBre04,Short06,SokAraKesGol07}.
A good review of the motivations, experimental studies, and modeling approaches to suspensions of swimming microorganisms is contained in the introduction to \cite{IshiSimPed06}.

The principal organizing role in the formation of large-scale patterns (e.g., \cite{MenBouWilAndWat99,GreCha04,DomCisChaGolKes04,SokAraKesGol07})
is believed to be played by hydrodynamic interactions between individual swimmers and the environment.
This includes the boundary effects as well as the hydrodynamic interaction with other swimmers
\cite{PedKes92,Kes00}.
These effects and interactions are also believed to set the spatial and temporal scales of the patterns.
At the same time, fundamental questions about the hydrodynamics of a single swimmer have
been studied by many researchers over several decades
(e.g., \cite{Taylor51,Lighthill76,BreWin77,NajGol04,Cisneros07} and references therein).
Here one of the central features is the very low Reynolds number $\mathcal{R}e$ of a typical microscopic swimmer --
$\mathcal{R}e \sim 10^{-4}$ -- $10^{-2}$ \cite{IshiSimPed06,BreWin77} -- making the governing dynamics Stokesian.
Since the Stokesian dynamics is time-reversible, the very possibility of propulsion at low
Reynolds numbers had to be clarified in general (see, e.g., \cite{NajGol04}), with some
of the early important contributions made by Purcell \cite{Purcell77}.
Specific studies of the propulsion mechanism of flagellates includes the work by Phan-Thien et al.
\cite{NasPhThien97,RamTulPhThien93}.
In particular, Ref.~\cite{NasPhThien97} in detail addresses the hydrodynamic interactions of two nearby microswimmers.
We also address the question of pairwise hydrodynamic interactions of swimmers.
Unlike in \cite{NasPhThien97}, however, our model abstracts from the method of propulsion
(e.g., rotating helix, water jet) and is applicable to a wider class of swimmers.
Also, being simpler structurally, our model allows us to perform simulations for
a larger collection of swimmers and relatively long time.


Studies of the fundamental interactions of small numbers (e.g., pairs) of swimming particles are
also important for validating mean-field theories of large-scale pattern formation
in active suspensions \cite{SimRama02}.  The continuum phenomenological  models proposed in such studies typically
rely on a two-phase formulation of the problem: the particle phase interacts with the fluid phase
via a postulated coupling mechanism.  It should be possible, at least in principle, to derive or
verify the proposed coupling mechanisms against the fundamental particle-particle
dynamics.
In this paper we focus on  the mechanisms of alignment
of a pair of elongated swimmers, with the aim of shedding light on a possible mechanism of
large-scale ordering in dilute bacterial suspensions.

The main emerging approaches to the modeling of microswimmers
\cite{SainShel07,IshiSimPed06,HerStoGra05,Cisneros07}
typically abstract away the details of the actual propulsion mechanism and use simple,
tractable, rigid geometries to model a swimmer.
In general ``higher-order'' effects such as signaling between bacteria and
chemotaxis are ignored, the emphasis being on the basic hydrodynamic interactions.
Both \cite{SainShel07} and \cite{HerStoGra05} model swimmers as elongated bodies; the former
employs slender body theory for cylindrical rods, while the latter models the elongated
body as a \emph{dumbbell} consisting of a pair of balls.
As our model is a modification of
\cite{HerStoGra05}, the following sections contain a more detailed description of
the dumbbell model.
In the aforementioned studies the self-propulsion mechanism is modeled by a prescribed \emph{force},
concentrated at a point inside the dumbbell ball, as in \cite{HerStoGra05}, or distributed over a part
of the surface of the body in the form of a specified tangential traction, as in \cite{SainShel07,Cisneros07}.
Both of these studies ultimately rely on numerical simulations with the goal of studying
the emergence of large-scale coherent patterns predicted by continuum theories such as \cite{SimRama02,AraSokGolKes07}
or observed in experiments \cite{DomCisChaGolKes04,SokAraKesGol07}.
By contrast, the main tool of our work  is an asymptotic analysis followed by straightforward numerical simulations.

A different structural and dynamic approach is taken by Pedley and coworkers (see \cite{IshiSimPed06,IshiPed07}).  Here a basic swimmer is modeled as a squirming sphere, with a prescribed tangential
\emph{velocity} as the model of the propulsion by motile cilia
(short hair on the surface of the cell beating in the same direction).
As with the dumbbell model, a spherical squirmer
allows using some fundamental solutions and relations (e.g., the Stokeslet solution and the Fax\'en relations)
to approximate the dynamics of the swimmers.
The work \cite{IshiSimPed06} is closer to ours in its goal
of quantifying the interaction of a \emph{pair} of swimmers, rather than a large collection of swimmers,
while the model in \cite{HerStoGra05} is closer to ours in the structural model of a swimmer.

In our work, self-propulsion is modeled by prescribed propulsion forces
(as in \cite{HerStoGra05} vs. prescribed velocities on the boundary, as in \cite{IshiSimPed06};
see also recent work \cite{HaiAraBerKar08} on the rheology of bacterial suspensions),
and an elongated body of bacterium is modeled by a dumbbell as in \cite{HerStoGra05}.
Our model is consistently derived from the equations of Stokesian fluid dynamics.
In particular, we model self-propulsion by a point force, whose location can vary.
This allows us to investigate the dependence of mutual dynamics (swim in/off)
of neighboring bacteria on the position of this force (which roughly can be interpreted as the effect of the shape of the microorganism and the way of  propulsion or distribution of cilia).
We study the hydrodynamic interaction for well-separated swimmers.
For this reason we say that two swimmers \emph{swim in}, starting from a given mutual orientation,
if at some point the distance between them decreases to the order of their size,
(i.e., they become not well separated).
Two swimmers, starting from a given mutual orientation, \emph{swim off} if the distance
between them increases to infinity without \emph{swim in} happening first.
We observe that positioning the propulsion force between the dumbbell balls results in
attractive behavior of swimmers.
On the other hand, positioning the propulsion force outside the dumbbell results in
repulsive behavior.
For comparison, in the earlier work \cite{HerStoGra05} the position of the propulsion force was fixed
(center of a ball in a dumbbell).


In this work we study the \emph{alignment} of a coplanar pair of three-dimensional
elongated swimmers, which propel themselves by ``pushing'' or, ``pulling'',
mimicking
a variety of self-propelled microorganisms, from sperm cells and bacteria to algae.
We derive asymptotic expressions for the dynamics
of the pair, which, complemented by numerical experiments, indicate that
the tendency of bacteria to swim in or swim off strongly depends on the position of the
propulsion force.
In particular, we observe that positioning of the propulsion force inside the dumbbell
results in the qualitative agreement with the dynamics observed in experiments \cite{AraSokGolKes07}.
We also observe that the dynamics of bacteria in a thin film
(with no-slip boundary conditions on the top and bottom) is qualitatively similar  to that for the whole space.
%
%

One of our objectives
is to develop a well-posed PDE model of an active suspension
derived from first principles (unlike many engineering models that use ad hoc assumptions).
Our proposed model is simple enough to allow for theoretical analysis (asymptotics)
yet captures basic features observed in experimental studies.

The paper is organized as follows.
In section \ref{sect:Model} we derive a full PDE model for the dynamics of swimmers
based on Stokesian hydrodynamics.
The well-posedness of the problem was demonstrated and can be found in
Appendix \ref{app:existence}.

In section \ref{sect:Asymptotic reduction} we introduce an asymptotic reduction of the PDE
model in the dilute limit of swimmer concentration.
Here we show how to solve the reduced model numerically and give an analytic (asymptotic)
solution for a pair of swimmers.
Then we analyze two basic (physically interesting) configurations of swimmers
based on the asymptotic formulas and numerical calculations.
In section \ref{sect:Conclusions} we make some concluding remarks and indicate areas for
future study.
Appendix \ref{sect:Stokes} contains the fundamental solutions to the Stokes problem,
which are used extensively in section \ref{sect:Asymptotic reduction}.
Appendix \ref{sect:formulas} contains some technical asymptotic formulas and calculations.
Appendix \ref{subsect:stability of MI} contains the stability analysis of certain
configuration of bacteria.

\section{The model}\label{sect:Model}

While the modeling of suspensions of passive particles is a well-established area,
the mathematical study of suspensions of active swimmers is still an open area without
universally accepted models that can serve as benchmarks for analytic and numerical
studies.
In this paper we propose a model that can serve as
a tractable reference case for the mathematical questions
such as existence and uniqueness of solutions.
We analyze this model by considering its
asymptotic reduction in the far-field regime and comparing its predictions
with experimental results for bacterial suspensions.
The main goal is to establish a model amenable to an analytic treatment that
still captures the main physical effects, such as alignment and
the emergence of large-scale coherent structures; here we focus on the question
of the pairwise dynamics of swimmers.

We set up a PDE model for a collection of
swimmers and then consider an asymptotic ODE reduction for a pair of swimmers that is suitable for numerical
analysis.
We address the question of the solvability of the PDE system and carry out a numerical
study of the ODE model.


\subsection{Structure and dynamics of a single swimmer}

\begin{figure}[h!]
\begin{center}
  \includegraphics[width=9cm]{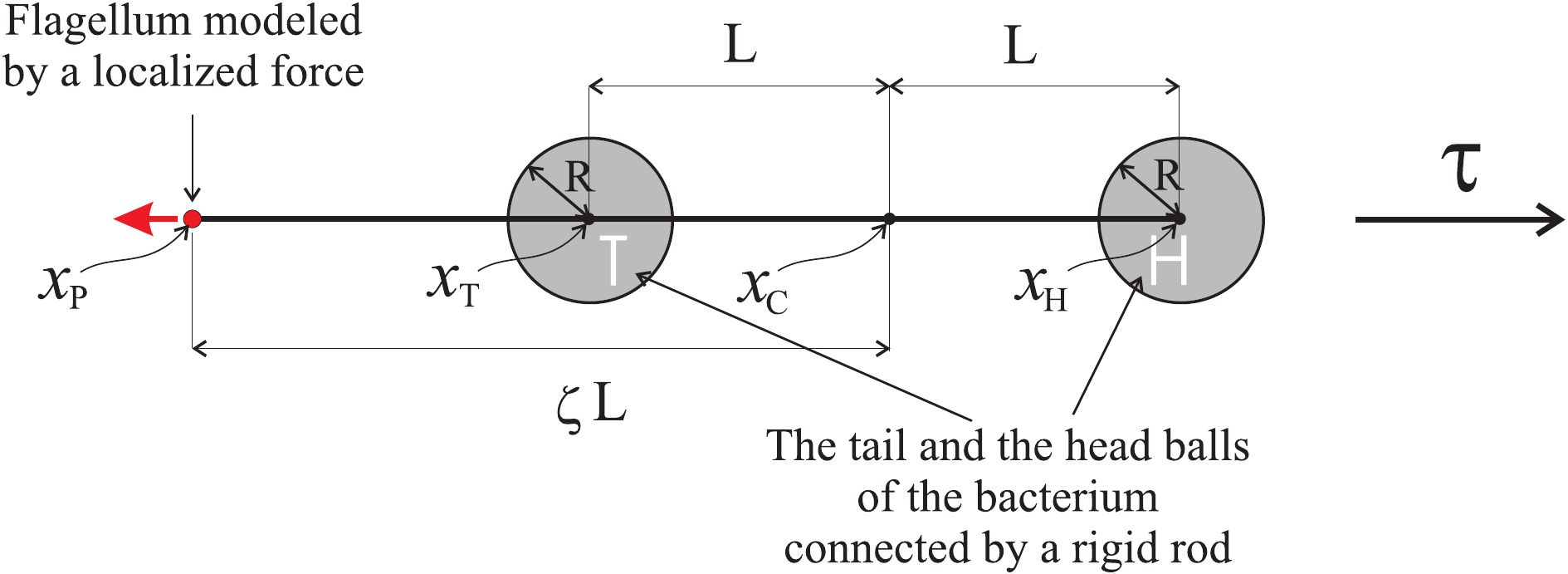}\\
  \caption{Model of a single bacterium:
  two balls (labeled head and tail) and the propulsion force (red ball with arrow) connected by a rigid rod (that does not interact with the fluid).}\label{fig:notation}
\end{center}
\end{figure}

Structurally, we model as a \emph{dumbbell} (see Fig.  \ref{fig:notation})
the elongated body of a bacterium (e.g., \emph{Bacillus subtilis}), referred to as the swimmer below.
Similar approximation  was used in \cite{HerStoGra05}.
It consists of two balls of equal mass $m$ and radius $R$ rigidly connected to one another a distance $2L$ apart.
Assuming a high aspect ratio of a swimmer, we have $2L \gg R$.
The balls are denoted $B_{_H}$ (head) and $B_{_T}$ (tail), with their centers located
at $\Bx_{_H}$ and $\Bx_{_T}$, respectively.
The unit vector directed from $\Bx_{_T}$ to $\Bx_{_H}$ is denoted by
\mbox{$\Btau = (\Bx_{_H} - \Bx_{_T})/|\Bx_{_H} - \Bx_{_T}|$}, indicating the direction of the swimmer's motion, and
the line connecting the centers of the balls is referred to as the {\em dumbbell axis}.
The action of the flagellum -- the  bacterial propulsion apparatus --
is represented either by a smooth volume force density $\BF$ supported in a
ball $B_P$ of small radius $\varrho$ ($\varrho \ll 1$), with the center
$\Bx_{_P}$ located on the dumbbell axis, or by a delta function concentrated at $\Bx_{_P}$.
This force and its support will sometimes be referred to as the \emph{propulsion force}.
The location $\Bx_{_P}$ of the propulsion force relative to the positions
$\Bx_{_H}$ and $\Bx_{_T}$ of the balls in the dumbbell is defined by the parameter
\begin{equation}\label{eq:zeta}
    \zeta = \frac{(\Bx_{_P}-\Bx_{_C})\cdot \tau}{L}, \qquad
    \text{where }\ \Bx_{_C} = \frac{\Bx_{_H}+\Bx_{_T}}{2}.
\end{equation}
For instance $\zeta=\mp1$ corresponds to the the center of the tail or head ball, respectively, and
$\zeta=0$ corresponds to the the center of the dumbbell.

Depending on the value of $\zeta$, swimmers are classified into pushers/pullers and inner/outer swimmers.
A swimmer is called a \emph{pusher} if $\zeta<0$ and \emph{puller} if $\zeta>0$.
A swimmer is called \emph{outer} if $|\zeta|>1$ and \emph{inner} if $|\zeta|<1$.

The total force exerted in the ball $B_{_P}$ has magnitude $f_p = \texttt{const}$
and is directed along the axis $\Btau$, that is,
\begin{equation}
        \int_{B_P} \BF(\Bx) d\Bx = -\BF_{_P}= -f_p \Btau.
        \label{eq:total_force}
\end{equation}
Here $-\BF_{_P}$ is the force of the flagellum on the fluid, and
$\BF_{_P}$ is the force of the fluid onto the flagellum, that is the force that propels the dumbbell.
For simplicity of presentation, we assume that the propulsion force is a delta function
$-\delta(\Bx-\Bx_{_P}) \BF_{_P}$.

We now discuss the dynamics of a single swimmer immersed in a fluid.
The kinematic constraints resulting from rigid connections between the balls and the location of the force can be
implemented mathematically by the equations
\begin{gather}
        \left(\Bv_{_H} - \Bv_{_T} \right) \cdot \Btau = \left(\Bv_{_P} - \Bv_{_T}\right)\cdot \Btau = 0
        \label{eq:swimmer_linear_rigidity}\\
        \omega_{_H} = \omega_{_T} = \omega_{_C}=\omega,
        \label{eq:swimmer_angular_rigidity}
\end{gather}
where $\Bv_{_H}$, $\Bv_{_T}$, and $\Bv_{_P}$ are the linear velocities of the head, the tail, and the force.
The second constraint \eqref{eq:swimmer_angular_rigidity} expresses the assumption that the balls
do not rotate with respect to the axis of the swimmer: $\omega_{_H}$, $\omega_{_T}$, and $\omega_{_C}$ are the
angular velocities of the balls and the axis respectively (hence we use the notation $\omega$ without a subscript).
The constraints (\ref{eq:swimmer_linear_rigidity}) and (\ref{eq:swimmer_angular_rigidity})
can be thought of as implemented by a
rigid rod (connecting $\Bx_{_H},\Bx_{_T}$, and $\Bx_{_P}$) of negligible thickness and mass, hence of negligible drag and inertia.

From this point on we will use the symbol $*$ to indicate $H,T$, or $P$.
When the meaning of $*$ can be ambiguous, we will explicitly mention the values it is allowed to take.

The motion of a point $\Bx$ on the surface of the ball $B_*$ ($*=H,T$) can be described in
two equivalent forms:
\begin{equation}\label{eq:surface velocity}
    \Bv(\Bx) = \Bv_{_C}+(\Bx-\Bx_{_C})\omega
    \qquad
    \text{or}
    \qquad
    \Bv(\Bx) = \Bv_{_*}+(\Bx-\Bx_{_*})\omega.
\end{equation}
Under different circumstances it is convenient to use one or another of these forms.
The connection between two forms is given by
\begin{equation}\label{eq:connection between forms of v}
    \Bv_{_H} = \Bv_{_C}-L\tau\times\w,\qquad
    \Bv_{_T} = \Bv_{_C}+L\tau\times\w.
\end{equation}

Note that any values of $\Bv_{_C}$ and $\omega$ define a rigid motion of a swimmer.
However, $\Bv_{_H}$, $\Bv_{_T}$, and $\omega$ must satisfy the additional
rigidity constrains (\ref{eq:rigidity},\ref{eq:rotational consistency}) to define a rigid motion.
The first rigidity constraint is the distance between $B_H$ and $B_T$ balls being preserved:
\begin{equation}\label{eq:rigidity}
    \tau\cdot(\Bv_{_H}-\Bv_{_T}) = 0.
\end{equation}
The second rigidity constraint is the consistency of the rotation
defined by $\Bv_{_H}$, $\Bv_{_T}$, and $\omega$:
\begin{equation}\label{eq:rotational consistency}
    (\Bv_{_H}-\Bv_{_T}) = -2L \tau\times \omega.
\end{equation}
The second rigidity constraint shows that the linear velocities $\Bv_{_H}$ and $\Bv_{_T}$
of the balls $B_{_H}$ and $B_{_T}$
define the angular velocity $\w$ of a dumbbell
(up to a rotation around the dumbbell axis $\tau$).

From Newton's second law of motion
\begin{gather}
        m \dot \Bv_{_C} = \BF_H + \BF_T + \BF_p,
        \label{eq:swimmer_linear_dynamics}\\
        I \dot \omega = \BT_H + \BT_T.
        \label{eq:swimmer_angular_dynamics}
\end{gather}
Here $\Bv_{_C} = \dot \Bx_{_C}$ is the velocity of the center of mass $\Bx_{_C}$ of the dumbbell,
$I$ is the moment of inertia for the dumbbell with respect to $\Bx_{_C}$,
and $\BF_{_P} = f_p\,\Btau$ is the reaction force onto the dumbbell from the point force pushing on the fluid
(modeling the action of flagellum).
The forces $\BF_{_H}$ and $\BF_{_T}$ are
due to the viscous drag exerted onto the head and the tail by the fluid; likewise,
$\BT_{_H}$ and $\BT_{_T}$ are the torques
due to the viscous force (the applied force $\BF_{_P}$ acts along the dumbbell's axis and results in zero torque).
All forces are applied at the center of mass $\Bx_{_C}$ and all torques are calculated
with respect to $\Bx_{_C}$.
The forces $\BF_*$ and the torques $\BT_*$ ($*=H,\,T$) are given by
\begin{gather}\label{eq:def:forcesNtorques}
        \BF_* = \int_{\partial B_*} \sigma(\Bu) \cdot \hat \Bn(\Bx) dS(\Bx),\quad
        \BT_* = \int_{\partial B_*} \left(\Bx - \Bx_{_C}\right) \times \sigma(\Bu) \cdot \hat \Bn(\Bx) dS(\Bx),
\end{gather}
where the stress tensor $\sigma(\Bu)$ is defined in terms of the strain rate (symmetrized gradient) $D(\Bu)$,
$$
   \sigma(\Bu) = 2\,\mu D(\Bu) - p \BI,\quad 2\,D(\Bu) = \D \Bu + \left(\D \Bu\right)^T,
$$
and $\hat \Bn(\Bx)$ is the unit inward normal to $\partial B_*$.

\subsection{Discussion of the model}

Many bacteria swim by rotating the flagellum, driven by the torque generating motors located within the bacterial membrane.
The rotation of the flagellum causes the body of the bacterium to rotate in the opposite direction around its axis of symmetry.
This rotation exerts a hydrodynamic torque on the fluid.
Experiments for which this torque is important are those where bacteria are close to one another or
to the wall of container.
For example it is known \cite{DilTurMayGarWeiBerWhi05} that a bacterium swimming next to a solid wall will swim in circles in the plane parallel to the wall.
The direction of the swimming is determined by the chirality of the flagellum.

On the other hand, for well-separated swimmers the effects of the torque around the axis of symmetry can be neglected
because the disturbance due to it decays as $r^{-3}$, faster than the decay $r^{-2}$ of the disturbance due to self-propulsion.
Ignoring the motor torque simplifies the model, making it readily amenable to analytical treatment while still
capturing the key features of the experimental observations in the dilute limit.

Note that our model is in fact a model of a \emph{self-propelled} swimmer as opposed to
a body propelled by an external force, such as gravity or the magnetic field.
The propulsion of our swimmer is due to the point force
that models the effective action of flagellum on the fluid.
This force is balanced -- see equations (\ref{eq:force_balance})-(\ref{eq:torque_balance}) ) -- by an equal and opposite force of the fluid onto the flagellum
transmitted to the body of the swimmer.
Therefore the propulsion force is not external
since all external forces (e.g., gravity) are unbalanced.

Another issue of concern at very small scales is Brownian motion.
How reasonable is it to ignore such motion?
The rotational diffusion coefficient $D_\text{rot}$ for an ellipsoid of length $l$ and
diameter $d$ (see \cite{AraTsi05,DoiEdw86}) is
\begin{equation*}
    D_\text{rot} = \frac{12}{\pi}\ \frac{k_{_B} T}{l^3 \ln\left( l/d\right)\eta}\ ,
\end{equation*}
where $k_{_B}$ is the Boltzmann constant
\[
k_{_B} = 1.3806503 \times 10^{-23} \frac{\text{m}^2 \cdot \text{kg}}{\text{s}^{2}\cdot \text{K}},
\]
$T$ is temperature, and
$\eta$ is the viscosity of the fluid.

Computing the value of the diffusion coefficient $D_\text{rot}$ for
\emph{Bacillus subtilis} ($l\approx 4 - 5\,\mu$m and $d\approx 0.7 - 1.0\,\mu$m) in water at near-room temperature
($T\approx 300$ K and $\eta=0.8\times 10^{-3} \text{N}\cdot \text{s}\cdot \text{m}^{-2}$),
we obtain
\begin{equation*}
    D_\text{rot}
    \approx
    \frac{12}{\pi}\ \frac{1.3806503 \times 10^{-23} \cdot 300}
    {(5\cdot 10^{-6})^3 \cdot\ln(5)\cdot 8\cdot 10^{-3}} \
    \frac{\text{m}^2 \cdot \text{kg}\cdot \text{K}\cdot\text{m}^2}
    {\text{s}^{2}\cdot \text{K}\cdot \text{m}^3\cdot \text{N}\cdot \text{s}}
    \approx
    10^{-2} \text{s}^{-1}.
\end{equation*}

An isolated swimmer in the absence of Brownian motion will swim in a straight line.
In the presence of Brownian motion,
the expected time $T(\theta)$ to deviate by an angle $\theta$ from a given orientation
can be computed by using the appropriate first passage time
as
\begin{equation*}
    T(\theta) = \theta^2/D_\text{rot}.
\end{equation*}

%

Thus, for ``interior'' swimmers in the ``mirror image'' configuration (see section~\ref{subsect:Mirror Image})
the expected time to leave the basin of attraction of the ``swim in'' configuration
due to the described thermal effects becomes comparable to the interaction time (the ``swim in" time)
at distances of $100\,\mu$m and larger.
At smaller separations we can, therefore,
ignore thermal effects, at least at the qualitative level.

%


\subsection{PDE and ODE models for a swimmer in a fluid}

In the preceeding section we described the dynamics of a swimmer.
In this section we derive the PDE and ODE models for a collection of swimmers interacting with a fluid.
We consider several neutrally buoyant swimmers (indexed by a superscript $i=1,\dots,N$)
immersed in a Newtonian fluid (water) that occupies the domain $\Omega$.
We are concerned with instantaneous velocities of swimmers and fluid due to propulsion forces.
The head and the tail balls of the $i$th bacterium
are denoted by $B^i_{_H}$ and $B^i_{_T}$, respectively.
The corresponding coordinates, velocities, forces, and torques are labeled accordingly.
The swimmers occupy domain $\Omega_B = \bigcup_{i,*} B^i_*$ while the fluid occupies domain $\Omega_F = \Omega\setminus\Omega_B$.

Based on the typical swimming velocities and sizes of swimming bacteria, the Reynolds number of the fluid
flow induced by the motion of the swimmers is usually less than $10^{-2}$ (see, e.g., \cite{IshiSimPed06}).
Then inertia forces on the fluid elements are entirely dominated by viscous forces.
Ignoring the inertial effects of the fluid in our model is reflected by reducing Navier-Stokes equation to Stokes
equation, described below.

In the Stokesian framework, the Stokes drag law is applicable, stating that the viscous drag
on each of the balls is proportional to the radius $R$, while the mass of the ball is proportional to $R^3$.
Indeed, a neutrally-buoyant swimmer has the density of the surrounding fluid $\rho$, so its mass is $\frac{4}{3}R^3 \rho$.
For sufficiently small $R$ and finite density $\rho$ and viscosity $\mu$, the inertial terms $m \Bv_{_C}$ and $I \omega$ in (\ref{eq:swimmer_linear_dynamics}-\ref{eq:swimmer_angular_dynamics}) can therefore be neglected, so that
(\ref{eq:swimmer_linear_dynamics},\ref{eq:swimmer_angular_dynamics})
reduce to the balance of forces and torques:
$$
        0 = \BF_{_H} + \BF_{_T} + \BF_{_P},\qquad
        0 = \BT_{_H} + \BT_{_T}.
$$

%

As discussed above, at any time the fluid obeys the steady Stokes
equation on the domain $\Omega_F$ determined by the instantaneous configuration of the dumbbells.
The fluid is at rest at the outer boundary $\partial \Omega$ (the container) and is
coupled to the dumbbells only through the no-slip boundary
conditions on the surface of the swimmers.
Therefore, given the fixed magnitude $f_p$ of the force $\BF_{_P}=\tau f_p$ and the instantaneous positions
$\Bx^i_{_H}$ and $\Bx^i_{_T}$ of the balls $B^i_{_H}$ and $B^i_{_T}$, their instantaneous velocities
$\Bv^i_{_H}=\dot{\Bx}^i_{_H}$ and $\Bv^i_{_T}=\dot{\Bx}^i_{_T}$ are related to the fluid velocity $\Bu(\Bx)$ through the incompressible Stokes equation
\begin{equation}
    \label{eq:stokes}
    \left\{
    \begin{array}{l}
      \mu \lap \Bu = \D p + \rho \sum_i \delta(\Bx-\Bx^i_{_P}) \BF_{_P}^i \\
      \diiv(\Bu)=0
    \end{array}
    \right.
    \qquad \text{in }\ \Omega_F=\Omega\setminus\Omega_B,\quad \Omega_B = \bigcup_{i;\ *=H,T} B^i_*
\end{equation}
subject to the boundary and balance conditions
\begin{eqnarray}
    \Bu(\Bx) = 0, \qquad \Bx \in \partial \Omega, & \qquad & \text{\emph{container at rest}},  \label{eq:decay}\\
    \Bu = \Bv^i_{_C} + \omega^i \times \left(\Bx - \Bx^i_c) \right), \qquad \text{on }\d B^i_{_*}, &\qquad &\text{\emph{no-slip}, * = H,T}, \label{eq:noslip}\\
    \BF^i_{_{H}}+\BF^i_{_{T}} = -\BF^i_{_{P}},&\qquad &\text{\emph{balance of forces}},\label{eq:force_balance}\\
    \BT^i_{_{H}}+\BT^i_{_{T}} = 0, &\qquad &\text{\emph{balance of torques}}\label{eq:torque_balance}.
\end{eqnarray}
where $\mu$ is the viscosity of the fluid.
\def\full_system{(\ref{eq:stokes})-(\ref{eq:torque_balance})}
\def\stokes{(\ref{eq:stokes})-(\ref{eq:noslip})}
\def\balance{(\ref{eq:force_balance})-(\ref{eq:torque_balance})}

The system \full_system implicitly defines an ODE initial-value problem for the swimmers.
Indeed, given the instantaneous positions $\Bx^i_{_C}$ and orientations $\Btau^i$, the corresponding velocities $\Bv^i_{_C}$ and $\omega^i$ and the fluid velocity field $\Bu$ can be simultaneously determined from (\ref{eq:stokes})-(\ref{eq:torque_balance}).

The well-posedness of \full_system is established in
Appendix \ref{app:existence}.
Heuristically, the first equations (\ref{eq:stokes})-(\ref{eq:noslip}) can be
solved for $\Bu$ as a linear function of the velocities ${\Bv^i_{_C},\,\omega^i}$.
Doing so eliminates $\Bu$ from the remaining equations where it enters through the definitions
\eqref{eq:def:forcesNtorques} of forces $\BF^i_*$ and torques $\BT^i_*$.
Since \eqref{eq:def:forcesNtorques} is linear in $\Bu$, $\BF^i_*$ and $\BT^i_*$ also depend linearly
on $\Bv^i_{_C}$ and $\w^i$.
Hence, equations {\balance} provide a nonhomogeneous linear system for the unknown
velocities $\Bv^i_{_C}$ and $\w^i$ in terms of positions $\Bx^j_*$
and given intensity $f_p$ of the propulsion forces $\BF^i_{_P} = \tau^i f_p$.

Solving this system for the velocities $\Bv^i_*$ and $\w^i$
and using \eqref{eq:connection between forms of v}, we arrive at an ODE system for $\Bx^i_*$
\begin{equation}\label{eq:ODE system}
    \dot\Bx^i_* = \mathcal{V}(f_p,\Bx^j_{_H},\Bx^j_{_T},\Bx^j_{_P}),\qquad j=1,\dots,N.
\end{equation}

In the remainder of the paper we investigate the dynamics of swimmers, deriving appropriate approximate
ODEs (the fluid is acting only as a mediator of hydrodynamic interactions between the swimmers).

%

\section{Asymptotic reduction of the PDE model}\label{sect:Asymptotic reduction}

In the dilute limit the problem of determining the drag forces and velocities on the individual balls can be effectively approximated
by using three classical relations:
the Stokes drag law and
the basic solutions for the flow due to a point force
and for the flow due to a moving sphere.
We assume that the bacteria are
sufficiently long ($2L\gg R$)
and far apart ($|\Bx^i_c-\Bx^j_c|\gg 2L,\ i\neq j$)
so that all the balls and the point of application of propulsion forces  are well separated ($|\zeta-1|,|\zeta+1|\sim 1$).

At a point $\Bx$
the flow due to an isolated  propulsion force from the rotation of flagellum $-\BF_p$ is
given by $\Bu(\Bx) = -G\cdot\BF_p$.
Here $G$ is the Oseen tensor; see \eqref{eq:sol:point force} in Appendix~A.
The velocity of the fluid due to a ball moving with translational velocity
$\Bv$ in an unbounded fluid domain is given by $\Bu(\Bx) = H\cdot\Bv$.
Also, the drag force onto a ball moving with velocity $\Bv$ is given by Stokes formula
\begin{equation}\label{eq:Stokes formula for drag}
    \BF = \gamma_0 \Bv, \qquad \gamma_0 = 6 \pi R \mu
\end{equation}
(see Appendix \ref{sect:Stokes} 
for the definition of $G$ and $H$).

Furthermore, at distances large compared to the radius of the ball $R$ we have 
\begin{equation}\label{eq:H through G}
    H(\Bx) \approx \gamma_0 G(\Bx)\qquad \text{for }|\Bx|\gg R.
\end{equation}
%
Using \eqref{eq:H through G}, both the flows due to a point force $-\delta(\Bx^i_p) \BF^i_{_P}$
and due to translating spheres $B^i_*$, ($*=H,T$) we can write in the same form
\eqref{eq:flow due to a force}
in terms of the Oseen tensor $G$ and the forces $\BF^i_*$ ($*=H,T,P$) exerted by the fluid
\begin{equation}\label{eq:flow due to a force}
    \Bu^i_*(\Bx) = -G(\Bx-\Bx^i_*)\BF^i_*,\qquad *=H,T,\text{ or } P.
\end{equation}

We use \eqref{eq:Stokes formula for drag}
to relate the drag force $\BF^i_*$ (here $*=H,T$) on the ball $B^i_*$
($-\BF^i_*$ is the force onto the fluid)
to the velocity of the ball $B_*^i$ relative to the flow $\bar \Bu^i_*$,
given by \eqref{eq:def:external field},
\begin{equation}\label{eq:approx:force}
    -\BF^i_* = \gamma_0 \left(\Bv^i_* - \bar\Bu^i_*\right), \qquad *=H,T.
\end{equation}
Hence,
\begin{equation}
    \Bv^i_* = \bar\Bu^i_*+\frac{1}{\gamma_0}\BF_*^i,\qquad *=H,T.
    \label{eq:approx:velocity}
\end{equation}
For a given configuration of balls and point forces, remove one ball
$B^i_*$ and replace it by fluid.
Then the velocity of the the center of the ``fluid ball'' $B^i_*$
is denoted by $\bar \Bu^i_*$.
For instance, for the ball $B^i_H$ (i.e., $*=H$) $\bar \Bu^i_{_H}$ is given
as
\begin{eqnarray}\label{eq:def:external field}
    \bar\Bu^i_{_H} &=& \sum_{j\neq i}^N \bigg(\Bu^j_{_H}(\Bx^i_{_H})+\Bu^j_{_T}(\Bx^i_{_H})+\Bu^j_{_P}(\Bx^i_{_H})\bigg)+
    \Bu^i_{_T}(\Bx^i_{_H})+\Bu^i_{_P}(\Bx^i_{_H})
    = \\ \nonumber
    & = & -\sum_{j\neq i} \bigg( G(\Bx^i_{_H}-\Bx^j_{_H})\BF^j_{_H}+
    G(\Bx^i_{_H}-\Bx^j_{_T})\BF^j_{_T}+
    G(\Bx^i_{_H}-\Bx^j_{_P})\BF^j_{_P}\bigg)-\\
    \nonumber
    & &-
    G(\Bx^i_{_H}-\Bx^i_{_T})\BF^i_{_T}-
    G(\Bx^i_{_H}-\Bx^i_{_P})\BF^i_{_P}.
\end{eqnarray}
(Similarly we express $\bar \Bu^i_{_T}$.)

Note that $\bar\Bu^i_{_H}$ is defined only in terms of the forces $\BF^i_*$ and positions $\Bx^i_*$.
Hence, if we know all the forces $\BF^i_*$, equation \eqref{eq:approx:velocity}
together with \eqref{eq:def:external field} will give us all the velocities $\Bv^i_*$.
In the remaining part of this section we explain how to obtain a closed system of $6N$
equations for the unknown forces $\BF^i_*$ ($*=H,T$).
Since the number of unknown components of forces $\BF^i_{_H}$ and $\BF^i_{_T}$
is $2\times 3\times N=6N$, the system is, indeed, closed.

The relation \eqref{eq:force_balance} consists of $3N$ equations.
Since 
$R\ll L$, we have $(\Bx-\Bx_{_C})\approx(\Bx-\Bx_{_C})=\pm L\tau$.
Hence the torques $\BT^i_*$ on the balls $B^i_*$ can be approximated by the moments of the forces on them,
so that the torque balance relation \eqref{eq:torque_balance} becomes
\begin{equation}\label{eq:new torque zero}
    \tau^i\times(\BF^i_{_{H}}-\BF^i_{_{T}}) = 0.
\end{equation}
Since linear operator $A_\tau F:=\tau\times F$ has a 1D kernel, the
relation \eqref{eq:new torque zero} gives us another $2N$ equations.
To obtain the remaining $N$ equations, substitute the equations \eqref{eq:approx:velocity} into the rigidity of bacteria
equations \eqref{eq:rigidity},
and use expression \eqref{eq:def:external field}
for $\bar\Bu^i_{_H}$ and $\bar\Bu^i_{_T}$.
The obtained $6N$ equations schematically are denoted by
\begin{equation}\label{eq:schematic balance}
    \mathcal{L}_1\left(\Bx^i_*, \BF^i_{_H}, \BF^i_{_T}\right)=\mathcal{L}_2\left(\Bx^i_*, \BF^i_{_P}\right),
\end{equation}
where $\mathcal{L}_1$ and $\mathcal{L}_2$ depend on all positions $\Bx^i_*$ and
forces $\BF^i_*$ ($*=H,T,P,\ i=1,...,N$).
Here $\mathcal{L}_1$ and $\mathcal{L}_2$ are nonlinear in positions $\Bx^i_*$,
due to the terms of the form $G(\Bx^i_*-\Bx^j_*)$, and linear in forces $\BF^i_*$.

We solve \eqref{eq:schematic balance} for $\BF^i_{_H}$ and $\BF^i_{_T}$
in terms of the known $\BF^i_{_P}$ and $\Bx^i_*$.
Then $\Bv^i_{_H}$ and $\Bv^i_{_H}$ are defined by
(\ref{eq:approx:velocity}) and (\ref{eq:def:external field}),
\begin{equation}\label{eq:v defined}
    \Bv^i_* = \mathcal{V}(\BF^j_{_P},\Bx^i_{_H},\Bx^i_{_T},\Bx^i_{_P}),
\end{equation}
and we obtain an ODE system ($\dot\Bx^i_*=\Bv^i_*$):
\begin{equation}\label{eq:ODE 1}
    \dot\Bx^i_* = \mathcal{V}(\BF^j_{_P},\Bx^i_{_H},\Bx^i_{_T},\Bx^i_{_P})
    \qquad i,j=1,\dots,N.
\end{equation}

In particular, for a single swimmer
($N=1$) simple computations show
\begin{equation}\label{eq:sol:single bacteria}
    \Bv_{_H}=\Bv_{_T}= v_0 \tau,\qquad
    v_0 = \frac{f_p}{8\pi\mu L}\left[\frac{1}{2}+\frac{4 L}{3 R}+\frac{1}{|\zeta-1|}+\frac{1}{|\zeta+1|}\right],
\end{equation}
where $\tau$ is the unit vector along the axis of the bacterium.

\begin{figure}
  $
  \begin{array}{ccc}
    \includegraphics[width=6.5cm]{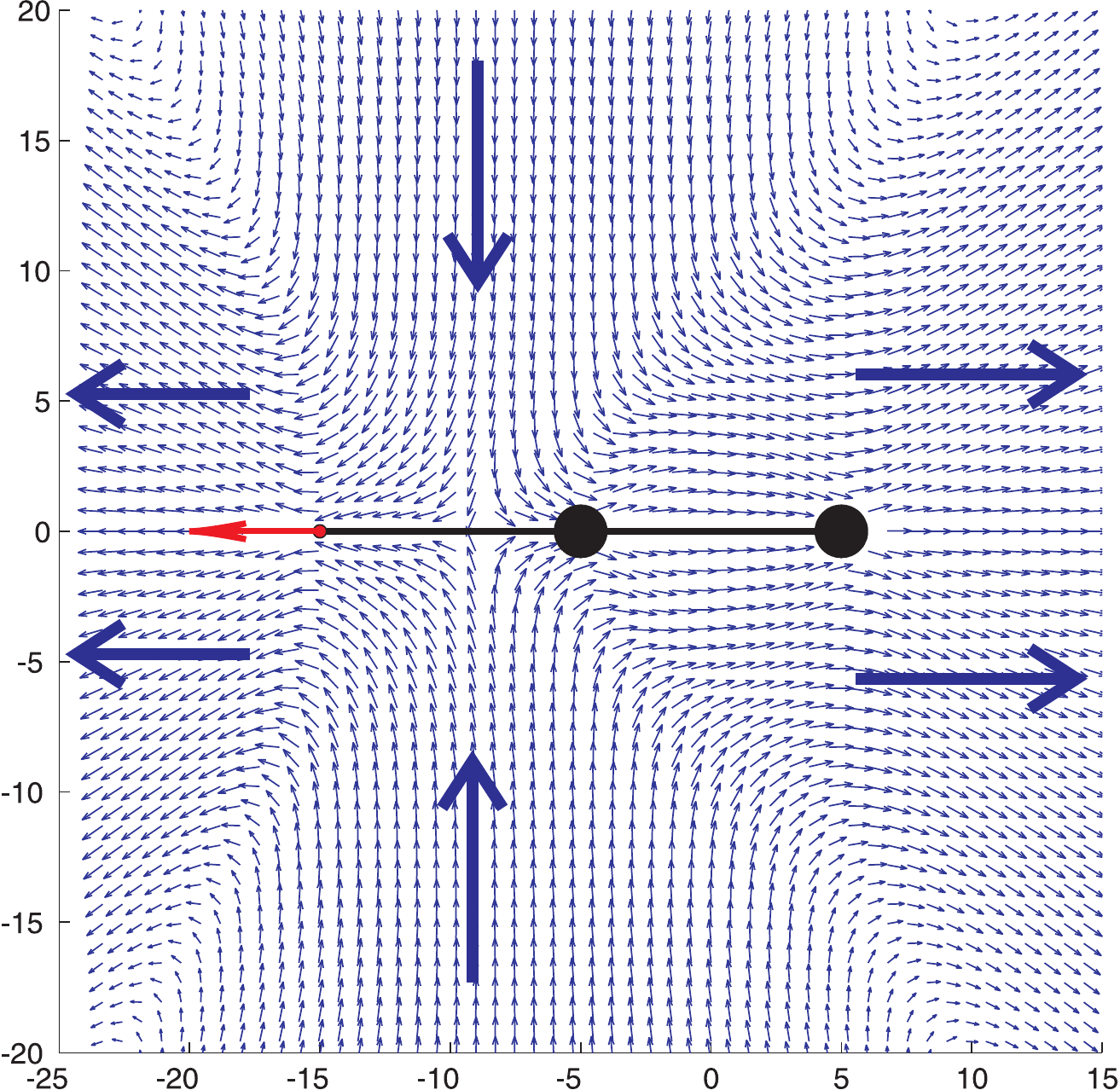}
    &
    \hspace{.2cm}
    &
    \includegraphics[width=6.5cm]{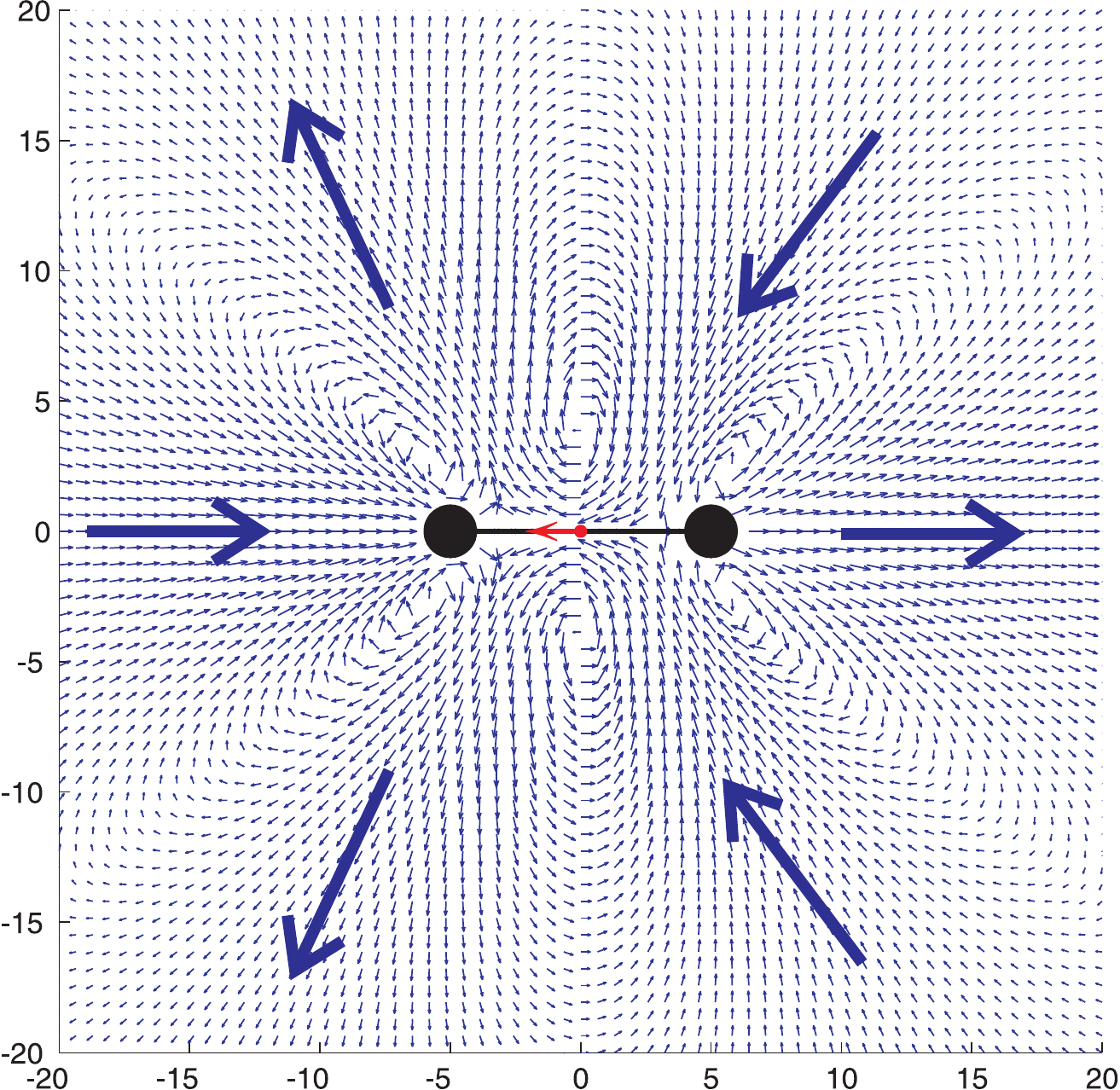}\\
    \text{(a)} && \text{(b)}
  \end{array}
  $
  \caption{
  (a) Velocity field of the fluid, computed from \eqref{eq:def:approx:u},
  around a swimmer with $\zeta=-3$.
  Heuristics: The bacterium is  moving left to right; hence the head ball pushes the fluid to the right.
  The force of the flagellum is pushing the fluid to the left.
  Because of incompressibility, the fluid is forced toward the bacterium from top and bottom.
  (b) Velocity field of the fluid, computed from \eqref{eq:def:approx:u},
  around a swimmer with $\zeta=0$.
  }
  \label{fig:bacteria field 1}
\end{figure}

Also, when all the forces have been found, one can find the velocity field of the fluid, using the
approximation \eqref{eq:H through G}:
\begin{eqnarray}
    \label{eq:def:approx:u}
    \Bu(\Bx) &=&  \sum_{i=1}^N \left(\Bu^i_{_H}(\Bx)+\Bu^i_{_T}(\Bx)+\Bu^i_{_P}(\Bx)\right)=\\
    &=&  -\sum_{i=1}^N \bigg( G(\Bx-\Bx^i_{_H})\cdot\BF^i_{_H}+
    G(\Bx-\Bx^i_{_T})\cdot\BF^i_{_T}+
    G(\Bx-\Bx^i_{_P})\cdot\BF^i_{_P}\bigg).
    \nonumber
\end{eqnarray}
This flow for a single bacterium is illustrated in Fig.~\ref{fig:bacteria field 1}.

\subsection{Bacterium as two force dipoles}

The above system of $6 N$ equations can be reduced to a smaller system for $N$ unknowns if one makes the following observation. From the balance equations \eqref{eq:force_balance},\eqref{eq:new torque zero}
and the form of the propulsion force $\BF^i_{_P} = f_p \tau^i$
it follows that
$\BF^i_{_H}$, $\BF^i_{_T}$, and $\BF^i_{_P}$ are all collinear.
Indeed, $\BF^i_{_P}$ is parallel to $\tau^i$ by definition,
and the balance equations \eqref{eq:force_balance},\eqref{eq:new torque zero}
imply that both the sum and the difference
of $\BF^i_{_H}$ and $\BF^i_{_T}$ are collinear with $\tau^i$.
Hence, $\BF^i_{_H}$ and $\BF^i_{_T}$ themselves are collinear with $\tau^i$.

\begin{figure}[h!]
  \begin{tabular}{l|r}
    \includegraphics[width=6cm]{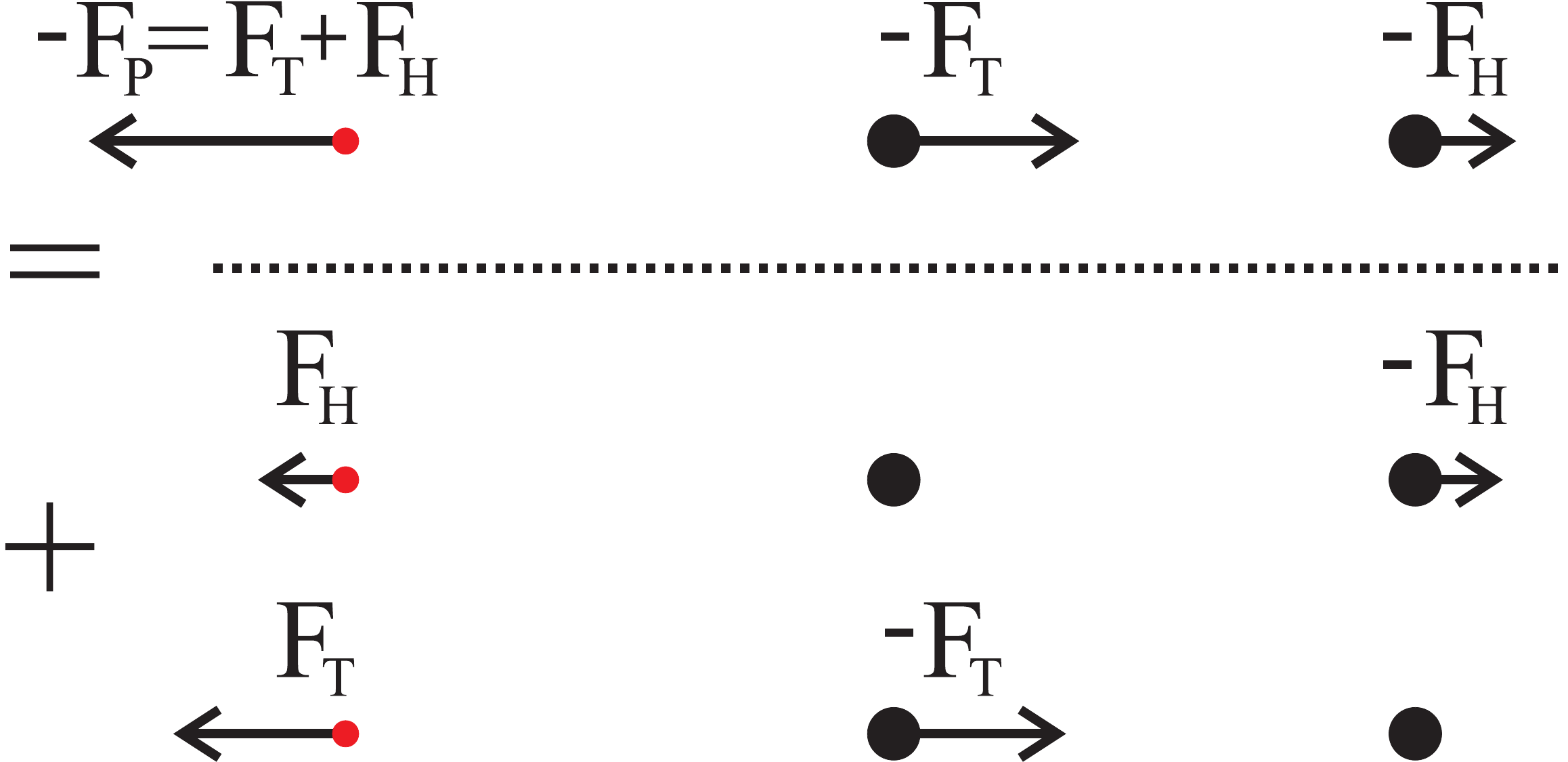} \hspace{1cm}
    &
    \hspace{1cm}\includegraphics[width=6cm]{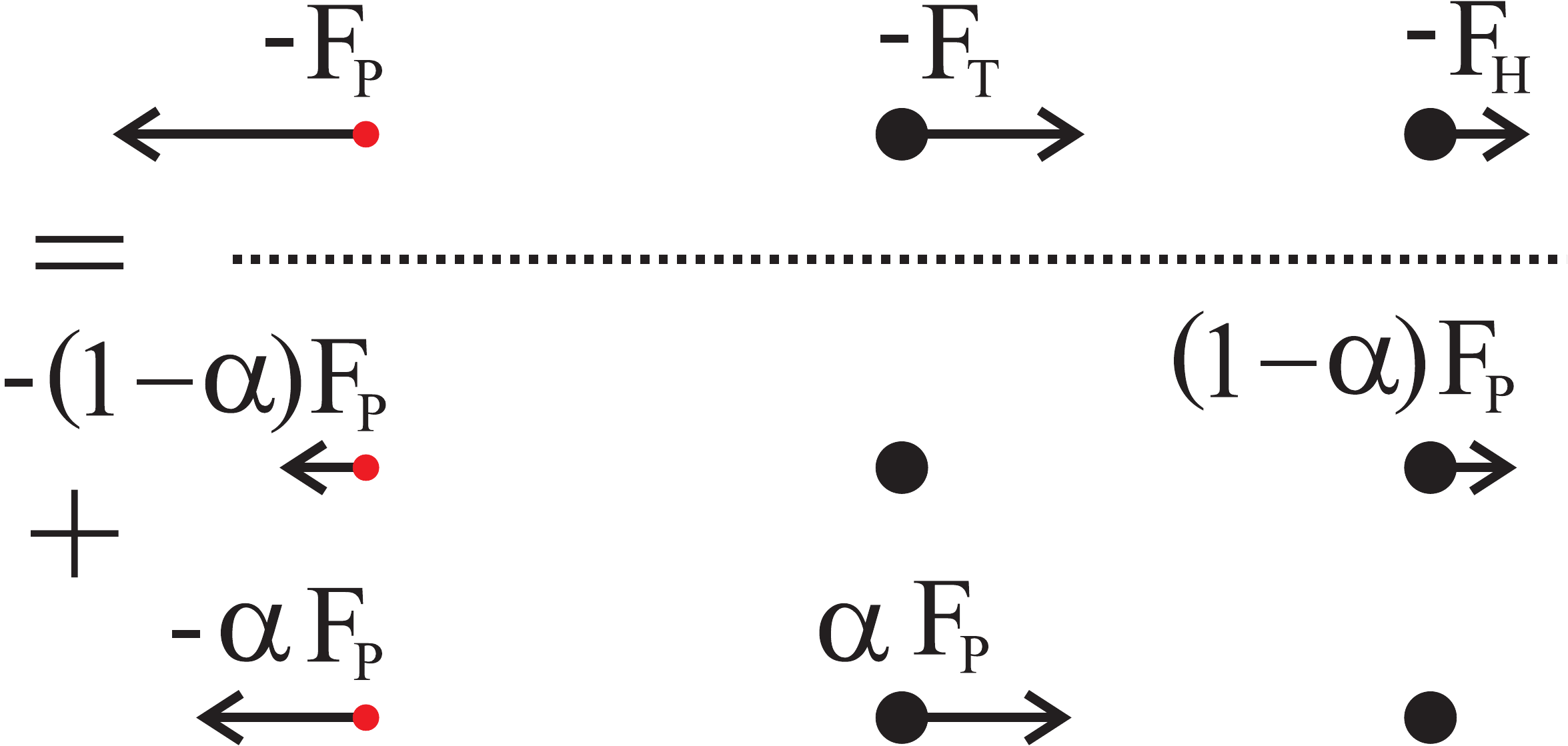} \\
  \end{tabular}
  \caption{Schematic representation of a swimmer as two ``force dipoles''
  in terms of $\BF_{_H}$ and $\BF_{_H}$ (left) and
  in terms of $\BF_{_P}$ and $\alpha$ (right).
  Arrows on the top line are the sum of arrows on the second and third line.
  }\label{fig:two dipoles}
\end{figure}

Therefore, for each bacterium there exists a scalar parameter $\alpha^i$ such that
\begin{equation}
      \BF^i_{_T}=\alpha^i \BF^i_{_P} = -\alpha^i f_p \tau^i,\quad
      \BF^i_{_H}=(1-\alpha^i) \BF^i_{_P} = -(1-\alpha^i) f_p \tau^i.
    \label{eq:def:alpha}
\end{equation}
Effectively, the parameter $\alpha^i$ groups the forces $\BF^i_*$ in two ``force dipole'' pairs
as illustrated in Fig.~\ref{fig:two dipoles}. Since $\BF^i_{_P}$ are given,
$\BF^i_{_H}$ and $\BF^i_{_T}$ are completely determined by the scalar $\a^i$.

We obtain a linear system for $\a^i$ by substituting \eqref{eq:def:external field}
into \eqref{eq:approx:velocity}
and using the rigidity constraint \eqref{eq:rigidity} for $\Bv^i_*$:
\begin{equation}\label{eq:rigidity system}
\begin{split}
    \a^i
    (\tau^i)^T\bigg[& G(\Bx^i_{_H}-\Bx^i_{_T})+G(\Bx^i_{_T}-\Bx^i_{_H})-\frac{2}{\gamma_0}\bigg]\tau^i+\\
    +
    \sum_{j\neq i} \a^j (\tau^i)^T\bigg[
    & -G(\Bx^i_{_H}-\Bx^j_{_H})+G(\Bx^i_{_T}-\Bx^j_{_H})
    +G(\Bx^i_{_H}-\Bx^j_{_T})-G(\Bx^i_{_T}-\Bx^j_{_T})
    \bigg]\tau^j=\\
    =
    \sum_{j\neq i}
    (\tau^i)^T\bigg[
    & -G(\Bx^i_{_H}-\Bx^j_{_H})+G(\Bx^i_{_T}-\Bx^j_{_H})
    +G(\Bx^i_{_H}-\Bx^j_{_P})-G(\Bx^i_{_T}-\Bx^j_{_P})
    \bigg]\tau^j+\\
    +
    (\tau^i)^T\bigg[
    & G(\Bx^i_{_T}-\Bx^i_{_H})+
    G(\Bx^i_{_H}-\Bx^i_{_P})-G(\Bx^i_{_T}-\Bx^i_{_P})
    -\frac{1}{\gamma_0}
    \bigg]\tau^i.
\end{split}
\end{equation}
Equation (\ref{eq:rigidity system}) all the coordinates $\Bx^i_*,\Bx^j_*$ and the directors $\tau^i,\tau^j$ are known;
$(\tau^i)^T$ stands for transpose of the vector $\tau^i$.
The system is linear for $\a^i$, and the coefficients depend
on the positions $\Bx^i_*,\Bx^j_*$ in a nonlinear way.


%

We now return to the ODE, equations \eqref{eq:ODE 1}.
As noted before, to find the explicit form for
$\mathcal{V}$, one has to find $\BF^i_{_H}$ and $\BF^i_{_T}$ in terms of $\BF^i_{_P}$
and $\Bx^i_*$.
Thus, the explicit form of $\mathcal{V}$ is determined by
(\ref{eq:def:alpha})-(\ref{eq:rigidity system}); see Appendix~\ref{sect:formulas}.

\subsection{Perturbation analysis of the ODE  for two bacteria}

\textbf{\underline{Notations}}

We next consider  \eqref{eq:ODE 1} for the particular case $N=2$.
Specifically, we focus on the dynamics of a coplanar pair of bacteria, when
$\tau^i,\ i=1,2$, lie in the $xy$-plane.
\begin{figure}[h!]
  \includegraphics[width=11cm]{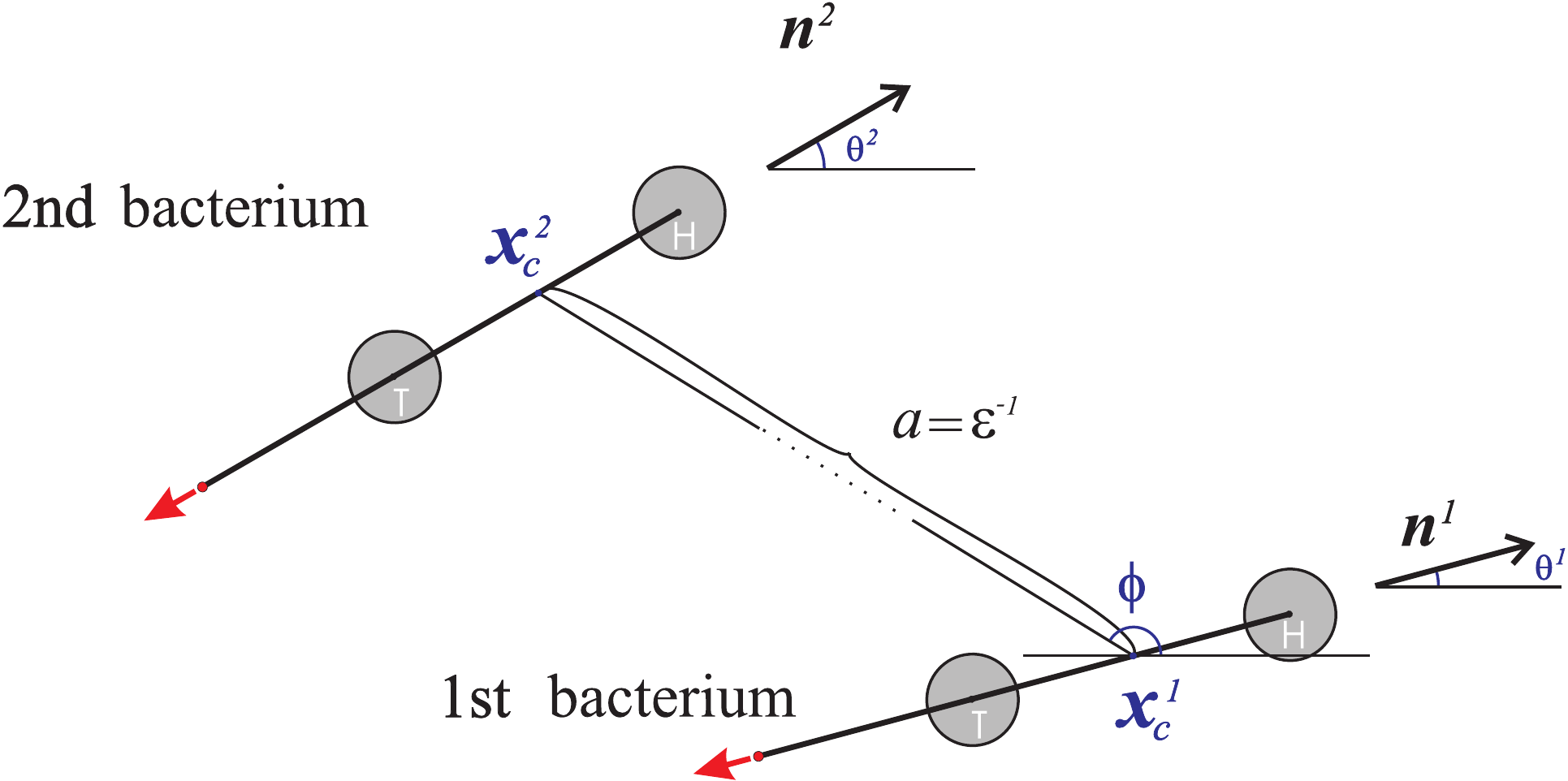}\\
  \caption{Coplanar pair of bacteria at a distance $a=\ve^{-1}$ apart.
  Here $\phi:=\angle\{(\Bx^2_c-\Bx^1_c),x\text{-axis}\}$ is the angle between $(\Bx^2_c-\Bx^1_c)$ and the $x\text{-axis}$;
  $\theta^1:=\angle\{\tau^1,x\text{-axis}\}$, $\theta^2:=\angle\{\tau^2,x\text{-axis}\}$.
  }
  \label{fig:two bacteria}
\end{figure}
To describe the relative positions of the two bacteria, we use the angles $\theta^1$, $\theta^2$, and $\phi$ shown in Fig.~\ref{fig:two bacteria}.
The angular velocity in this case is characterized by a single scalar
(no rotation around the axis $\tau^i$ of the $i$th bacterium) and is completely determined by
the velocities of the balls of the dumbbell:
\begin{equation}
    \omega^i = \dot \theta^i = \frac{\Bv^i_{_H}-\Bv^i_{_T}}{2L}\cdot \tau^i_{\perp},
\label{eq:def:w}
\end{equation}
where $\tau^i_{\perp}$ is obtained from $\tau^i$ by a $90^o$ in-plane rotation.
The velocity of the center of mass is
\begin{equation}
    \Bv^i_{_C} = \frac{\dot \Bx^i_{_H} + \dot \Bx^i_{_T}}{2} = \frac{\Bv^i_{_H}+\Bv^i_{_T}}{2},
\label{eq:def:vc}
\end{equation}
and, as mentioned before, the dynamics of the swimmer pair is completely determined by
$(\Bv^i_c,\,\omega^i)$.

\vspace{.5cm}

\textbf{\underline{Asymptotic expressions for velocities}}

As explained above, the ODE system \eqref{eq:ODE 1} can be written
explicitly ($\mathcal{V}$ can be expressed in terms of $\a^i,\Bx^i_*$, and $f_p$).
However, this system is too cumbersome for direct analysis.


But, since we confine ourselves to the dilute limit regime, we can satisfy ourselves
with asymptotic expressions in terms of the natural small parameter $\ve = a^{-1}$,
where $a$ is the distance between the centers of the two bacteria: $a=|\Bx^2_{_C}-\Bx^1_{_C}| \gg 1$.

We consider the asymptotic expansion \eqref{eq:anzats a} for $\a^u$.
After substituting it into equation \eqref{eq:rigidity system} for $\a^i$
and equating the terms at the same orders of $\ve=|\Bx^2_{_C}-\Bx^1_{_C}|^{-1}$
(see Appendix \ref{subsect:alpha expansion}), we obtain
\begin{equation}\label{eq:expansion alpha}
    \alpha^i = \a_0+O(\ve^2),\qquad \text{where }
    \a^i_0
    =
    \frac{1}{2}\
    \bigg(
    1+z(\zeta)\frac{R}{L}+z(\zeta)\left(\frac{R}{L}\right)^2+\dots
    \bigg),
\end{equation}
where $z(\zeta)$ is defined by \eqref{eq:z}.

To this end, we consider the asymptotic expansions
\begin{eqnarray}
    \label{eq:def:exp:vc}
    \Bv^i_c &=& \Bv^i_1  + \Bv^i_1 \ve + \Bv^i_2 \ve^2 + \Bv^i_3 \ve^3 + \dots,\\
    \label{eq:def:exp:w}
    \omega^i &=& \omega^i_0  + \omega^i_1 \ve + \omega^i_2 \ve^2 + \omega^i_3 \ve^3 + \dots
\end{eqnarray}
and substitute them into the LHS of (\ref{eq:def:w})-(\ref{eq:def:vc}).
In the RHS of (\ref{eq:def:w})-(\ref{eq:def:vc})
we express $\Bv^i_*$ using (\ref{eq:approx:velocity})-(\ref{eq:def:external field})
in terms of the
drag forces (or, equivalently, in terms of $\alpha^i$)
and expand the obtained formulas in $\ve$ (see Appendix \ref{sect:formulas}).

Equating the terms at every order of $\ve$ results in the following expressions for $\Bv^i_c,\,\omega^i$:
\begin{eqnarray}
  \label{eq:order 0}
  O(\ve^0)\ :     &\qquad & \omega^i_0=0, \qquad \Bv^i_0=v_0 \tau^i;\\
  \label{eq:order 1}
  O(\ve^1)\ :   &\qquad & \omega^i_1=0, \qquad \Bv^i_1=\mathbf{0};\\
  \label{eq:order 2}
  O(\ve^2): &\qquad & \omega^i_2=0, \qquad \Bv^i_2=A^i(f_p,L,R,\mu,\zeta,\alpha_0)\mathbf{B}^i(\theta^1,\theta^2, \phi);\\
  \label{eq:order 3}
  O(\ve^3): &\qquad & \omega^i_3=A^i(f_p,L,R,\mu,\zeta,\alpha_0)C^i(\theta^1,\theta^2, \phi),\\
  \label{eq:order 4}
  O(\ve^4): &\qquad & \omega^i_4=D^i(f_p,L,R,\mu,\zeta,\alpha_0) E^i(\theta^1,\theta^2, \phi),
\end{eqnarray}
where $j\neq i$.

At the leading order ($\ve^0$) $i$th bacterium swims straight
along its axis $\tau^i$ with a constant velocity $v_0$ -- as if there were no other bacterium.
Also, each bacterium can be viewed as two ``force dipoles''; see Fig.~\ref{fig:two dipoles}.
The disturbance due to a point force, given by $G(\Bx)\tau$, decays as $|\Bx|^{-1}$.
Hence, the disturbance due to a ``force dipole,'' given by
$\left(G(\Bx+\tau L)-G(\Bx-\tau L)\right)\tau,$
which is like a derivative of $G(\Bx)$, decays as $|\Bx|^{-2}$.
Therefore, the first nonzero correction in \eqref{eq:def:exp:vc} is $\Bv^i_2\ve^2$.
The first nonzero correction to rotational velocity $\w^i$ appears only at order
$\ve^3$.
Heuristically, this can be seen from \eqref{eq:def:w}, since the RHS of \eqref{eq:def:w}
is like a finite-difference derivative of the vector field, decaying as $\ve^2=|\Bx|^{-2}$.


Notice that all corrections starting from $\ve^2$ are in a separable form
\begin{equation}\label{eq:separable form}
    \text{Mat}(f_p,L,R,\mu,\zeta,\alpha_0)\text{Trig}(\theta^1,\theta^2, \phi).
\end{equation}
Here the function $\text{Mat}(f_p,L,R,\mu,\zeta,\alpha_0)$ is determined by the properties
($f_p$, $L$, $R$, $\mu$, $\zeta$, $\alpha_0$) and $\alpha_0=\alpha(L,R,\zeta)$, given by
(\ref{eq:rs o0 5})-(\ref{eq:z}), of bacteria and viscosity $\mu$ of the fluid
(material properties); the function $\text{Trig}(\theta^1,\theta^2, \phi)$
depends only on the mutual orientations ($\theta^1,\theta^2,\phi$) of bacteria.

The separable form \eqref{eq:separable form} allows us to study separately two questions:
\emph{(a)} For given material properties, how does the dynamics depend on the initial orientations of bacteria?
Here we show that swim in or swim off is determined by the sign of $\text{Trig}(\theta^1,\theta^2, \phi)$.
\emph{(b)} For given orientations, how does the dynamics depends on bacterial structure
(primarily, the position $\zeta$ of the propulsion force)?
Here we show that swim in or swim off is determined by the sign of $\text{Mat}(f_p,L,R,\mu,\zeta,\alpha_0)$.
In particular,
\begin{equation}\label{eq:properties of A}
\begin{split}
    A^i(f_p,L,R,\mu,\zeta,\alpha_0)>0\qquad &\text{for pushers }(\zeta<0),\\
    A^i(f_p,L,R,\mu,\zeta,\alpha_0)<0\qquad &\text{for pullers }(\zeta>0)
\end{split}
\end{equation}
(see Appendix \ref{subsect:A sign}).

For bacterium 1 we have
\begin{eqnarray}
    \label{eq:A=Algebraic}
    A^1(f_p,L,R,\mu,\zeta,\alpha_0)&=& \frac{f_p L}{32 \pi\mu} (1-\zeta-2\alpha_0)>0
    ,\\
    \label{eq:B=unit vectors}
    \mathbf{B}^1(\theta^1,\theta^2, \phi)&=&
    -2\big(1+3\cos(2(\theta_2-\phi))\big)\left[
    \begin{array}{c}
      \cos(\phi) \\
      \sin(\phi)
    \end{array}
    \right],
\end{eqnarray}
and
\begin{eqnarray}
    \label{eq:C=Trigonometric}
    C^1(\theta^1,\theta^2, \phi)
    &=&
    3\sin(\theta_1-\phi)  \big[
    5\cos(\theta_1+2\theta_2-3\phi)+2\cos(\theta_1-\phi)+\\
    \nonumber
    & &+
    \cos(\theta_1-2\theta_2+\phi)
    \big],\\
    D^1(f_p,L,\mu,\zeta)
    &=&
    \frac{3 f_p L^2(\zeta^2-1)}{256 \pi\mu},\\
    \label{eq:E}
    E^1(\theta^1,\theta^2, \phi)
    &=&
    35\sin(2\theta^1+3\theta^2-5\phi)
    +5\sin(2\theta^1+\theta^2-3\phi)+\\
    \nonumber
    & & +
    5\sin(2\theta^1-\theta^2-\phi)
    -4\sin(\theta^2-\phi)-\\
    \nonumber&-&
    20\sin(3\theta^2-3\phi)
    +3\sin(2\theta^1-3\theta^2+\phi).
\end{eqnarray}
To obtain the corresponding expressions
for bacterium 2, we simply switch the indexes 1,2 and replace $\phi$ with $(\pi+\phi)$, since $\phi$ is the angle
between $(\Bx^2_c-\Bx^1_c)$ and the $x$-axis.
Note that from \eqref{eq:A=Algebraic}, the sign of $C^i$ and $\mathbf{B}^i$, given by
(\ref{eq:C=Trigonometric},\ref{eq:B=unit vectors}), will give the sign of the first-order corrections to $\Bv^i_c$ and
$\omega^i$.

\subsection{Dynamics of two bacteria}

\begin{figure}[h!]
\[
\begin{array}{ccc}
  \includegraphics[width=6.5cm]{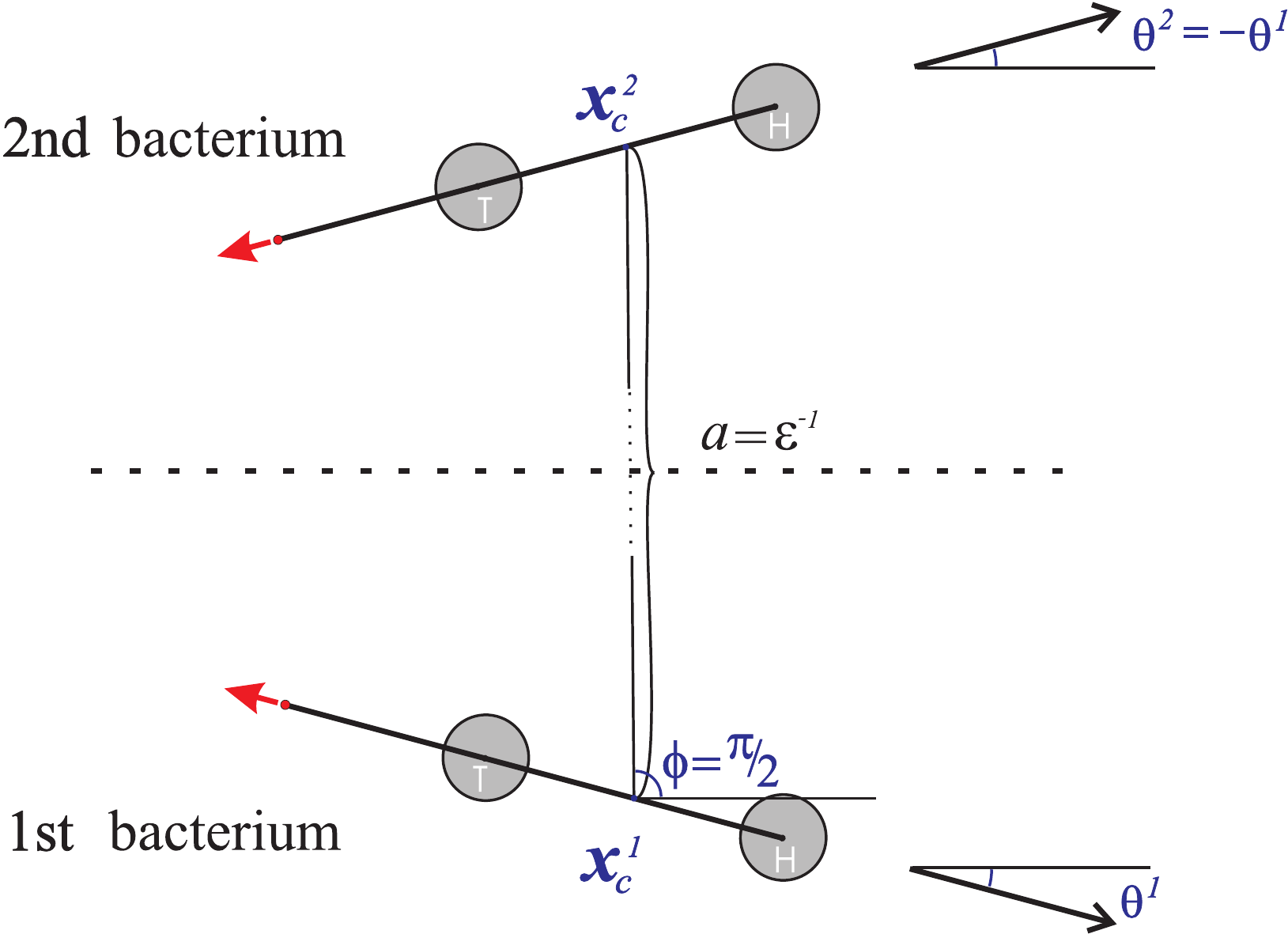} & &
  \includegraphics[width=6.5cm]{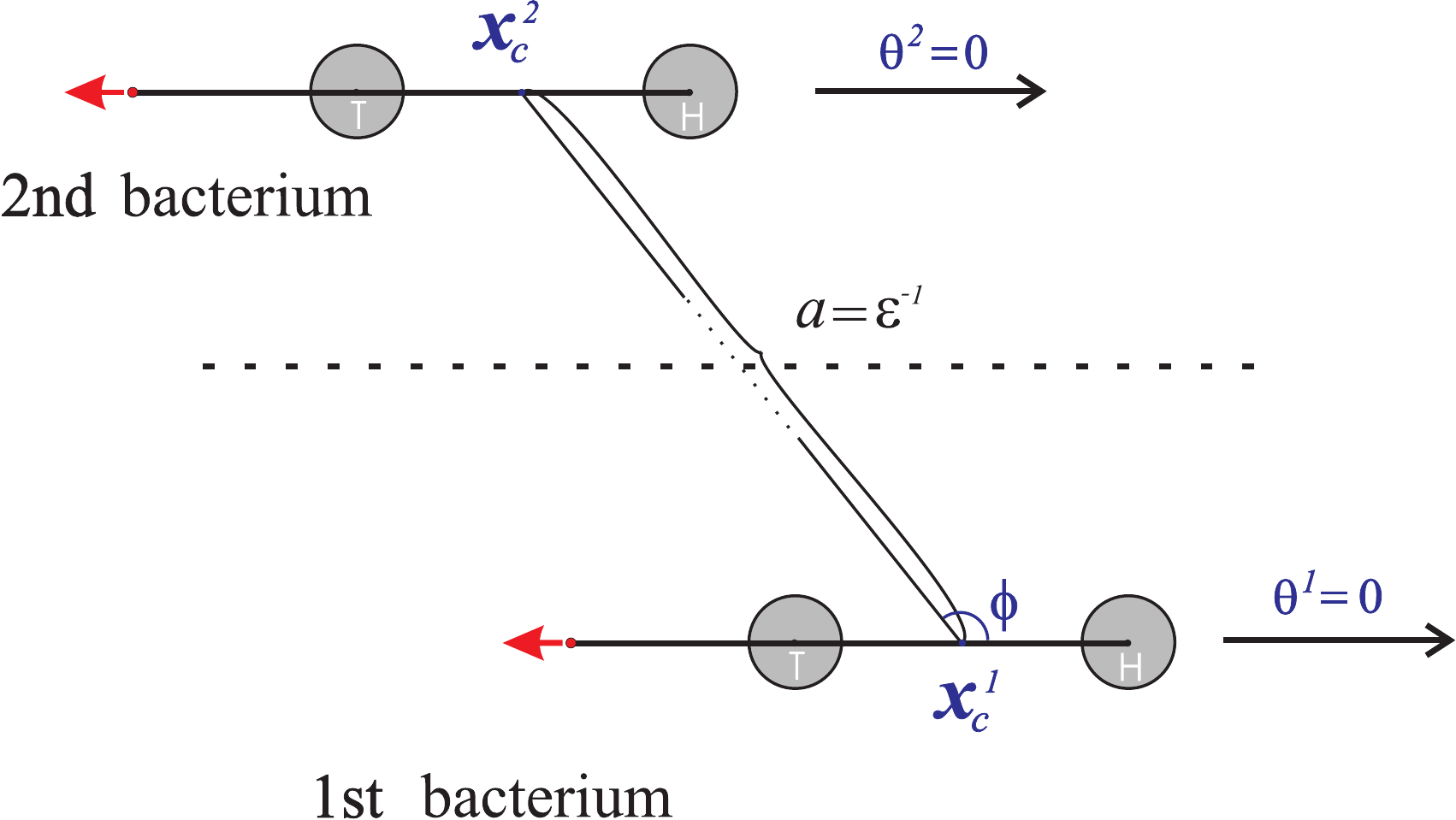} \\
  \text{(a) ``Mirror image'' configuration} & \qquad & \text{(b) ``Parallel'' configuration}
\end{array}
\]
\caption{Two basic configurations for the relative position of two bacteria.
Configuration 1 is called the ``mirror image'' configuration, because bacteria are symmetric relative to $x$-axis.
Configuration 2 is called the ``parallel'' configuration, because bacteria are parallel to one another.}\label{fig:two cases}
\end{figure}

The asymptotic formulas (\ref{eq:order 0})-(\ref{eq:order 4}) describe the dynamics
of a well-separated pair of bacteria.
The difficulty with interpreting these equations is the number of independent parameters
($\theta^1$, $\theta^2$, $\phi$), which does not allow having a single, comprehensive graph for the
trajectories of two bacteria.
Therefore, we consider two basic, yet representative, configurations (see Fig.~\ref{fig:two cases})
where there is only one free parameter and the remaining parameters are fixed.

The motivation for the choice of these basic configurations is twofold.
First, the evolution of a simple symmetric states, such as ``mirror image,'' provides insights
into the  behavior of the pair of bacteria in the course of collisions.
Second, these configurations allow for at least qualitative comparison with experimental data
(on swim in/swim off of a pair of bacteria as shown in \cite{AraSokGolKes07}).

\subsubsection{``Mirror image''}\label{subsect:Mirror Image}

We first consider the case with the two bacteria positioned symmetrically with respect to the $x$-axis
(see Fig.~\ref{fig:two cases}(a)).
Because of the symmetry, the positions of the bacteria will remain symmetric relative to the $x$-axis at all times.
%
Then the factors (\ref{eq:B=unit vectors},\ref{eq:C=Trigonometric},\ref{eq:E})
in the asymptotic expressions (\ref{eq:order 2})-(\ref{eq:order 4}) become
\begin{eqnarray}\label{eq:case1:C}
    C^1(\theta^1,\theta^2, \phi)
    &=&
    C^1(\theta^1,-\theta^1, \pi/2)
    =
    -3\sin(2\theta^1) [3- \cos(2\theta^1)],\\
    \nonumber
    C^2(\theta^1,\theta^2, \phi)
    &=&
    C^1(-\theta^1,\theta^1, -\pi/2)
    =
    3\sin(2\theta^1) [3- \cos(2\theta^1)]
    =
    -C^1(\theta^1,-\theta^1, \pi/2),
\end{eqnarray}
\begin{eqnarray}
    \label{eq:case1:B}
    \mathbf{B}^1(\theta^1,\theta^2, \phi)
    &=&\mathbf{B}^1(\theta^1,-\theta^1, \pi/2)
    =\left[
    \begin{array}{c}
      0 \\
      \ \ \,
      -2\big(1 - 3 \cos(2\theta^1)\big)
    \end{array}
    \right], \\
    \nonumber
    \mathbf{B}^2(\theta^1,-\theta^1, \pi/2)
    &=&
    \mathbf{B}^1(\theta^1,-\theta^1, \pi/2)
\end{eqnarray}
and
\begin{eqnarray}
    \label{eq:case1:E1}
    E^1(\theta^1,-\theta^1, \pi/2)
    &=&
    \cos(\theta^1)
    \left[2-56 \cos(2 \theta^1)+6\cos(4 \theta^1)\right],\\
    E^2(\theta^1,-\theta^1, \pi/2)
    &=&
    -
    E^1(\theta^1,-\theta^1, \pi/2).
\end{eqnarray}

\vspace{.5cm}

\underline{\emph{Analysis of steady states}}

If a steady (invariant) configuration of two bacteria exists (determined by $\theta^1$), it has to be
rotationally steady,
\begin{equation}\label{eq:rotationally steady}
    \w^1(\theta_1)=0=\w^2(\theta_1),
\end{equation}
and translationally steady,
\begin{equation}\label{eq:translationally steady}
    \Bv^1(\theta^1)=\Bv^2(\theta^1).
\end{equation}
We show below that no value of $\theta^1$ satisfies both
\eqref{eq:rotationally steady} and \eqref{eq:translationally steady}
(while each of these conditions can be satisfied separately).

\begin{figure}[!h]
  $
  \begin{array}{cc}
     \includegraphics[width=6cm]{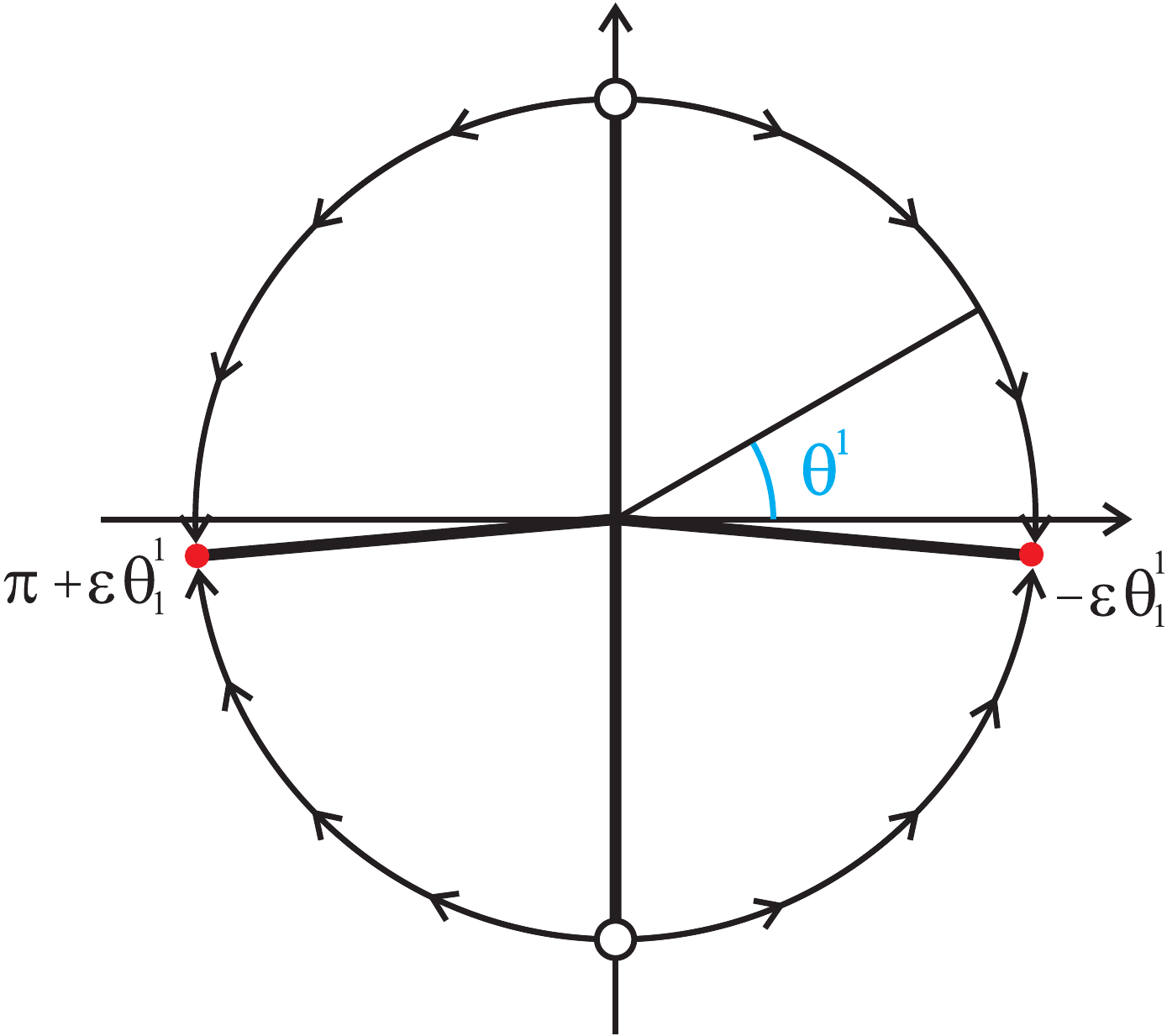}
    &
     \includegraphics[width=8cm]{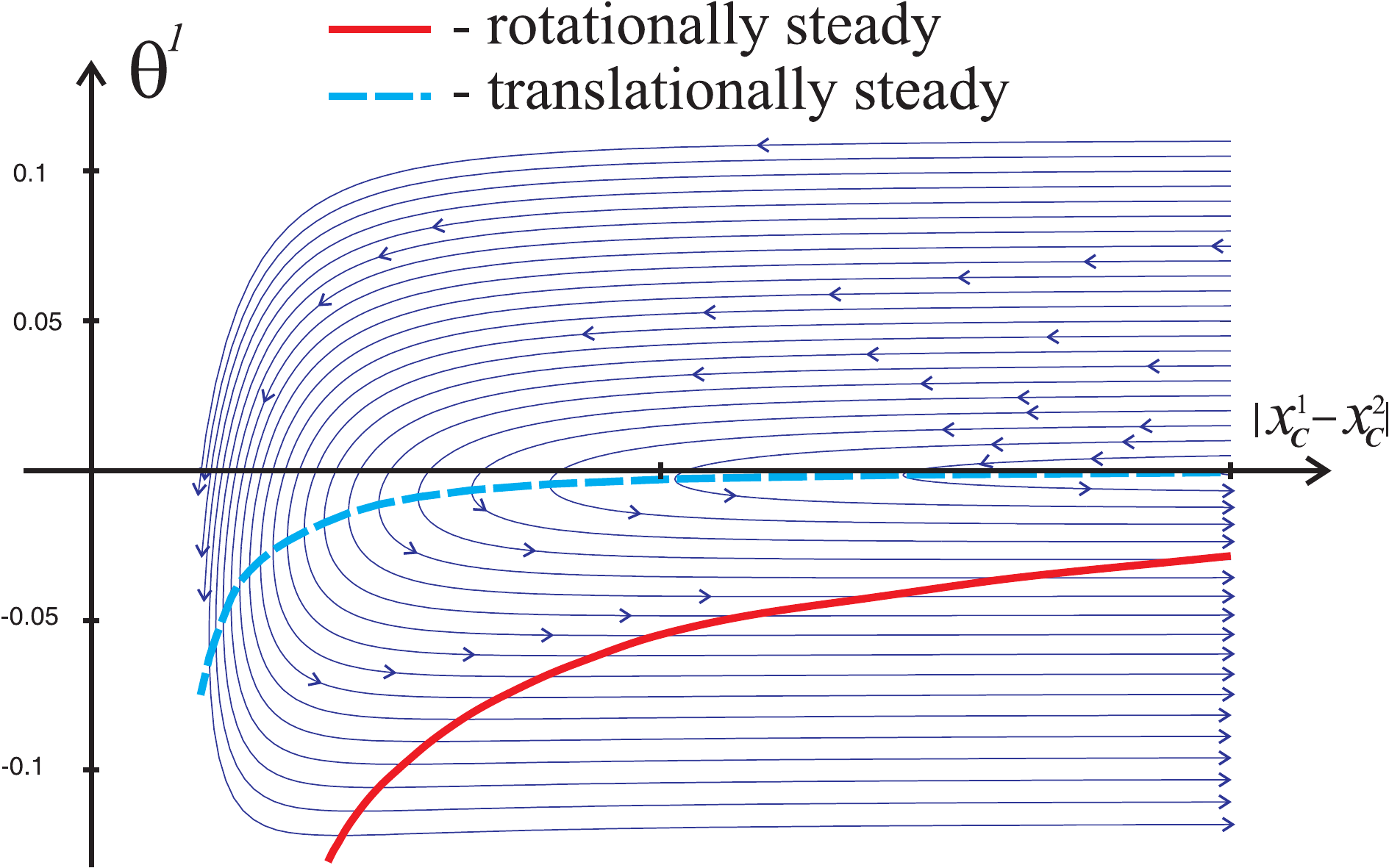}\\
    \text{(a)}
    &
    \text{(b)}
  \end{array}
  $
  \begin{small}
  \caption{
  Figures (a) and (b) correspond to the bacteria
  for which the position of the propulsion force is given by
  $\zeta=-2$.
  \emph{(a)} Schematic illustration of the dynamics of
  the angle $\theta^1$ (angle between the axis $\tau^1$ of the first bacterium and the $x$-axis)
  in the ``mirror image'' configuration.
  Arrows on the unit circle
  indicate the direction of change
  of the angle $\theta^1$.
  The states $\theta^1=\pm\pi/2$ (empty circles) are the unstable steady states for pushers ($\zeta<0$)
  and stable for pullers ($\zeta>0$), from \eqref{eq:properties of A}.
  The states $\theta^1=-\ve\theta^1_1$ and $\theta^1=\pi+\ve\theta^1_1$ (solid circles)
  are the stable steady states for pusher ($\zeta<0$) and
  unstable for pullers ($\zeta>0$), from \eqref{eq:properties of A}.
  \emph{(b)}
  The angle $\theta^1$ is plotted against the distance $|\Bx^1_{_C}-\Bx^2_{_C}|$ between bacteria.
  For $\theta^1>0$ (bacteria oriented inward) the dynamics indicated by arrows shows that
  bacteria move toward each other ($|\Bx^1_{_C}-\Bx^2_{_C}|$ decreases)
  and rotate outward ($\theta^1$ decreases).
  This action corresponds to the first part ($T_0<t<T_1$) of trajectories in
  Fig.~\ref{fig:case1:par traj}(a).
  The bold red curve indicates the rotationally steady states,
  obtained from \eqref{eq:rotationally steady}
  and (\ref{eq:order 3})-(\ref{eq:order 4}).
  The dashed blue curve indicates the translationally steady states, obtained from
  \eqref{eq:translationally steady} and (\ref{eq:order 0})-(\ref{eq:order 2}).
  These curves never intersect --
  no state is rotationally and translationally steady at the same time.
  }
  \end{small}
  \label{fig:phase plot}\label{fig:case1:rot stability}
\end{figure}

The trivial rotationally steady states are $\theta^1=\pm \pi/2$
(respectively, bacteria moving toward or away from one other on a vertical line).
These configurations are rotationally steady, since the configuration and the PDE
\eqref{eq:stokes}
are invariant under reflection across the $yz$-plane:
\begin{equation}\label{eq:symmetry 1}
    \left[
      \begin{array}{c}
        u^x(-x,y,z) \\
        u^y(-x,y,z) \\
        u^z(-x,y,z) \\
      \end{array}
    \right]
    =
    \left[
      \begin{array}{c}
        -u^x(x,y,z) \\
        u^y(x,y,z) \\
        u^z(x,y,z) \\
      \end{array}
    \right],\qquad p(-x,y,z)=p(x,y,z).
\end{equation}
Hence the trajectories and orientations of bacteria will also be invariant under this reflection.
Thus, bacteria starting with their centers $\Bx^i_{_C}$ on a vertical line ($x=z=0$) and oriented vertically
($\theta^1=\pm\pi/2$) will move vertically on that line.
Hence, these configurations are, indeed, rotationally steady.

But, since the distance between bacteria is not preserved,
these configurations are not translationally steady.
Thus, these configurations are not steady.


Two other rotationally steady angles (both are stable under variations of $\theta^1$ for pushers and unstable for pullers, due to \eqref{eq:properties of A}) are
\begin{equation}\label{eq:rotationally steady 2}
    \theta^1 = 0-\theta^1_1\ve+O(\ve^2),
    \text{ and }
    \theta^1 = \pi+\theta^1_1\ve+O(\ve^2),
\end{equation}
where
\begin{equation*}
    \theta^1_1 = 4\frac{A^1(f_p,L,\mu,\zeta,\alpha^2)}{D^1(f_p,L,\mu,\zeta)}=\frac{-3 L(\zeta+1) (1-\zeta)}{2(1-\zeta-2\alpha^2)}>0.
\end{equation*}

The angles \eqref{eq:rotationally steady 2} are found by setting
\begin{equation*}
    \w^1 = \ve^3 \w^1_3 + \ve^4 \w^1_4 + O(\ve^5)= O(\ve^5).
\end{equation*}

For these angles $\theta^1$, the vertical component of the translational velocity has the form
\begin{equation*}
    \Bv^1
    \left[
      \begin{array}{c}
        0 \\
        1 \\
      \end{array}
    \right]
    =
    |\Bv^1_0|\sin(\theta^1)+ O(\ve^2)
    =
    -|\Bv^1_0| \frac{2A^1(f_p,L,\mu,\zeta,\alpha^2)}{D^1(f_p,L,\mu,\zeta)}\ve+O(\ve^2)<0
\end{equation*}
for $\ve\ll 1$.
Thus bacteria are moving apart and the states are not translationally steady.

Therefore, there is no steady ``mirror image'' configuration of bacteria under the assumptions of the model.

\vspace{.5cm}

\noindent\underline{\emph{Dependence of the dynamics of bacteria on the position $\zeta$ of the propulsion force}}

Next, we plot the trajectories of two bacteria in the ``mirror image'' configuration.

\begin{figure}[h!]
  $
  \begin{array}{ccc}
    \includegraphics[width=4.5cm]{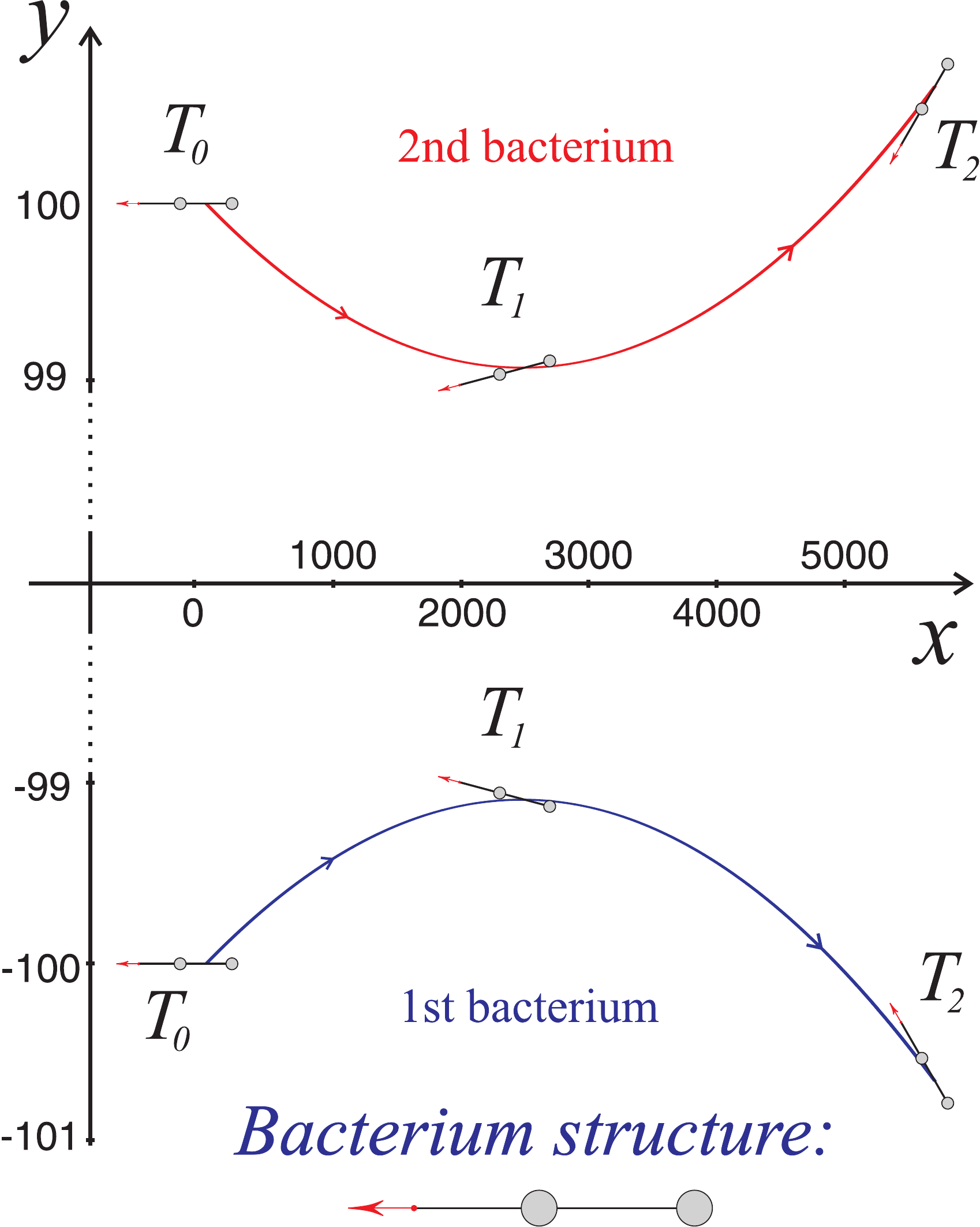}
    &
    \includegraphics[width=4.7cm]{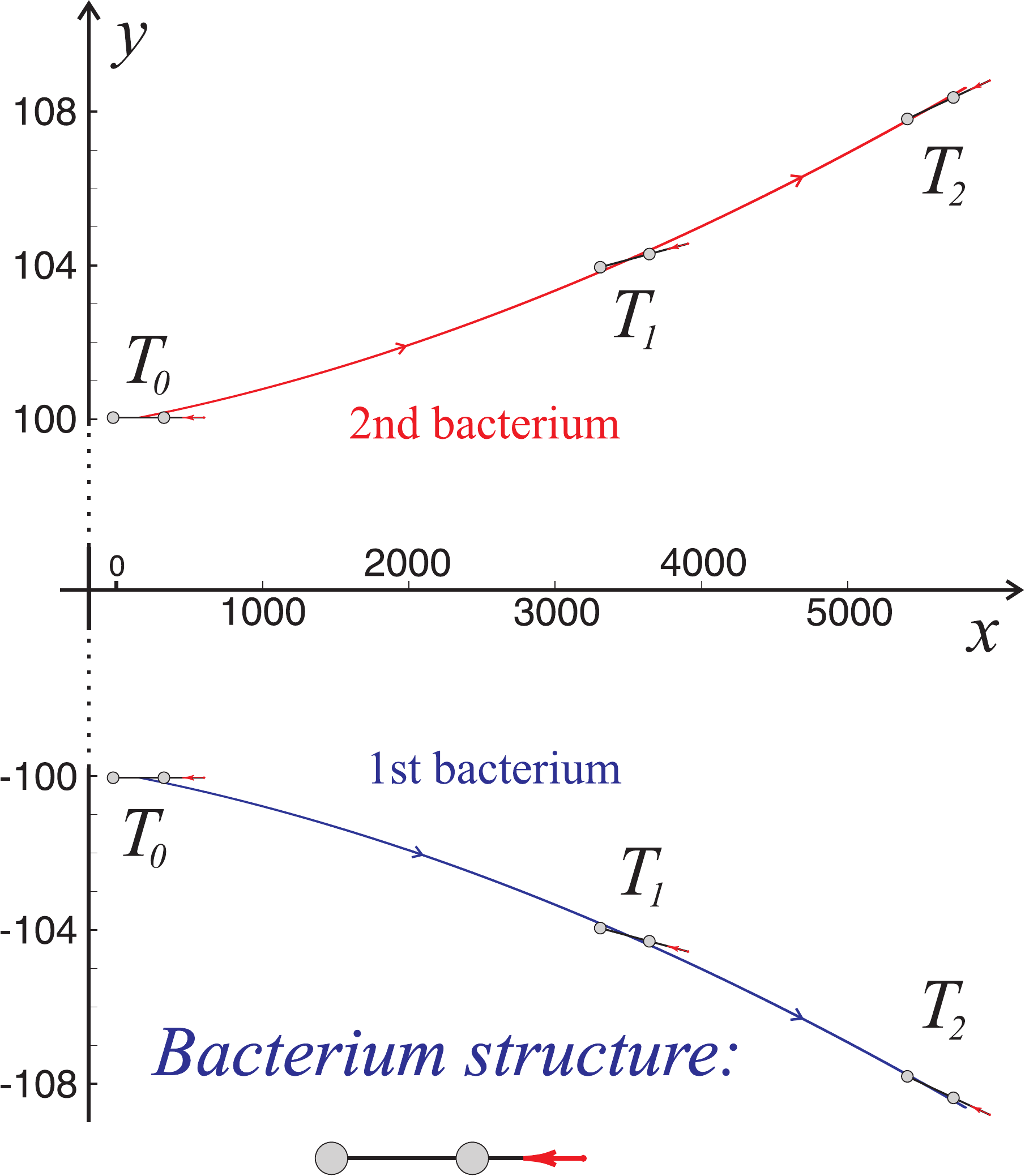}
    &
    \includegraphics[width=4.5cm]{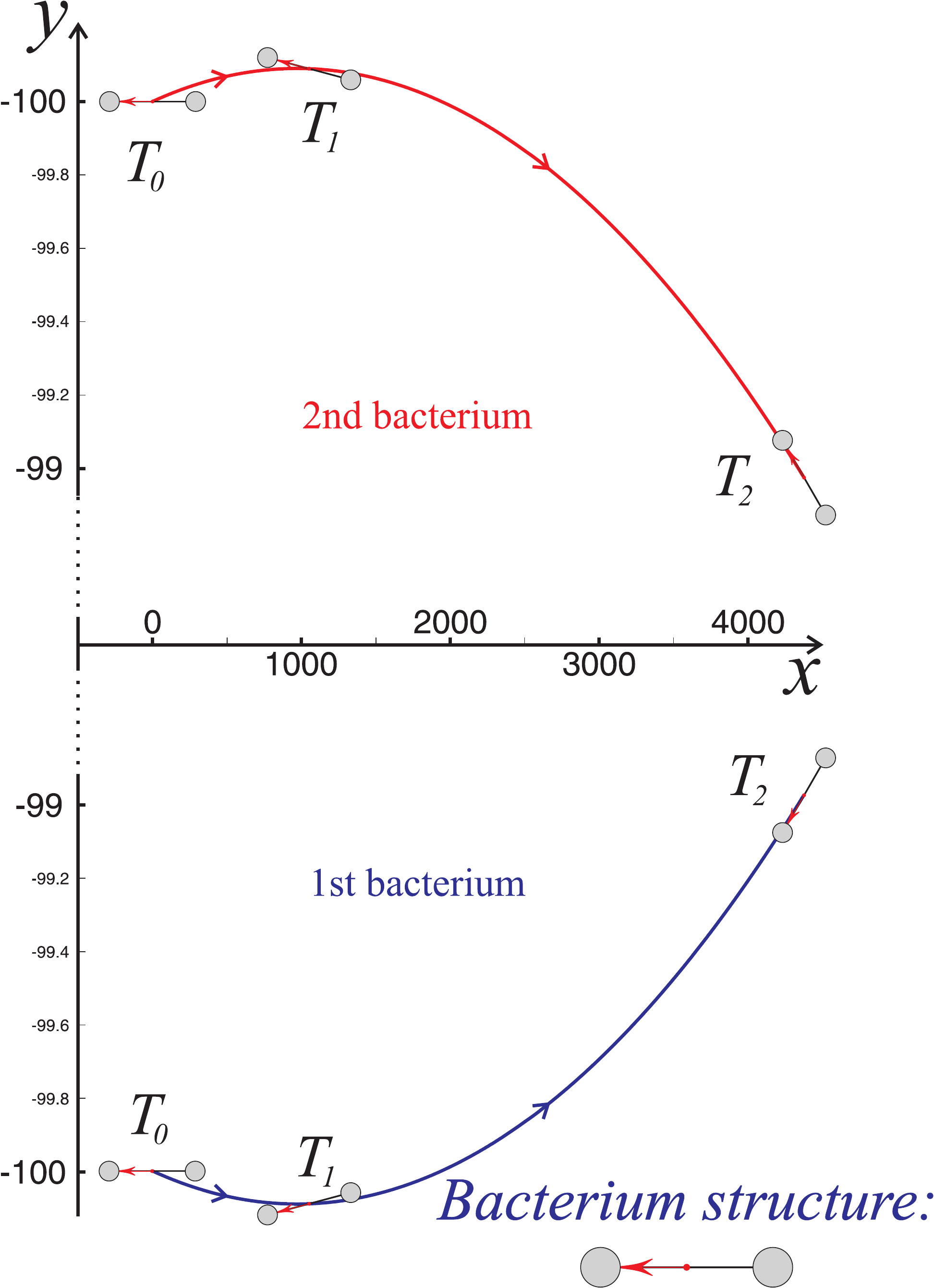}\\
    \text{(a)}
    &
    \text{(b)}
    &
    \text{(c)}
  \end{array}
  $
  \caption{Trajectories of initially parallel bacteria in the ``mirror image'' configuration
  (radius of dumbbell balls $R=1$).
  (a) External pushers ($\zeta<-1$) at first ($T_0<t<T_1$) attract and rotate outward.
  When rotated sufficiently outward ($t>T_1$), the bacteria swim off.
  (b) External pullers ($\zeta>1$) swim off.
  (c) Internal swimmers ($|\zeta|<1$) at first ($T_0<t<T_1$) repel and rotate inward.
  When rotated sufficiently inward ($t>T_1$), the bacteria swim in.}
  \label{fig:case1:par traj}
\end{figure}

We choose an initial orientation of bacteria parallel to the $x$-axis ($\theta^1=0$).
The trajectories of the centers of bacteria are shown in Fig.~\ref{fig:case1:par traj}.

We observe that when the propulsion force is applied between the dumbbell balls ($-1<\zeta<1$),
bacteria initially move apart and rotate inwards (see Fig.~\ref{fig:case1:par traj}.c).
After time $t_0$ (when bacteria have rotated sufficiently inwards) they start approaching each other
and swim in. Eventually, the distance between the bacteria decreases, and the assumptions about well-separated bacteria become
invalid, so more accurate representations of the drag forces and velocities are needed to address evolution of the pair
in this state. Remarkably, this behavior is consistent with the experimentally observed attraction  between two nearby
bacteria; see Fig.~1  in \cite{AraSokGolKes07}.
While experiments suggest the existence of long-living states of a close pair of bacteria swimming on parallel tracks,
it is likely that this state cannot be properly captured in the asymptotic far-field approximation for the velocity fields of moving spheres  used in our paper.

The detailed explanation of this behavior is as follows.
Initially, at the leading order ($\ve^0$) the translational motion is along the
$x$-axis ($\Bv^i_0=v_0 \tau \parallel ox$).
The next correction $\ve^2\Bv^i_2$ is directed outwards; hence, initially the bacteria move apart.
This can be seen by substituting
(\ref{eq:order 0})-(\ref{eq:order 2})
and (\ref{eq:A=Algebraic},\ref{eq:case1:B}) into the expansion \eqref{eq:def:exp:vc}.
As the bacteria move apart they are rotating inwards due to the $\ve^3\w^i_3$ term in \eqref{eq:def:exp:w}.
The rotation changes the orientation of $\tau^i$ and hence of the leading-order translational motion
$\Bv^i_0=v_0 \tau$.
At $t=t_0$, which solves
\begin{equation}\label{eq:t0}
    \left(\Bv_0^i(t_0)+\ve^2 \Bv_2^i(t_0)\right)\cdot e_2 = 0,
\end{equation}
the bacteria rotated sufficiently inwards that the terms
$\Bv^i_0$ and $\ve^2\Bv^i_2$ balance each other.
After this moment ($t>t_0$) the contribution of $\Bv^i_0$
to the motion along $oy$-axis dominates $\ve^2\Bv^i_2$.
Hence, bacteria start approaching each other (swim in).

We also observe that when the propulsion force is not between the dumbbell balls
($|\zeta|>1$), the bacteria swim off (see Fig.~\ref{fig:case1:par traj}(a) and \ref{fig:case1:par traj}(b)).
These observations emphasize the fact that the dynamics of the pair of bacteria depends sensitively on the
position of the propulsion force and, consequently, on the shape of the microorganisms and the structure of its propulsion.

\noindent\underline{\emph{Stability of the ``mirror image'' configuration:}}\\
The ``mirror image'' configuration is a reduction that allows us to describe the state of the swimmer pair
with only two parameters ($a$ and $\theta^1$).  How generic is this subset within the space of all configurations?
Appendix~\ref{subsect:stability of MI} address this question in some detail
and shows that the ``outward'' configuration ($\theta^1 < 0$) is stable whereas the ``inward'' configuration ($\theta^1 >
0$) is unstable, so that nearby configurations in the ``general position'' tend to approach the ``outward mirror image''
state but not the ``inward mirror image,'' at least when the interbacterial distance is large: $a \gg 1$.
Since the ``inward mirror image'' is central to our description of asymptotic scattering of swimmers, we briefly comment on its validity.
We regard this configuration as representative of the general asymptotic dynamics in that if
a ``swim off'' (see below) occurs for the interaction of swimmers in the ``inward mirror image'' position,
it will certainly occur in the ``general position'' case.
At the same time, if a ``swim in'' occurs in the ``inward mirror image'' situation
(as is shown below for specific choices of the force location),
it is likely to occur for the nearby ``general position'' configurations,
since the crucial ``inward'' character of the configuration is robust to perturbations.
We plan to investigate this matter more closely in the future
by considering a wider subspace of swimmer pair configurations.

\begin{figure}[!h]
\[
\begin{array}{ccc}
  \includegraphics[width=4.5cm]{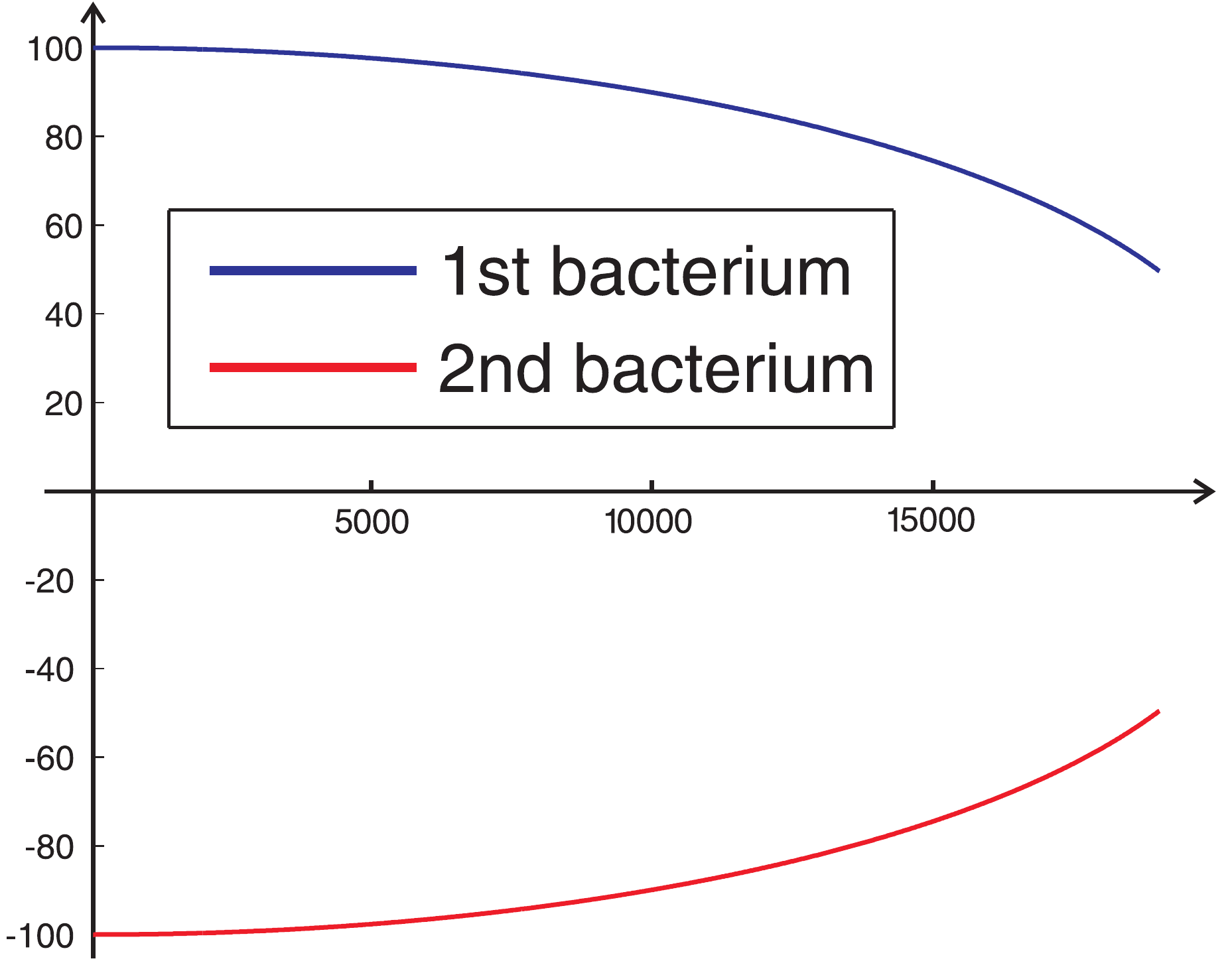}
  &
  \includegraphics[width=4.5cm]{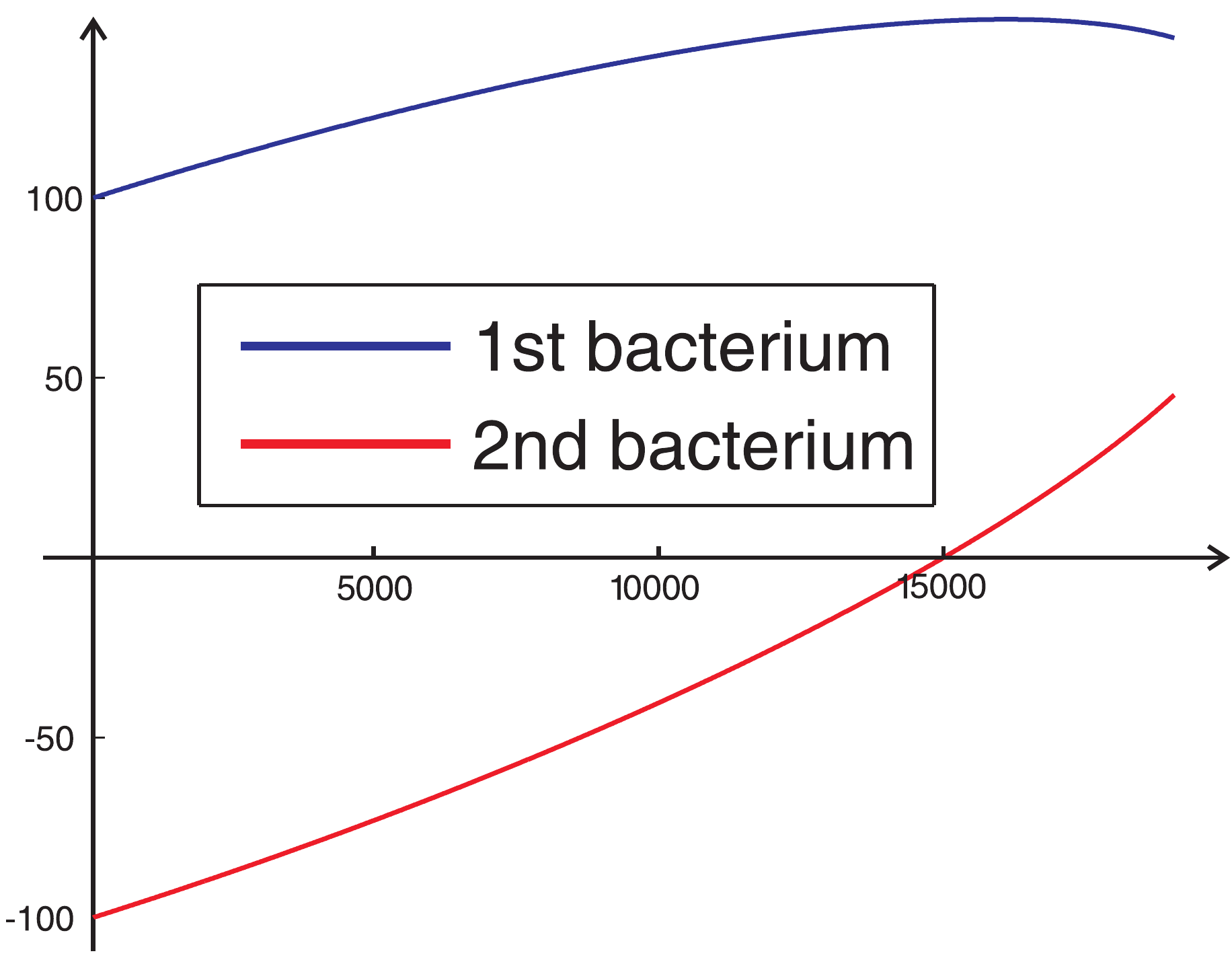}
  &
  \includegraphics[width=4.5cm]{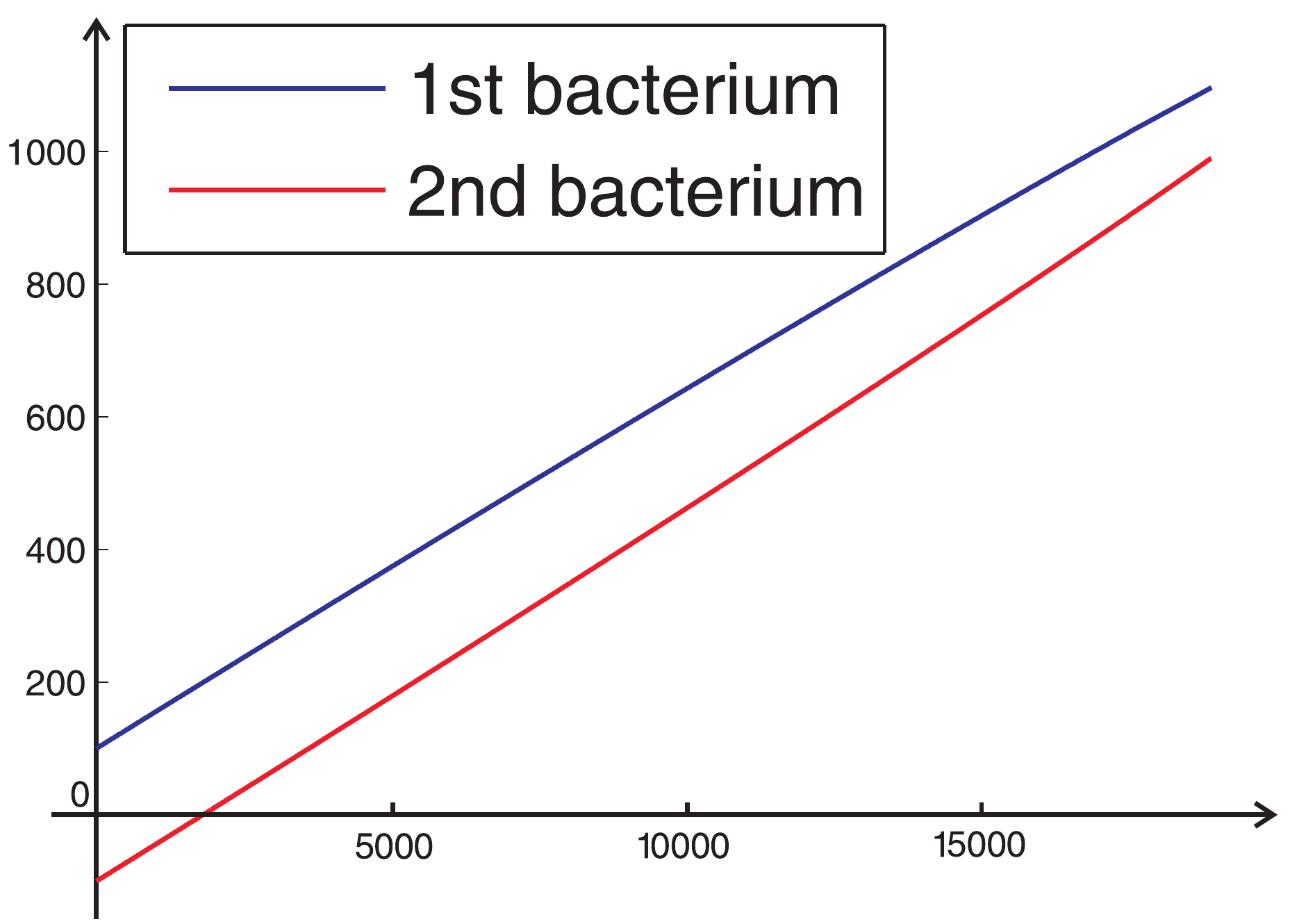}  \\
  (a) & (b) & (c)
\end{array}
\]
  \caption{Pictures \emph{(a)-(c)} show the trajectories for inner swimmers ($\zeta=0$)
  starting from a perturbed parallel ``mirror image'' configuration.
  The measure of perturbation $\delta$ is defined by \eqref{eq:perturbation of the mirror image}.
  The picture \emph{(a)} corresponds to $\delta(0)=0$, that is, the unperturbed ``mirror image.''
  The picture \emph{(b)} corresponds to $\delta(0)=0.01$, and
  the picture \emph{(c)} corresponds to $\delta(0)=0.1$.
  The unit of length here is the radius $R=1$ of a ball in the swimmer dumbbell.
  }
  \label{fig:delta perturbation}
\end{figure}

\begin{figure}[!h]
  \begin{center}
    \includegraphics[width=10cm]{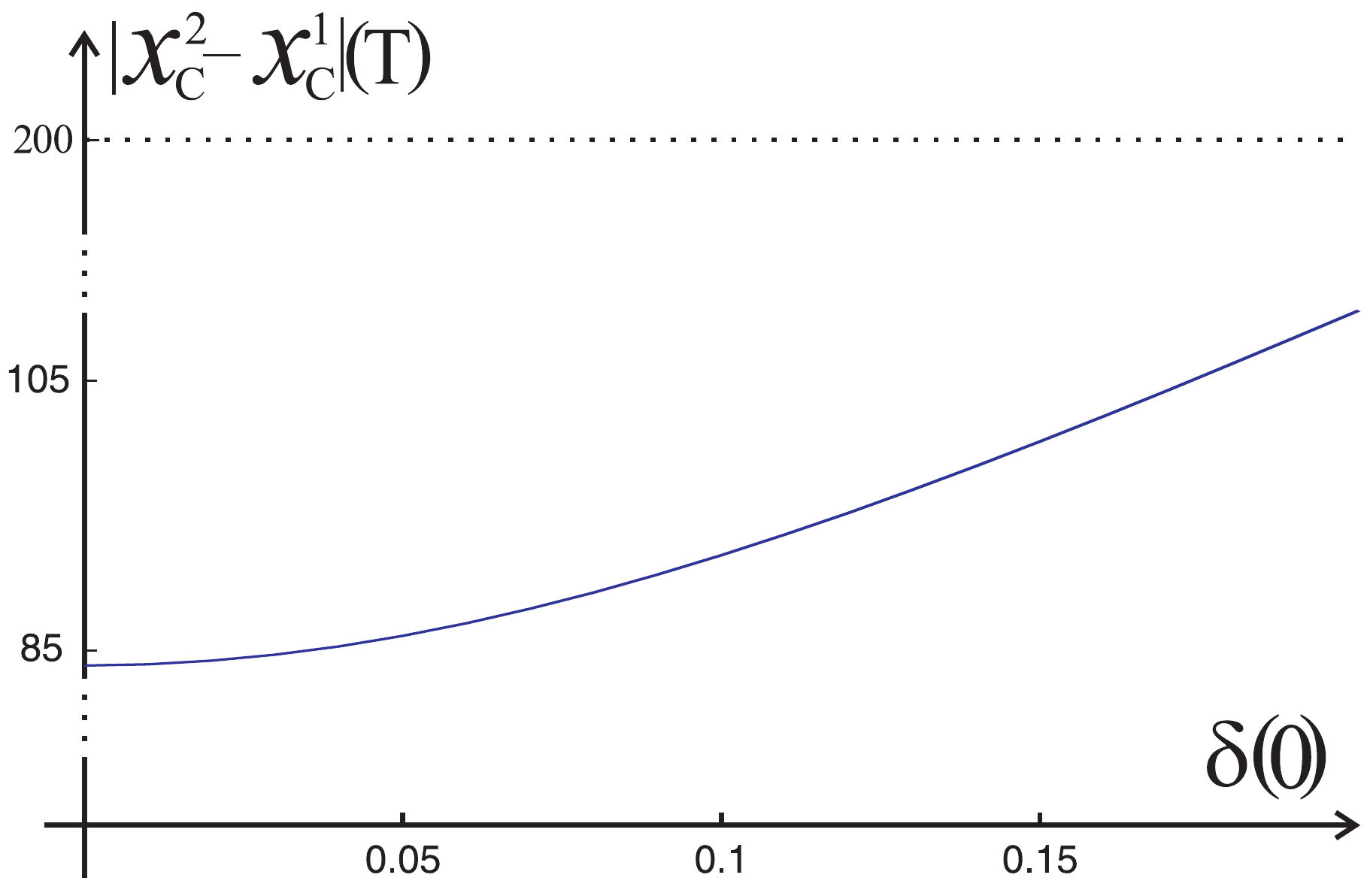}
  \end{center}
  \caption{Distance between two swimmers, starting from a perturbed parallel ``mirror image'' configuration distance 200 apart, at time $T$ ($T\approx 100$ seconds) as a function of the initial perturbation $\delta(0)$.
  The unit of length here is the radius $R=1$ of a ball in the swimmer dumbbell.}
  \label{fig:distance with delta}
\end{figure}

The ``mirror image'' configuration is stable under small perturbations in orientations of the swimmers when they are rotated outward from one another.
The ``mirror image'' configuration is unstable under small perturbations in orientations of the swimmers when they are rotated inward to one another (see Appendix \ref{subsect:stability of MI}).

Nevertheless, the swim in or swim off of swimmers in the ``mirror image'' configuration
(as can be seen from Fig.~\ref{fig:delta perturbation} and Fig.~\ref{fig:distance with delta}) is representative of their mutual dynamics, resulting from hydrodynamic interactions.

\subsubsection{``Parallel'' configuration}

Here we consider a pair of bacteria in a ``parallel'' configuration (see Fig.~\ref{fig:two cases}):
one located ahead of the other parallel to one another ($\theta^1=\theta^2$)
and initially parallel to the $x$-axis.
The parameter $\phi$ measures the angle between $(\Bx^2_{_C}-\Bx^1_{_C})$ and the $x$-axis.
As bacteria may not be aligned with the ${x}$-axis for $t>0$,
we introduce another parameter $\tphi$ --
the angle between $(\Bx^2_{_C}-\Bx^1_{_C})$ and the axis of the first bacterium $\tau^1(t)$.
Thus, $\tphi(t) = \phi(t)-\theta^1(t)$.


The factors (\ref{eq:order 2})-(\ref{eq:order 4})
in the asymptotic expressions
(\ref{eq:def:exp:vc})-(\ref{eq:def:exp:w}), when written in terms of $\tphi$ instead of $\phi$, become
\begin{eqnarray}
    \label{eq:case2:C1}
    C^1(\theta^1,\theta^2, \tphi)
    &=& C^1(0,0, {\tphi})
    = 3 \sin(2 {\tphi})\left[1 - 5 \cos(2 {\tphi})\right]   \\
    \label{eq:case2:C2}
    C^2(\theta^1,\theta^2, \tphi)
    &=& C^1(0,0, \tphi+\pi)
    = 3 \sin(2 {\tphi})\left[1 - 5 \cos(2 {\tphi})\right]
    = C^1(0,0, \tphi)
\end{eqnarray}
and
\begin{eqnarray}\label{eq:case2:B1}
    \mathbf{B}^1(\theta^1,\theta^2, \tphi)
    &=&
    \mathbf{B}^1(0,0, \tphi)
    =
    -2 \big(1 + 3 \cos(2{\tphi})\big)
    \left[
    \begin{array}{c}
      \cos({\tphi}) \\
      \sin({\tphi})
    \end{array}
    \right],\\
    \label{eq:case2:B2}
    \mathbf{B}^2(\theta^1,\theta^2, \tphi)
    &=&
    \mathbf{B}^1(0,0,\tphi+\pi)
    =
    2 \big(1 - 3 \cos(2{\tphi})\big)
    \left[
    \begin{array}{c}
      \cos({\tphi}) \\
      \sin({\tphi})
    \end{array}
    \right],
\end{eqnarray}
\begin{eqnarray}
    \label{eq:case2:E1}
    E^1(\theta^1,\theta^1, {\tphi})
    &=&
    -\frac{1}{2}\sin({\tphi})
    \left[9+20 \cos(2 {\tphi})+35\cos(4 {\tphi})\right],\\
    \label{eq:case2:E2}
    E^2(\theta^1,\theta^1, {\tphi})
    &=&
    -
    E^1(\theta^1,\theta^1, {\tphi}).
\end{eqnarray}

Since for a general angle $\tphi$ term $E^1(0,0, \tphi)\neq 0$,
equation \eqref{eq:case2:E2} implies that
\begin{equation}\label{eq:blahblah 1}
    \w^1({\phi})\neq \w^2({\phi}).
\end{equation}
This means that a pair of bacteria in the ``parallel'' configuration may not remain
in the ``parallel'' configuration at some later time.
In other words, the ``parallel'' configuration may not be preserved in time
(unlike the ``mirror image'' configuration, which is preserved in time).

The only angles $\tphi$ for which bacteria remain in the ``parallel'' configuration
are $\tphi=0$ and $\tphi=\pi$.
These angles correspond to a pair of bacteria one following another on the same straight line;
we call this the \emph{head-to-tail} configuration.
The difference between $\tphi=0$ and $\tphi=\pi$ is only in assigning numbers to bacteria
($\tphi=\pi$ means that the leading bacteria is called the first,
while
$\tphi=0$ means that the trailing bacteria is called the first).
Next, w.l.o.g. we consider the case $\tphi=0$.

\underline{``Head-to-tail'' configuration}\\
From the top-bottom symmetry, it follows that for $\tphi=0$ the rotational corrections at all orders vanish.
For instance, to the order $\ve^4$ this can be checked by plugging $\tphi=0$ into
(\ref{eq:order 3})-(\ref{eq:order 4}) using (\ref{eq:A=Algebraic}),(\ref{eq:C=Trigonometric}) and (\ref{eq:E})

The stability (under variations in $\tphi$) of the ``head-to-tail'' configuration of bacteria
is determined by the sign of the leading-order correction terms $C^1$ and $C^2$ in the
rotational velocities $\w^1$ and $\w^2$; see (\ref{eq:case2:C1},\ref{eq:case2:C2}).

Take $\theta^1(0)=\theta^2(0)$ (bacteria initially aligned with ${x}$-axis)
and $\tphi(0)=\phi(0)$ small positive (the second bacterium is ahead and slightly above the first one).
Then $C^1=C^2<0$, which means $\w^1,\w^2<0$: the  bacteria are rotating clockwise.
The angle  $\tphi=\phi-\theta^1$ increases.
Similarly, take $\tphi(0)=\phi(0)$ small negative.  Then $C^1=C^2>0$,
which means $\w^1,\w^2>0$:
the bacteria are rotating counterclockwise and $\tphi=\phi-\theta^1$ decreases.

Therefore, from \eqref{eq:properties of A},
for pushers ($\zeta<0$) the ``head-to-tail'' configuration is unstable and
for pullers ($\zeta>0$) it is stable. This result is in fact consistent with the simulations of \cite{HerStoGra05}
indicating formation of close ``head-to-tail'' pairs of puller dumbbells.

\vspace{.5cm}

\noindent\underline{\emph{Dependence of the dynamics of bacteria on the position $\zeta$ of the propulsion force}}

We study the dependence of the dynamics of two swimmers in the ``head-to-tail'' configuration,
depending on the position $\zeta$ of the propulsion force.
Since the two swimmers are positioned on the same line (there is no preferred direction
other than this line), will stay on this line and can either get closer together or
get farther apart as they move on this line.

\begin{figure}[h!]
    \begin{center}
    \includegraphics[width=10cm]{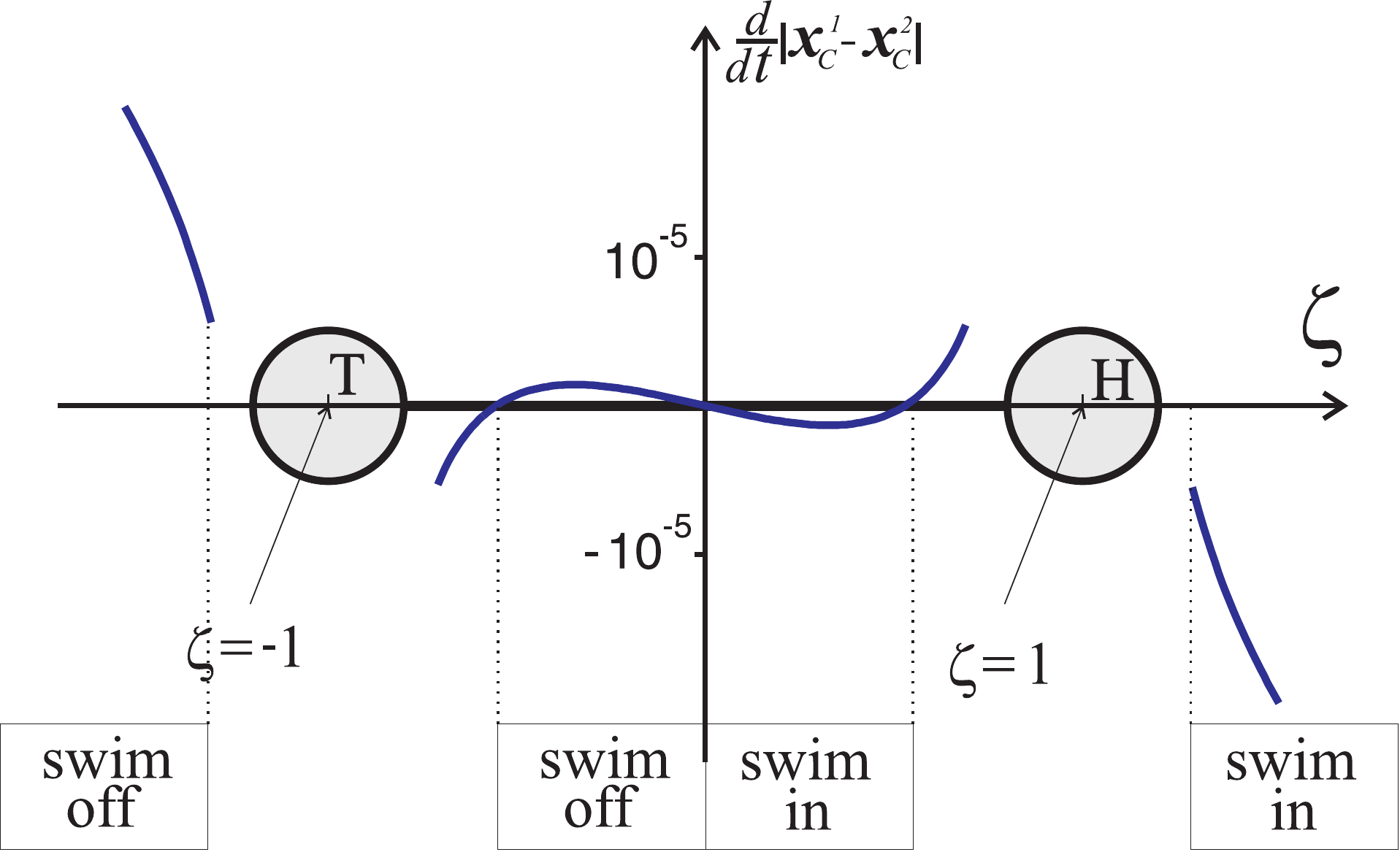}
    \end{center}
  \caption{
  Dependence of the relative velocities of two swimmers in the ``head-to-tail'' configuration
  on position $\zeta$ of the propulsion force.
  Pushers ($\zeta<0$) swim off; pullers ($\zeta>0$) swim in.
  The asymptotic technique used by us cannot be applied
  to the uncharacterized regions (between dashed lines) of $\zeta$ close to $\pm 1$.
  }
  \label{fig:case2:stag distance}
\end{figure}

We observe that pushers ($\zeta<0$) always swim off and pullers ($\zeta>0$)
always swim in (see Fig.~\ref{fig:case2:stag distance}).

\subsection{Quasi-two-dimensional model}

In this section we consider two bacteria swimming in a thin film (quasi-two-dimensional fluid, abbreviated Q2D).
The interest in studying this case is due to a number of physical experiments
(e.g., \cite{WuLib00,DomCisChaGolKes04,AraSokGolKes07,SokAraKesGol07}) observing the motion
of bacteria in a thin film (in particular, in experiments in \cite{SokAraKesGol07} the thickness of the film was of the
same order as the thickness of the bacteria).
The thin film allows us to focus a microscope on individual bacteria and track their motion with time.

The modeling in a thin film (of thickness $2h$) differs from the above model in the whole space
because the boundary conditions on the top and bottom of the thin film must be taken into account.
While the experiments in  \cite{WuLib00,SokAraKesGol07}
were performed with free-standing fluid film, suggesting free slip boundary conditions on the interfaces,
the experiment in \cite{SokAraKesGol07} indicates formation of thin, solid-like walls on the fluid-air interfaces
due to the byproducts of bacteria metabolism.
Therefore, in fact, the correct boundary conditions for the in-plane  velocities are no-slip.

Hence, instead of the fundamental solution $G(\cdot)$ of the Stokes equation in the whole space,
we use
its Q2D analog --  the Green's function
$\tilde{G}(\cdot)$ with no-slip boundary conditions on the horizontal walls ($z=\pm h$):
\begin{equation}\label{eq:BC tilde G}
    \tilde{G}(x,y,h)=\tilde{G}(x,y,-h)=\mathbf{0}.
\end{equation}
The series expansion for the velocity of the fluid due to a point force $\delta(\Br)e_1$
is obtained in \cite{LirMoc76a}:
\begin{equation}\label{eq:tilde u}
    \Bu(\Br)=\tilde{G}(\Br) e_1.
\end{equation}
Taking the leading term in this series (as $|\Br|\to\infty$),
we get an approximation
\begin{equation}\label{eq:Q2D:point force}
    \Bu(\Br)=
    \left[
      \begin{array}{c}
        u_x(\Br) \\
        u_y(\Br) \\
        u_z(\Br) \\
      \end{array}
    \right]
    \approx
    \frac{f(z)}{|\Br|^4}
    \left[
      \begin{array}{c}
        x^2-y^2 \\
        2xy \\
        0 \\
      \end{array}
    \right],
\end{equation}
where $f(z)$ is a known function (see \cite{LirMoc76a} and Appendix \ref{sect:Stokes:Q2D}).


Analogously to the 3D approximation \eqref{eq:H through G},
we want to approximate the fluid flow due to a sphere moving
(in the $xy$-plane) midway between the walls by
\begin{equation}\label{eq:Q2D:H approximation}
    \Bu(\Br) \approx -\gamma_0 \tilde{G}(\Br) \BF,
\end{equation}
where $\BF$ is the drag force on the sphere.
The approximation \eqref{eq:Q2D:H approximation} is
valid when $R\ll h$.
%
%
It applies here, because we are concerned with the following scaling regime:
$R\ll h\ll L\ll \ve^{-1}$, where
$\ve^{-1}=|\Bx^2_{_C}-\Bx^1_{_C}|$ is the distance between the two bacteria.

\begin{figure}[h!]
  \begin{center}
  $
  \begin{array}{cc}
    \includegraphics[width=7cm]{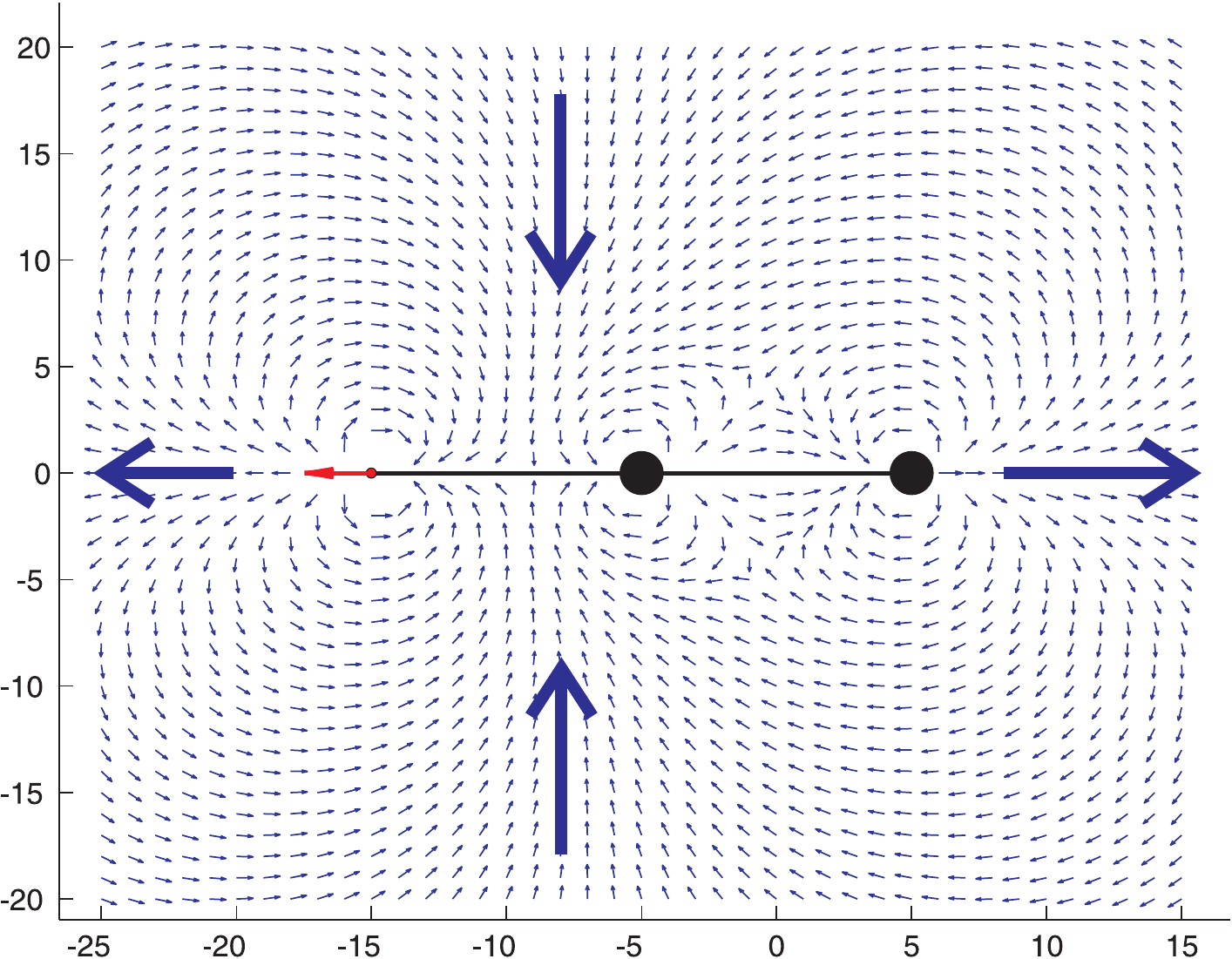}
    &
    \includegraphics[width=7cm]{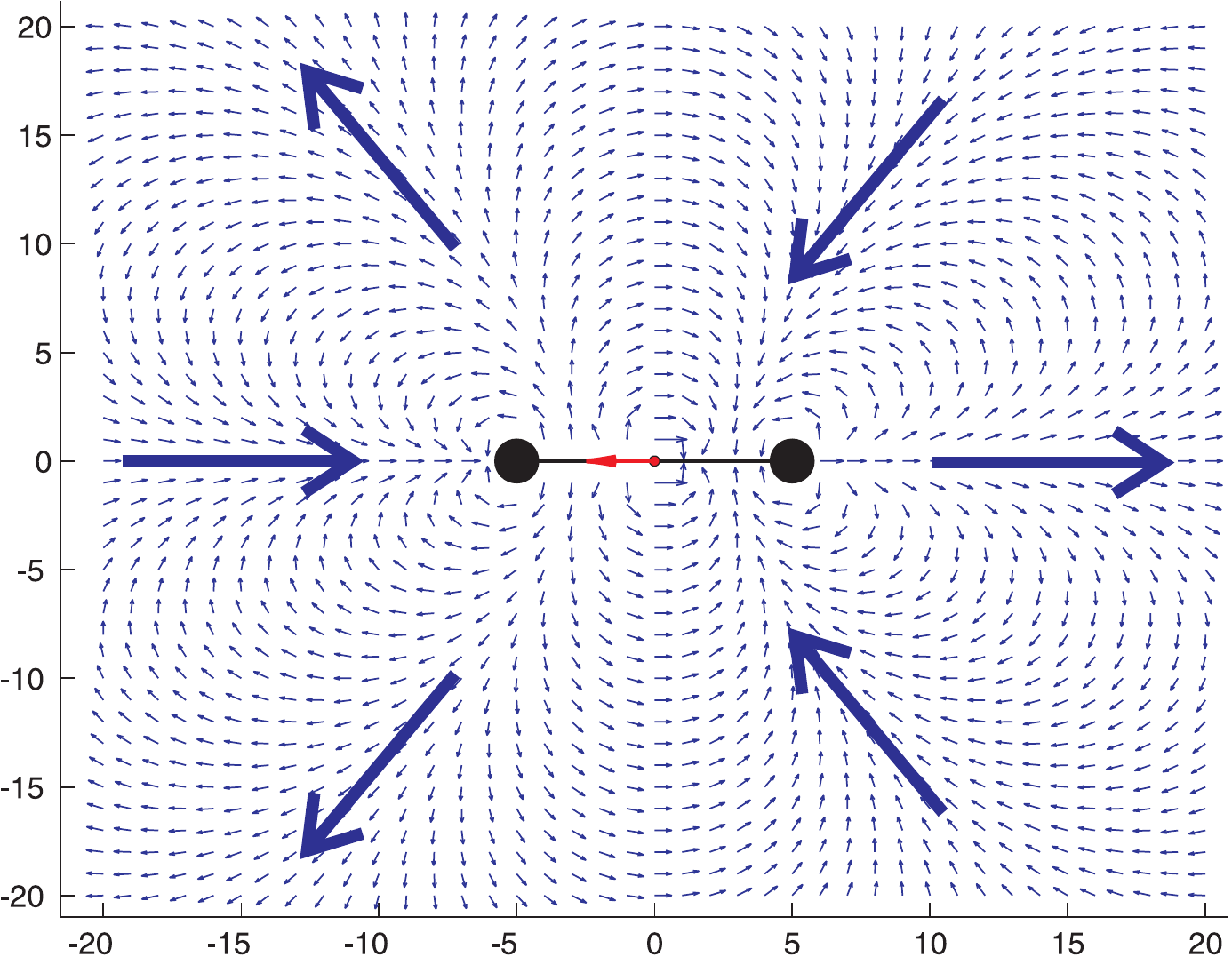}\\
    \text{(a)} & \text{(b)}
  \end{array}
  $
  \end{center}
  \caption{ Velocity field of the Q2D fluid due to a single swimmer:
  \emph{(a)} pusher and
  \emph{(b)} mid-swimmer.
  }
  \label{fig:Q2D:bacteria velocity field}
\end{figure}

%
The solution procedure is exactly the same as for the 3D fluid, except that
$G(\Br)$ is replaced by $\tilde G(\Br)$.
Using a Q2D analog of \eqref{eq:def:approx:u}, we obtain the
velocity field of the Q2D fluid due to a swimming bacterium (see Fig.~\ref{fig:Q2D:bacteria velocity field}).

Note that the asymptotic Green's function $\tilde G(\cdot)$ for the Q2D fluid
is qualitatively different from the Green's function $G(\cdot)$ in a 3D fluid.
For instance, they have different rates of decay:
$G(\Br)\sim |r|^{-1}$ and $\tilde G(\Br)\sim |r|^{-2}$. In addition,  since the shear modes in the Q2D geometry decay
exponentially with the decay rate determined by the spacing between the wall $2 h$, (see, e.g.,  \cite{DiaCuiLinRic05}),
only curl-free ``pressure modes'' decay powerlike survive far away from the origin.
But, most important, $\tilde{G}$ has negative coupling,
$e_1^T\tilde{G}(e_2) e_1<0$.
This means that by applying force to the Q2D fluid in the positive direction along
the $x$-axis some of the fluid will actually be moving in the negative direction
(unlike in 3D, where the coupling is positive and all fluid moves in the positive direction).
In spite of these qualitative differences, the velocity of the fluid
due to a swimming bacterium in Q2D and 3D fluids have similar structures
(compare the bold arrows on Figs.~\ref{fig:bacteria field 1} and \ref{fig:Q2D:bacteria velocity field}).
This similarity of the velocity fields suggests that the dynamics of bacteria may
also be similar for 3D and Q2D fluids. Indeed, we find this to be the case.

Note that the velocity field \eqref{eq:Q2D:point force} due to the point force is curl-free.
Therefore, the velocity fields (a) and (b) are also curl-free as superpositions of
velocity fields of the form \eqref{eq:Q2D:point force}.
It appears that circulation of the velocity fields in Figs.~\ref{fig:Q2D:bacteria velocity field}(a)
and (b)  along closed contours passing through dumbbell balls is not zero, since vector
fields point counterclockwise along some curves.
This situation leads to an apparent contradiction with the Stokes formula.
However, all such curves pass through a singular point of the vector field in the center of the ball,
and the Stokes theorem does not apply (compare to classical electrostatics where all field lines pass through point charges).


\vspace{.5cm}

\textbf{\underline{Asymptotic expressions for velocities}}

Substitute the asymptotic expansion (\ref{eq:def:exp:vc})-(\ref{eq:def:exp:w}) into
the LHS of (\ref{eq:def:w})-(\ref{eq:def:vc}).
Write the velocities of the balls in the RHS of (\ref{eq:def:w})-(\ref{eq:def:vc})
in terms of $\alpha_0$; see (\ref{eq:vH})-(\ref{eq:vT}) and  \eqref{eq:anzats a}.
Expand $G(\cdot)$ in powers of $\ve$ and solve the equations at like powers of $\ve$,
\begin{eqnarray}
  \label{eq:Q2D:order 0}
  O(1)\ :     &\qquad & \omega^1_0=0, \qquad \Bv^1_0=v_0 \tau^i,\\
  \label{eq:Q2D:order 1}
  O(\ve)\ :   &\qquad & \omega^1_1=0, \qquad \Bv^1_1=\mathbf{0},\\
  \label{eq:Q2D:order 2}
  O(\ve^2): &\qquad & \omega^1_2=0, \qquad \Bv^1_2=\mathbf{0},\\
  \label{eq:Q2D:order 3}
  O(\ve^3): &\qquad & \omega^1_3=0, \qquad \Bv^1_3=\frac{f_p L (1-\zeta -2\alpha_0)}{4 \pi \mu}
  \left[
    \begin{array}{c}
      -\cos(2\theta^2-3\phi) \\
      \sin(2\theta^2-3\phi) \\
    \end{array}
  \right],\\
  \label{eq:Q2D:order 4}
  O(\ve^4): &\qquad & \omega^1_4=\frac{3f_p L (1-\zeta -2\alpha_0)}{4 \pi \mu} \sin\left(2\theta^1+2\theta^2-4\phi\right),\\
  \label{eq:Q2D:order 5}
  O(\ve^5): &\qquad & \omega^1_5=\frac{3f_p L^2 (\zeta^2-1)}{2 \pi \mu}\sin\left(2\theta^1+3\theta^2-5\phi\right).
\end{eqnarray}

Next, we analyze the dynamics of two well-separated bacteria in the ``mirror image''
and ``head-to-tail'' configurations (see Fig.~\ref{fig:two cases}) in the Q2D fluid.
We observe that the dynamics of bacteria is qualitatively the same as that of a 3D fluid.
The robustness of the dynamics can be explained by the similarity between the velocity fields
due to swimming bacteria (compare Figs.~\ref{fig:bacteria field 1} and \ref{fig:Q2D:bacteria velocity field}).

\vspace{.5cm}

\noindent\underline{\emph{Dependence of the dynamics of bacteria
on the position $\zeta$ of the propulsion force}}\\
\underline{\emph{A. (``Mirror image'' configuration, Q2D fluid)}}

\begin{figure}[h!]
  \begin{center}
  $
  \begin{array}{ccc}
    \includegraphics[width=4.6cm]{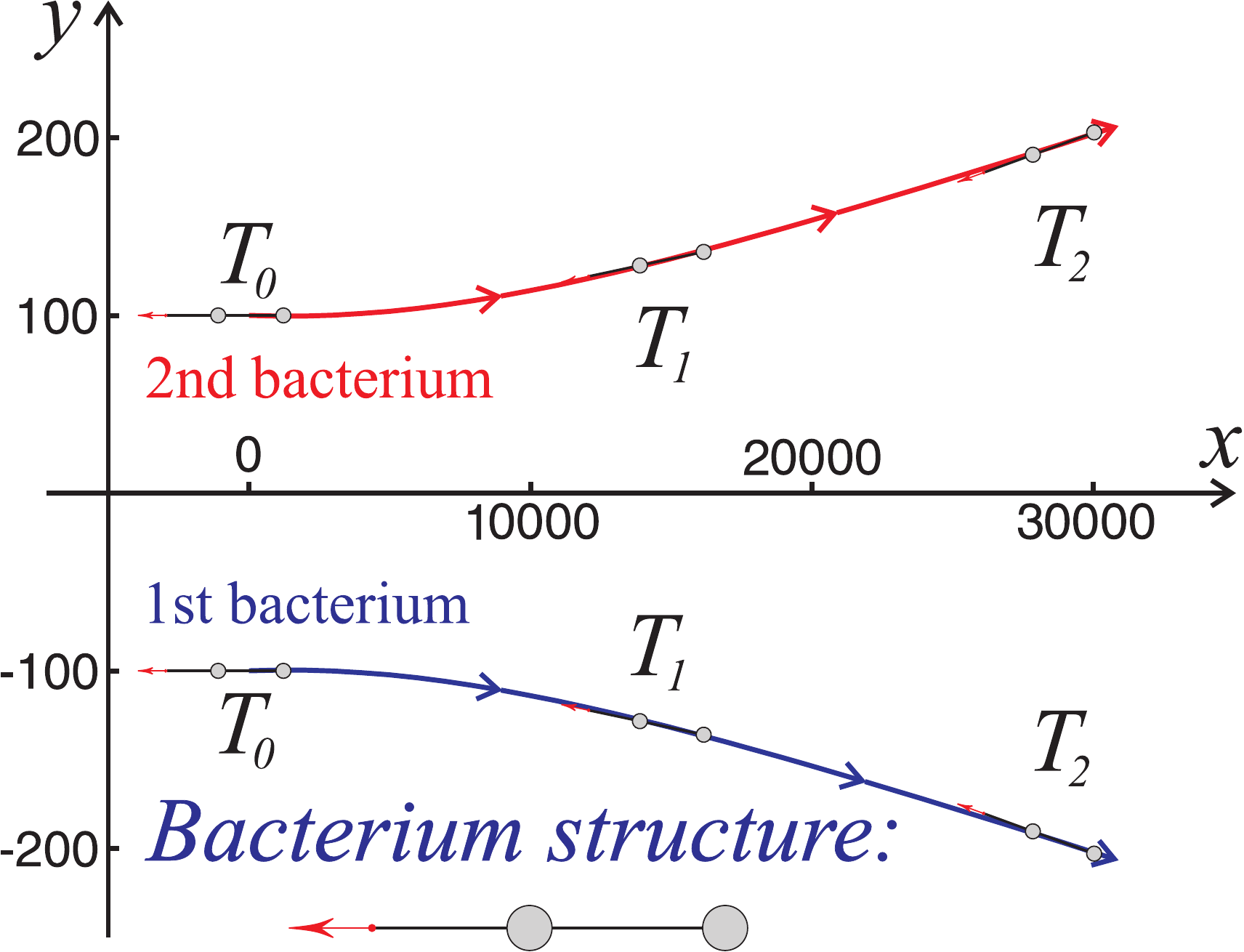}
    &
    \includegraphics[width=4.6cm]{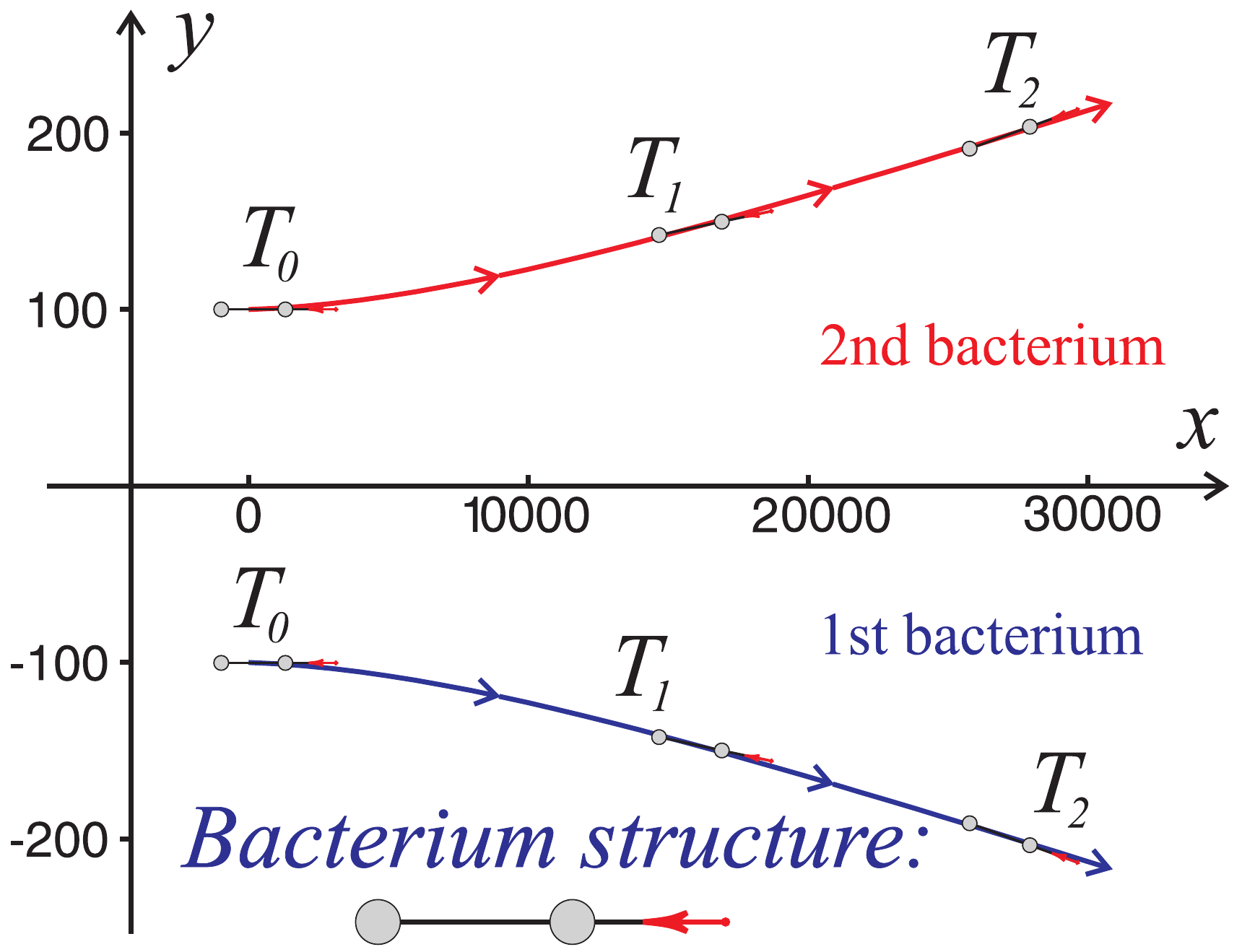}
    &
    \includegraphics[width=4.6cm]{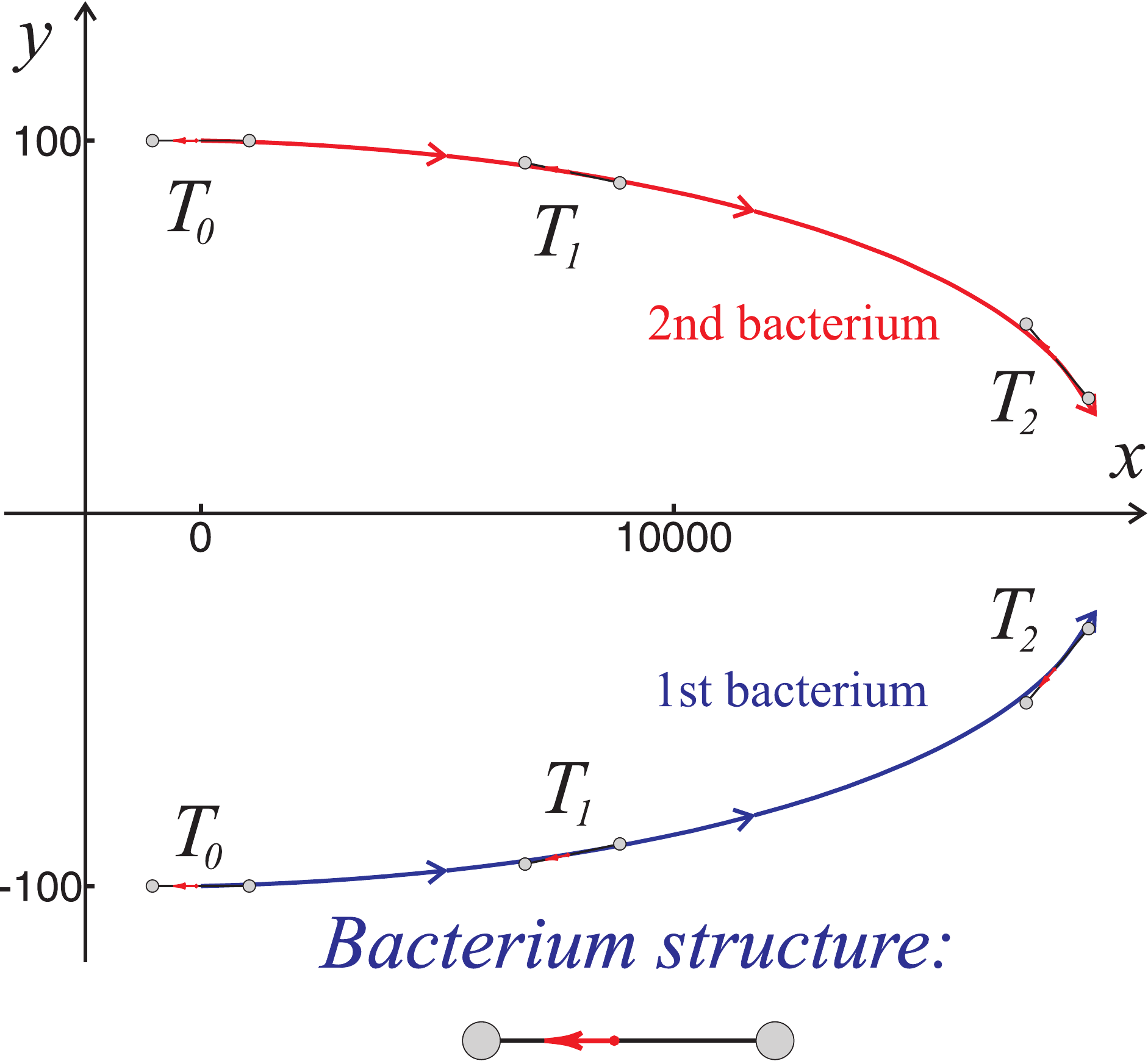}\\
    \text{(a)}
    &
    \text{(b)}
    &
    \text{(c)}
  \end{array}
  $
  \end{center}
  \caption{Trajectories of two swimmers in a Q2D fluid in the ``mirror image'' configuration,
  starting from ($T_0$) parallel orientation ($\theta_1=\theta_2=0$; radius of dumbbell balls $R=1$).
  Initially in (a) through (c) the bacteria are parallel to each other $\theta_1=\theta_2=0$:
  (a) $\zeta<-1$, the outer pushers swim off;
  (b) $\zeta>1$, the outer pullers swim off;
  (c) $|\zeta|<1$, the inner swimmers swim in.
  }
  \label{fig:Q2D:case1:pushers}
\end{figure}

We analyze the dynamics of bacteria depending on the position $\zeta$ of the propulsion force
for the ``mirror image'' configuration of bacteria.
We observe that (as in 3D, see Fig.~\ref{fig:case1:par traj}(c))
when the propulsion force is positioned between
the dumbbell balls ($|\zeta|<1$) the bacteria swim in (see Fig.~\ref{fig:Q2D:case1:pushers}(c)).

Also, (as in 3D, see Fig.~\ref{fig:case1:par traj}(a) and Fig.~\ref{fig:case1:par traj}(b))
when the propulsion force is positioned outside the dumbbell ($|\zeta|>1$)
the bacteria swim off (see Fig.~\ref{fig:Q2D:case1:pushers}(a) and Fig.~\ref{fig:Q2D:case1:pushers}(b)).
\begin{figure}[h!]
    \begin{center}
    \includegraphics[width=9cm]{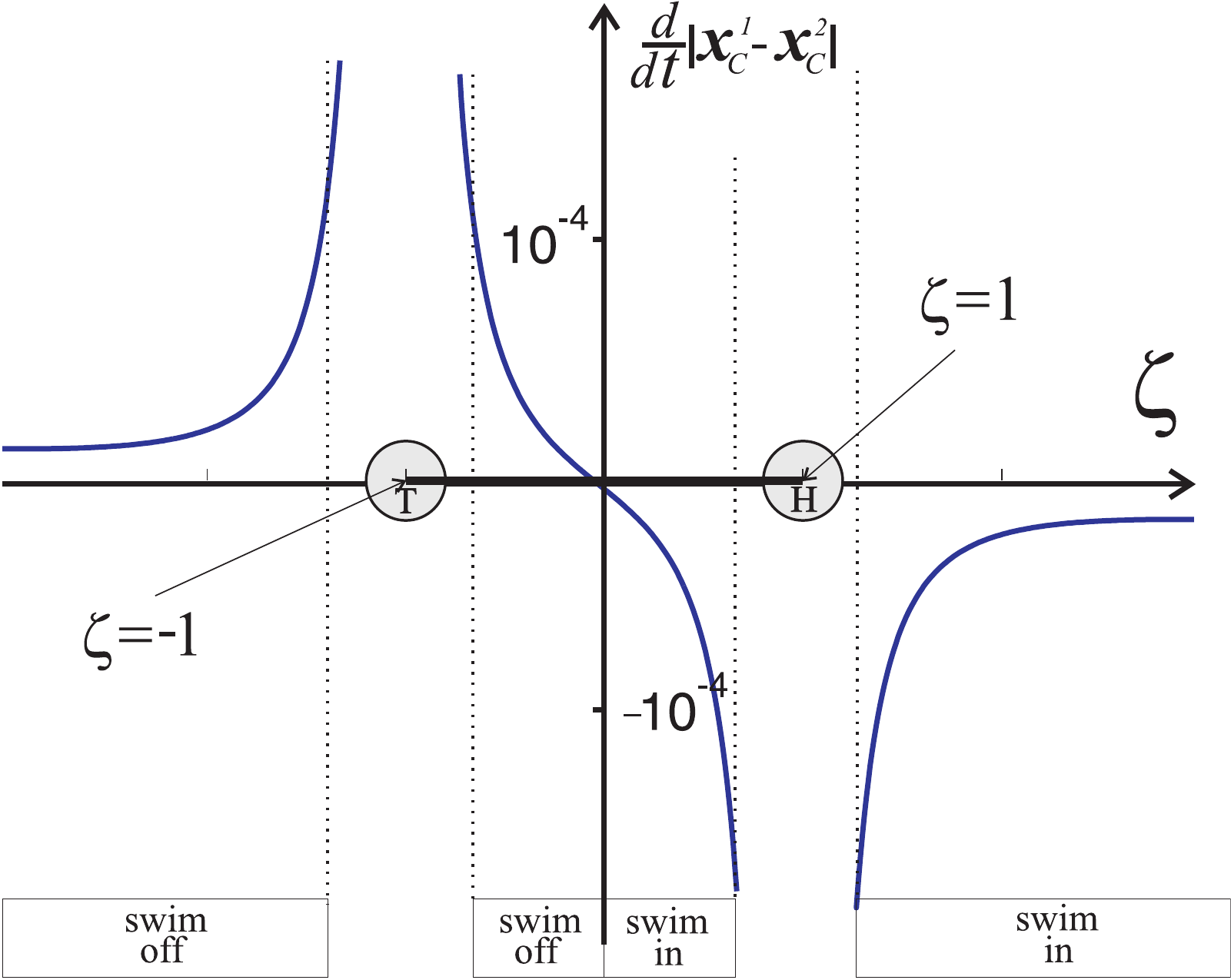}
    \end{center}
  \caption{
  Dependence of the relative velocities of two Q2D swimmers in the ``head-to-tail'' configuration
  on $\zeta$, which determines the position of the propeller.
  Pushers ($\zeta<0$) swim off, and
  pullers ($0<\zeta$) swim in.
  The asymptotic technique used by us cannot be applied
  to the uncharacterized regions (between dashed lines) of $\zeta$ close to $\pm 1$.
  }
  \label{fig:Q2D:case2:stag distance}
\end{figure}

\underline{B. \emph{``Head-to-tail'' configuration, Q2D fluid}}

For two bacteria in the ``head-to-tail'' configuration in the Q2D fluid we observe the same dynamics as
for the 3D fluid -- pushers ($\zeta<0$) swim off and pullers ($0<\zeta$) swim in
(see Fig.~\ref{fig:Q2D:case2:stag distance}).

\section{Conclusions}\label{sect:Conclusions}

In this paper we studied the hydrodynamic interaction between two microscopic swimmers, 
modeled as self-propelled dumbbells, in two distinct settings: 
three-dimensional and quasi-two-dimensional fluid domains.
The interaction in a three-dimensional fluid domain models the interaction of swimmers in the bulk
(away from the walls of the container), while
the interaction in a quasi-two-dimensional fluid domain models the interaction of swimmers 
in a thin film.  Qualitatively, models in both settings produced the same results, 
thus suggesting that the hydrodynamic interaction of a pair of swimmers is robust under the change of geometry of the fluid domain.

At the same time, 
the shape of the swimmer, that is, the position ($\zeta$) of the effective propulsion force, proved to have a critical
effect on the character of the hydrodynamic interaction of swimmers.
In the ``mirror-image'' configuration the dynamics of swimmers differentiates inner ($|\zeta|<1$) and outer swimmers ($|\zeta|>1$).
Inner swimmers ($|\zeta|<1$) in the ``mirror-image'' configuration experience a swim in, approaching each other and
perfectly matching the swim in experimentally observed in \cite{AraSokGolKes07} 
for the rod-shaped bacterium {\it Bacillus subtilis}, which has multiple flagella distributed over the cell surface.
Unlike inner swimmers, outer swimmers ($|\zeta|>1$) in the ``mirror-image'' configuration experience a swim off, due to outward rotation.

In the ``head-to-tail'' configuration the dynamics of the swimmers differentiates pushers ($\zeta<0$) and 
pullers ($\zeta>0$).
Pushers ($\zeta<0$) in the ``head-to-tail'' configuration experience a swim off;
that is, while they remain oriented along the same straight line, the distance between them gradually increases.
Unlike pushers, pullers ($\zeta>0$) in the ``head-to-tail'' configuration experience a swim in;
that is, while they remain on the same straight line, the distance between them gradually decreases. 
Moreover, for pushers the ``head-to-tail'' configuration is not stable, whereas for pullers it is stable. 
Thus, our model predicts a formation of ``head-to-tail'' structures by pullers and no such structures for pushers.

The surprising sensitivity of the observed hydrodynamic interaction of swimmers
to the flagellum position
(more generally to the structure of the propulsion apparatus) 
and, therefore, to the structure and the shape of the swimmer 
explains the wide range of behaviors exhibited by microorganisms 
(see, e.g., \cite{Riedel05} for a study of a sperm cell with a very long flagellum and
\cite{RufNul85} for a study of algae that pull themselves forward with flagella
positioned in the forward part of the body)
and different models of microscopic swimmers, such as dumbbells
\cite{HerStoGra05},
squirmers
\cite{IshiPed07},
self-locomoting rods
\cite{SainShel07},
and three-sphere swimmers \cite{Pooley07}.

Further refinements of our model are keenly needed. In particular, our calculations are conducted in the dilute limit, where the distance between the swimmers is large compared to their size.
However, as we demonstrated, pushers have a tendency to converge, thus eventually violating this approximation.
Therefore, nontrivial regularizations of the interaction at small distances using,  possibly, lubrication forces and
hard-core repulsion must be included into the model in order to obtain agreement with experiments and simulations.
Further,  at high concentration, deviations from the pairwise interaction may also become important, especially
because hydrodynamic forces decay very slowly in the three-dimensional geometry of the sample.


\section*{Acknowledgments}

The work of Igor Aranson and Dmitry Karpeev was supported by US DOE contract DE-AC02-06CH11357.
The work of Vitaliy Gyrya and Leonid Berlyand was supported by
DOE grant DE-FG02-08ER25862
and NSF grant DMS-0708324.

   %
\newpage
\appendix
\section{Basic Stokes solutions}\label{sect:Stokes}

\subsection{Point force}
\label{app:point_force}
  The velocity field due to a point force $\BF$ in an unbounded fluid domain is
  \begin{equation}\label{eq:sol:point force}
    \Bu(\Bx) = G(\Bx) \cdot \BF,\qquad
    G(\Bx) =\frac{1}{8\pi \mu |\Bx|}\left(\BI+\frac{\Bx \Bx^T}{|\Bx|^2}\right).
  \end{equation}
  Tensor $G$ (along with a suitable pressure tensor $P$)
solves the Stokes problem with a point force
\begin{equation*}
    \left\{
    \begin{array}{l}
      \mu \lap G = \D P - \delta(\Bx)\\
      \diiv(\Bu)=0
    \end{array}
    \right..
\end{equation*}
  Therefore, it is the fundamental solution to the above problem, given in components by
  $$
        G_{ij}(\Bx) = \frac{1}{8\pi \mu |\Bx|}\left( \delta_{ij} + \frac{x_i x_j}{|\Bx|^2}\right),
  $$
  with the corresponding pressure, a vector, given by (up to an additive constant)
  $$
        P_i(\Bx) = \frac{1}{4 \pi} \frac{x_j}{|\Bx|^3}.
  $$
  The stress tensor corresponding to $G$ and $P$ is a triadic $\Sigma$:
  $$
        \Sigma_{ijk} = -P_j \delta_{ik} + \frac{\mu}{2}\left(G_{ij,k} + G_{kj,i}\right) =
        -\frac{3}{4 \pi} \frac{x_i x_j x_k}{|\Bx|^5}.
  $$
  For more details see \cite{KimKar91}.

\subsection{Swimming ball}
\label{app:swimming_ball}

 A ball of radius $R$ moving with a constant velocity $\Bv$
  through an unbounded fluid domain creates the velocity field:
  \begin{eqnarray}\label{eq:sol:ball}
    &\Bu(\Bx) = H(\Bx; R)\, \Bv,\qquad \qquad
    H(\Bx;R) = \frac{3 R}{4 r} \left[\alpha \BI + \beta\, \Bn\Bn^T\right],&\\
    \nonumber &\alpha  = 1+\frac{R^2}{3 r^2},\qquad
    \beta = 1-\frac{R^2}{r^2},\qquad  r=|\Bx|,\qquad \Bn = \frac{\Bx}{r},&
  \end{eqnarray}
  where $\BI$ is the identity matrix and $\left(\Bn\Bn^T\right) \Bv = (\Bv\cdot\Bn) \Bn$ is the dyadic product.

  Away from the origin ($r\gg R$)
$$
        \qquad H(\Bx; R) \approx \gamma_0\ G(\Bx),\qquad \gamma_0 = 6\pi \mu R,
$$
where $\gamma_0$ is the inverse mobility of the ball, characterizing the applied force necessary
to generate a steady translational velocity of unit magnitude.


\subsection{Stokes law for drag}
\label{app:stokes_law}
  The drag force from the fluid of viscosity $\mu$ on a ball of radius $R$, moving with a velocity $\Bv$
  through unbounded fluid is
  \begin{equation}
    \BF = -\gamma_0 \Bv.
    \label{eq:drag}
  \end{equation}
  More generally, suppose that a ball is added to given an initial background flow $\overline \Bu$
  and that under the influence of external forces the ball undergoes a steady tranlational motion of the ball with velocity $\Bv$.
  The Stokes law for drag states that the accompanying drag force $\BF$ on the ball is proportional to the difference
  of the velocity of the ball and the velocity of the background flow, which would exist at the
  location of the ball $\overline \Bx$ in its absence:
  \begin{equation}
        \BF = -\gamma_0 \left(\Bv - \overline \Bu(\overline \Bx)\right).
        \label{eq:stokes_drag}
  \end{equation}
  Since in the Stokes framework the drag on the ball must be balanced by the applied forces on the particle,
  (\ref{eq:stokes_drag}) provides a means of calculating the net applied force that results in a given
  translation velocity $\Bv$.

  The Stokes law is an approximation to Fax\'en's first law \cite{KimKar91}:
  \begin{equation}\label{eq:faxen}
        \BF = -\gamma_0 \left(\Bv - \overline \Bu(\overline \Bx)\right) +
        \gamma_0 \frac{R^2}{6} \nabla^2 \overline \Bu(\overline \Bx).
  \end{equation}
  If the background flow is due to a point force or another translating sphere at $\Bx$, far from $\overline \Bx$,
  then it follows from (\ref{eq:sol:point force}) and (\ref{eq:sol:ball}) that the gradient is small -- $\sim
\frac{1}{|\Bx-\overline\Bx|^2}$.  In this case the Stokes law (\ref{eq:stokes_drag}) is a good approximation to (\ref{eq:faxen}).

\subsection{Q2D Green's function}\label{sect:Stokes:Q2D}


Take formula (51) in \cite{LirMoc76a},
\begin{eqnarray}
\label{eq:formula 51}
    \Bu^k_j &\approx &
    -\frac{3H}{\pi\mu}\frac{x_3}{H}\left(1-\frac{x_3}{H}\right)\frac{h}{H}
    \left(1-\frac{h}{H}\right)\frac{1}{\rho^2}
    \left[\frac{1}{2} \delta_{\a\beta}-\frac{r_\a r_\beta}{\rho^2}\right]
    \delta_{j\a}\delta_{k\beta}+\nonumber\\
    & & +
    \delta_{j3}\delta_{k3} O\left(\rho^{-\frac{1}{2}} e^{-\rho y_1/H}\right)
    +
    (\delta_{j3}\delta_{k\a} + \delta_{k3}\delta_{j\a})
    O\left(\frac{r_\a}{\rho} \rho^{-\frac{1}{2}} e^{-\rho y_1/H}\right)+\nonumber\\
    & & +
    \delta_{j\a}\delta_{k\beta}
    \left[
    O\left(\frac{r_\a}{\rho}\frac{r_\beta}{\rho} \rho^{-\frac{1}{2}} e^{-\rho y_1/H}\right)
    +
    O\left(\frac{r_\a}{\rho}\frac{r_\beta}{\rho} \rho^{-\frac{1}{2}} e^{-\rho \pi/H}\right)
    \right],
\end{eqnarray}
where $y_1\approx 4.2$, and rewrite it in our notations.
The point force is applied midway between the walls of the film.
Replace $h=\frac{1}{2}H$;
here $H$ is thickness of the film (replace by $h$).
Assume $k=1$, that is force is applied along $e_1$.
Here $\rho$ is the radius vector from point force (replace by $r=|\Br|$).
Replace $x_3$ by $z$.

Performing the above changes, we obtain
\begin{eqnarray}
    \label{eq:formula 51 adapted}
    \Bu^1_j &\approx&
    -\frac{3 z}{4\rho^2\pi\mu}\left(1-\frac{z}{2h}\right)
    \left[\frac{1}{2} \delta_{\a\beta}-\frac{r_\a r_\beta}{\rho^2}\right]
    \delta_{j\a}\delta_{1\beta}+\nonumber\\
    & & +
    \delta_{j3}\delta_{13} O\left(\rho^{-\frac{1}{2}} e^{-\rho y_1/(2h)}\right)
    +
    (\delta_{j3}\delta_{1\a} + \delta_{13}\delta_{j\a})
    O\left(\frac{r_\a}{\rho} \rho^{-\frac{1}{2}} e^{-\rho y_1/(2h)}\right)+\nonumber\\
    & & +
    \delta_{j\a}\delta_{1\beta}
    \left[
    O\left(\frac{r_\a}{\rho}\frac{r_\beta}{\rho} \rho^{-\frac{1}{2}} e^{-\rho y_1/(2h)}\right)
    +
    O\left(\frac{r_\a}{\rho}\frac{r_\beta}{\rho} \rho^{-\frac{1}{2}} e^{-2\rho \pi/(2h)}\right)
    \right].
\end{eqnarray}
Note that only the first term in \eqref{eq:formula 51 adapted} does not decay exponentially in $\rho$:
\begin{equation}\label{eq:formula 51 adapted 2}
\begin{split}
    \Bu^1_j &\approx
    -\frac{3}{4\pi\mu}z\left(1-\frac{z}{2h}\right)
    \frac{1}{\rho^2}
    \left[\frac{1}{2} \delta_{\a\beta}-\frac{r_\a r_\beta}{\rho^2}\right]
    \delta_{j\a}\delta_{1\beta}
    =\\
    &=
    \frac{3}{4\pi\mu}z\left(1-\frac{z}{2h}\right)
    \frac{1}{\rho^2}
    \left[\frac{r_\a r_\beta}{\rho^2}\delta_{j\a}\delta_{1\beta}
    -\frac{1}{2} \delta_{\a\beta}\delta_{j\a}\delta_{1\beta}\right]=\\
    &=
    \frac{3}{4\pi\mu}z\left(1-\frac{z}{2h}\right)
    \frac{1}{\rho^2}
    \left[\frac{r_j r_1}{\rho^2}-\frac{1}{2} \delta_{j1}\right].
\end{split}
\end{equation}
Rewriting $\Bu^3$ in components, we have
\begin{equation*}
    \Bu^1 \approx
    \left[
      \begin{array}{c}
        u_1 \\
        u_2 \\
        u_3 \\
      \end{array}
    \right]=
    \frac{3}{4\pi\mu}
    z\left(1-\frac{z}{2h}\right)
    \frac{1}{\rho^2}
    \left[
      \begin{array}{c}
        \frac{r_1 r_1}{\rho^2}-\frac{1}{2} \delta_{11} \\
        \frac{r_2 r_1}{\rho^2}-\frac{1}{2} \delta_{21} \\
        \frac{r_3 r_1}{\rho^2}-\frac{1}{2} \delta_{31} \\
      \end{array}
    \right]
    =
    \frac{3}{4\pi\mu}
    z\left(1-\frac{z}{2h}\right)
    \frac{1}{\rho^2}
    \left[
      \begin{array}{c}
        \frac{x^2}{\rho^2}-\frac{1}{2} \\
        \frac{x y}{\rho^2} \\
        \frac{x z}{\rho^2} \\
      \end{array}
    \right].
\end{equation*}
Note that $\rho^2\approx x^2+y^2$.
Hence, \eqref{eq:formula 51 adapted 2} takes the form
\begin{equation}\label{eq:formula 51 adapted 4}
    \Bu^1 \approx
    \frac{3}{8\pi\mu}
    z\left(1-\frac{z}{2h}\right)
    \frac{1}{\rho^4}
    \left[
      \begin{array}{c}
        x^2-y^2 \\
        2x y \\
        2 x z \\
      \end{array}
    \right].
\end{equation}

Since $|z|<h\ll 1$, one has an approximation for the Green's function in Q2D:
\begin{equation}\label{eq:Q2D G approximation}
    \tilde{G}(\Br) \approx
    f(z)\left[
          \begin{array}{ccc}
            x^2-y^2 & 2xy       & 0 \\
            2xy     & y^2-x^2   & 0 \\
            0       & 0         & 0 \\
          \end{array}
        \right],\qquad
    \text{where }\ f(z) = \frac{3}{8\pi\mu}
    z\left(1-\frac{z}{2h}\right).
\end{equation}
Here $f(z)$ satisfies no-slip boundary conditions: $f(z=0)=0$ and $f(z=2h)=0$.


\section{Existence and uniqueness}
\label{app:existence}
In this section we prove the existence and uniqueness of solutions to the system (\full_system)
under different assumptions on the regularity of the propulsion forces and the size of the container.
The most restrictive case of smooth forces and a bounded container yields the clearest proof
that is essentially classical, but contains a few novel features.  The other cases refine
the argument in the case of point forces and an unbounded container.

For the sake of clarity we consider the case of a single swimmer.  The extension to a multiswimmer
system is straight-forward.  We always assume that the swimmers do not overlap,
which in particular excludes the overlap of any propeller with any head or tail.

\subsection{Regular case}
To elucidate the main issues in the existence and uniqueness proof we initially
consider the case of regular problem data.
Assuming that the propeller force density is smooth and the container $\Omega$
is bounded, the argument is an adaptation of the classical
techniques based on coercivity to deduce the existence of weak solutions,
followed by an application of elliptic regularity results.
The novel feature is the presence of a boundary condititon on the forces and torques
(\balance)
and the consequent use of Korn's inequality in place of the usual Poincare's inequality.

\subsubsection{Space of admissible flows}
In order to obtain an appropriate weak formulation of the problem that incorporates the balance conditions
(\balance), we have to circumscribe the space in which the solutions will be sought.
Consider the following space
\begin{equation}
    \mathcal{V} = \Big\{\Bu\in \left(\mathcal{D}(\overline{\Omega}_F)\right)^3 \big|
    \ \Bu(\Bx) = \Bv^\Bu+\omega^\Bu\times(\Bx-\Bx_{_C}) \text{ for }\Bx\in\d B_{_H}\cup \d B_{_T},\quad
    \Bu|_{\partial \Omega} = 0, \Bv^\Bu, \omega^{\Bu}\in\R^3 \Big\},
    \label{eq:V_1}
\end{equation}
where $\mathcal{D}\left(\overline{\Omega}_F)\right)$ is, as usual, the calss of restrictions to $(\Omega_F)$
of $C^\infty_0(\R^3)$ -- smooth functions on $R^3$ with compact support.

An equivalent definition of $\mathcal{V}$ is through constrains
\begin{equation}
    \Bu(\Bx) = \Bv^\Bu_{_H}+\omega^{\Bu}\times(\Bx-\Bx_{_H})\qquad \Bx\in \d B_{_H},\qquad
    \Bu(\Bx) = \Bv^\Bu_{_T}+\omega^{\Bu}\times(\Bx-\Bx_{_T})\qquad \Bx\in \d B_{_T},
    \label{eq:V_2}
\end{equation}
where $\Bv^{\Bu}_{_H},\Bv^{\Bu}_{_T}\in\R^3$ are given in terms of $\Bv^{\Bu}_{_C}$ and
$\omega^{\Bu}_{_C}$ by
\begin{equation}
    \label{eq:remark uh through uc}
    \Bv^{\Bu}_{_H} = \Bv^{\Bu}_{_C}+\omega\times(\Bx_{_H}-\Bx_{_C}),\qquad
    \Bv^{\Bu}_{_T} = \Bv^{\Bu}_{_C}+\omega\times(\Bx_{_T}-\Bx_{_C}),\qquad
    \Bv^\Bu_{_C} = \frac{\Bv^\Bu_{_H} + \Bv^\Bu_{_T}}{2}.
\end{equation}
or after substituting $(\Bx_{_H}-\Bx_{_C})=L\tau$ and $(\Bx_{_H}-\Bx_{_C})=-L\tau$
\begin{equation*}
    \Bv^{\Bu}_{_H} = \Bv^{\Bu}_{_C}+L\omega\times\tau,\qquad
    \Bv^{\Bu}_{_T} = \Bv^{\Bu}_{_C}-L\omega\times\tau,
\end{equation*}
Now define \mbox{$V = \overline{\mathcal{V}}$} as
the closure of $\mathcal{V}$ in $H^1(\Omega_F)$.

It is important to note that $V$ is not empty: the boundary data of the form
\begin{equation}
        \Bw(\Bx) = \left\{
        \begin{array}{lr}
                \Bv^\Bw+\omega^\Bw\times(\Bx-\Bx_{_C}),& \Bx \in \d \Omega_B\\
                0,&\Bx \in \d \Omega
        \end{array}
        \right.
        \label{eq:boundary_data}
\end{equation}
can be continued to $\BW(\Bx)$ defined on all of $\Omega_F = \Omega \setminus \Omega_B$.
To show this we can apply the standard theorem found, for instance, in \cite{Galdi94},
as soon as we show that the boundary data of the form \eqref{eq:boundary_data}
satisfy the compatibility condition:
$$
        \int_{\partial \Omega_F} \Bw(\Bx) \cdot \Bn(\Bx) d\Bx =
        -\int_{\partial \Omega} \Bw(\Bx) \cdot \Bn(\Bx) d\Bx +
        \sum_{i,*}\int_{\partial B^i_*} \Bw(\Bx) \cdot \Bn(\Bx) d\Bx,
$$
where for consistency with the rest of the paper we denoted by $\Bn$ is the normal pointing
\emph{into the fluid} (i.e., the negative of the outward normal).
The first integral above is clearly $0$, while the second integral vanishes by the
divergence theorem.
Indeed, on the boundary $\partial B^i_H \bigcup \partial B^i_T$ of a fixed $i-th$ dumbbell
the field is $\Bw(\Bx) = \Bv^i+\omega^i\times(\Bx-\Bx^i_{_C})$, which clearly can
be continued into the interior of the dumbbell where it is divergence-free.

Note that $\Bv^\Bw,\ \Bomega^\Bw$ as well as $\Bv^\Bw_H$ and $\Bv^\Bw_T$ are uniquely
defined for any $\Bw$ in $V$.
For $\Bu \in \mathcal{V}$ that solves (\full_system) (a strong solution)
$\Bv^\Bu$ and $\Bomega^\Bu$ are
the linear and angular velocities $\Bv$ and $\Bomega$ entering into (\ref{eq:rigidity}).
Taking a fixed extension $\BW$ for each distinct boundary field defined by a pair
$(\Bv^\Bw,\,\Bomega^\Bw) \in \R^3 \oplus \R^3 = \R^6$
we can inn fact showt that $V$ is isomorphic to the direct sum $\R^6 \oplus V_0$.
Here $V_0$ is the closure in $H^1_0(\Omega_F)$ of divergence-free functions from
$\mathcal{D}(\Omega_F)$.  The isomorphism is established by sending $\Bw \in V$
to $\left((\Bv^\Bw,\,\Bomega^\Bw),\,\Bw - \BW\right) \in \R^6 \oplus V_0$.

\subsubsection{Weak form}
In this section we derive the weak form of (\full_system) on $V$.
To this end, assume there is a strong solution $\Bu \in \mathcal{V}$ to (\full_system).
Starting with the usual formulation of the Stokes' equation in terms of the Laplacian to derive
the weak form eventully leads to bilinear boundary terms, which are hard to bound from below
to show coercivity.  Instead, adding  $0=\mu\D (\Div \Bu)=\mu \Div (\D \Bu)^T$
to the first equation in (\ref{eq:stokes})
we obtain an equivalent formulation in terms of the symmetrized gradient $D(\Bu)$:
\begin{equation}
    2\mu \Div \big( D(\Bu)\big) = \D p - \BF.
    \label{eq:stokes_sym}
\end{equation}

Multiplying \eqref{eq:stokes_sym} by $\Bw\in \mathcal{V}$
\begin{equation*}
    \int_{\Omega_F} \left(2\mu \Div \big( D(\Bu)\big) + \BF - \D p\right)\cdot \Bw\ dx = 0
\end{equation*}
and integrating by parts yields
\begin{gather*}
  -2\mu \int_{\Omega_F} D(\Bu):\D\Bw\ dx + 2\mu\int_{\d \Omega_F}\Bn\cdot D(\Bu)\cdot \Bw\ dx +
  +\int_{\Omega_F}\BF \cdot \Bw\ dx+\int_{\Omega_F} p\Div (\Bw)\ dx - \int_{\d\Omega_F} p \Bw\cdot \Bn\ dx
  = 0.
\end{gather*}
Using the incompressibility condition, the symmetry of $D(\Bu)$ and rearranging terms we obtain
\begin{equation}
    4\mu \int_{\Omega_F} D(\Bu):D(\Bw)\ dx
    =
    \int_{\Omega_F}\BF \cdot \Bw\ dx
    +
    \int_{\d \Omega_F}\Bn\cdot \left[ 2\mu D(\Bu)-p \BI\right]\cdot \Bw\ dx.
    \label{eq:form_1}
\end{equation}

At this point the idea is to extend \eqref{eq:form_1} to $V$
thereby
obtaining a weak form of (\full_system).  However, the presence of a bilinear boundary integral term will
resist easy lower bounds needed to show the coercivity of the problem.
This difficulty is circumvented by using the balance conditions (\balance) to eliminate $\Bu$
from the integral, converting it into a linear functional of $\Bw$.
To this end, note that the boundary integral contains a product of $\Bw$ with the boundary tractions
$\Bn \cdot \sigma(\Bu)$ generated flow $\Bu$.

Rewrite the boundary integral in \eqref{eq:form_1} using
the boundary conditions \eqref{eq:V_1} and \eqref{eq:V_2} on $\Bw$ featuring in the definition of
$\mathcal{V}$:
\begin{eqnarray}
  \label{eq:boundary int:1}
  \int_{\d \Omega_F}\Bn\cdot \sigma(\Bu)\cdot \Bw\ dx
  &=&\\
  \label{eq:boundary int:2}
  =\int_{\d B_{_H}\cup\d B_{_T}}\Bn\cdot \sigma(\Bu)\cdot (\Bv^\Bw+
  \omega^\Bw\times(\Bx-\Bx_{_C}))\ dx
  &=&\\
  \label{eq:boundary int:3}
  =\int_{\d B_{_H}}\Bn\cdot \sigma(\Bu)\cdot (\Bv^\Bw+
  \omega^{\Bw}\times(\Bx_{_H}-\Bx_{_C})+
  \omega^{\Bw}\times(\Bx-\Bx_{_H}))\ dx
  &+&\\
  \label{eq:boundary int:4}
  +\int_{\d B_{_T}}\Bn\cdot \sigma(\Bu)\cdot (\Bv^\Bw+
  \omega^\Bw\times(\Bx_{_T}-\Bx_{_C})+
  \omega^\Bw\times(\Bx-\Bx_{_T}))\ dx
  &=&\\
  \label{eq:boundary int:5}
  =\BF_{_H}\cdot\Bv^\Bw_{_H} + \BT_{_H}\cdot \omega^\Bw+
  \BF_{_T}\cdot\Bv^\Bw_{_H}+  \BT_{_T}\cdot \omega^\Bw.
\end{eqnarray}
Rewrite \eqref{eq:boundary int:5} as
\begin{eqnarray}
  (\BF_{_H}+\BF_{_T})\cdot \left(\frac{\Bv^\Bw_{_H}+\Bv^\Bw_{_T}}{2}\right)+
  (\BF_{_H}-\BF_{_T})\cdot \left(\frac{\Bv^\Bw_{_H}-\Bv^\Bw_{_T}}{2}\right)+
  (\BT_{_H}+\BT_{_T})\cdot \omega^\Bv,
\end{eqnarray}
where $\BF_*$ and $\BT_*,\ * = H,T$ are the forces and torques associated with $\Bu$ as defined the previous
section.
We will show now that the first term defines a continuous linear functional on $\Bw$,
and the last two terms vanish due to the conditions on
$\BF_*,\BT_*$ and $\Bv^\Bw_*$.
Indeed, since $\tau\cdot(\Bv^\Bw_{_H}-\Bv^\Bw_{_T})=0$, (i.e., $(\Bv^\Bw_{_H}-\Bv^\Bw_{_T})\perp\tau$), we have
\begin{equation*}
    (\Bv^\Bw_{_H}-\Bv^\Bw_{_T})=\tau\times\left(\tau\times(\Bv^\Bw_{_H}-\Bv^\Bw_{_T})\right) .
\end{equation*}
Hence, from the relation
\eqref{eq:remark uh through uc}
and the
balance of torques \eqref{eq:torque_balance}
we have
\begin{gather*}
  (\BF_{_H}-\BF_{_T})\cdot \left(\frac{\Bv^\Bw_{_H}-\Bv^\Bw_{_T}}{2}\right)+
  (\BT_{_H}+\BT_{_T})\cdot \omega^\Bv =
  (\BF_{_H}-\BF_{_T})\cdot L\tau\times\omega^\Bv + (\BT_{_H}+\BT_{_T})\cdot \omega^\Bv = \\
  \label{eq:togivebalanceoftorques}
  \big[(\BF_{_H}-\BF_{_T})\times L\tau+\BT_{_H}+\BT_{_T}\big]\cdot \omega^\Bv = 0.
\end{gather*}
Finally, from the balance of forces
\eqref{eq:force_balance}
we have
\begin{equation}\label{eq:use balance}
    \left(\frac{\Bv^\Bw_{_H}+\Bv^\Bw_{_T}}{2}\right)\cdot (\BF_{_H}+\BF_{_T})=
    -\left(\frac{\Bv^\Bw_{_H}+\Bv^\Bw_{_T}}{2}\right)\cdot \BF_{_P} = -\Bv^\Bw_{_C} \cdot \BF_{_P}.
\end{equation}
which defines a continuous linear functional of $\Bw$ in terms of the fixed total force $\BF_P$.

Thus the solution of the Stokes equation \eqref{eq:stokes} satisfies the following variational problem
\begin{equation}\label{eq:variational 1}
    4a(\Bu,\Bw) = b(\Bw)\qquad \forall \Bw\in \mathcal{V},
\end{equation}
where
\begin{gather}
    a(\Bu,\Bw):= \mu\int_{\Omega_F} D(\Bu): D(\Bv)\ dx
    \label{eq:def:bilinear form a}\\
    b(\Bw):= \int_{\Omega_F}\left(\BF \cdot \Bw\ dx - \BF_{_P}\cdot\Bv^\Bw_{_C}\right).
    \label{eq:def:L}
\end{gather}
We want to emphasize that it was the use of the symmetrized gradient $D(\Bu)$ in place of the usual gradient
that lead to a boundary integral in terms of tractions, which, apart from having a clear physical meaning,
enabled elimination of $\Bu$ with the help of the force balance conditions.


The minimization problem corresponding to the variational problem \eqref{eq:variational 1} is
\begin{equation}\label{eq:minimization}
    \min_{\Bu\in\mathcal{Z}} E[\Bu],
\end{equation}
where the energy functional is
\begin{equation}\label{eq:Energy E}
    E[\Bu] = 2\,a(\Bu,\Bu)-b(\Bu).
\end{equation}
The quadratic term $a(\Bu,\Bu)$ is the usual viscous dissipation rate and the
linear term $b(\Bu)$ represents the work of the forces in the fluid and on its boundary.
The interpretation of $b(\Bu)$ as the work done by the forces
becomes clearer once we rewrite it as:
\begin{equation*}
    b(\Bu) = \int_{\Omega_F}\BF \cdot \Bv\ dx
    +
    \BF_{_H}\cdot\Bv^\Bw_{_H}+
    \BT_{_H}\cdot\omega^\Bw+
    \BF_{_T}\cdot\Bv^\Bw_{_T}+
    \BT_{_T}\cdot\omega^\Bw.
\end{equation*}

\subsubsection{Existence, uniqueness and regularity}
The existence and uniqueness of minimizers of \eqref{eq:minimization} is proved in a standard
way provided that the coercivity of the bilinear form $a(\cdot,\cdot)$ can be shown.
The coercivity proof, using Korn's inequality, is essentially contained in \cite{DuvLio76} as we now explain.
\begin{thm}
    The bilinear form $a(\cdot,\,\cdot)$ is coercive on $V$ with respect to the norm $||\cdot||_1$ induced from
$H^1(\Omega_F)$.  In particular, $a(\cdot,\cdot)$ defines an equivalent inner product on $V$.
\end{thm}

\begin{proof}
Coercivity of $a(\cdot,\,\cdot)$ relies in an essential way on Korn's inequality:
\begin{equation}
        a(\Bu,\,\Bu) + ||\Bu||^2 > c\,||\Bu||^2_1,
        \label{eq:Korn}
\end{equation}
for some $c > 0$ (here $||\cdot||$ denotes the $L_2$ norm).
The proof of \eqref{eq:Korn} found in \cite{DuvLio76} applies to the case for any subspace
$U \subset H^1(\Omega_F))$ consisting of functions with a zero trace on a part of the boundary with
nonzero two-dimensional measure.  This applies to $V$ as its elements vanish on $\partial \Omega$
 -- the no-slip boundary conditions on the outer boundary of $\Omega_F$.
In particular, $a(\cdot,\,\cdot)$ is nondegenerate, since the nontrivial kernel of $D(\Bu)$, consisting
of the rigid motions $\Bu(\Bx) = \Bu_0 + \Bomega_0 \times \Bx$, is excluded from $V$
due these boundary conditions.
The result \eqref{eq:Korn} is nontrivial, since the left-hand side contains only symmetric
combinations of the derivatives of $\Bu$.

The coercivity proof is completed by showing the existence of the following bound:
\begin{equation}
      a(\Bu,\,\Bu) > d ||\Bu||^2,
      \label{eq:pseudo_Poincare}
\end{equation}
for some $d > 0$.  This replaces Poincare's inequality in the case of the symmetrized
gradient.  It can be proved for $V$ as is done in \cite{DuvLio76},
using the compactness of the embedding $V \hookrightarrow L_2(\Omega_F)$.
This embedding is induced from the usual compact embedding $H^1(\Omega_F) \hookrightarrow L_2(\Omega_F)$,
since $V$, being a closed subspace of $H^1(\Omega_F)$ is also weakly
closed (see, e.g., \cite{Lax02}).
\end{proof}

With the coercivity of $a(\cdot,\,\cdot)$ proved, the existence of minimizers for \eqref{eq:minimization}
can be proved by standard techniques.  Since each minimizer satisfies \eqref{eq:variational 1},
the difference of any two of them is $a(\cdot,\cdot)$-orthogonal to a dense subset of $V$, hence is zero,
which proves uniqueness.

Finally, the unique field $\Bu$ that solves \eqref{eq:minimization} is a weak solution of the Stokes equation
on a regular bounded domain.  Therefore, once again by the standard theory  (e.g., \cite{Galdi94})
there exists a unique pressure field $p \in L_2(\Omega_F)$, which together with $\Bu$
satisfies the a priori $L_2$ estimates \cite{Galdi94}.  Since the boundary of $\Omega_F$ and the
righ-hand side of (\stokes) are smooth, these estamates imply that $(\Bu,\, p)$ are smooth too.
By reversing the steps leading to the weak formulation \eqref{eq:variational 1}, we now see that $(\Bu,\,p)$ form a strong solution of
(\full_system).

\subsection{Point forces}
The limit of point forces $\delta \rightarrow 0$ is useful because it simplifies
many concrete calculations in the asymptotic analysis of the model.  Heuristically,
smooth forces can be replace by point forces because the fluid velocity in $B_P$ is ill-defined anyway,
being a simplified representation of a complicated periodic action of the flaggelum.
Therefore, the precise value of the velocity and pressure near the propeller $B_P$ do not matter
and can be left undefined.
At the same time, away from the propeller $B_P$, both a smooth force density in $\mathcal{D}(B_P)$
and a point force density produce comparable results, as will be shown below.

The case of point forces, however, does not fit into the existence proof of the previous subsection
because $\BF$ is no longer in $H^1(\Omega_F)$.  In this case, however, there is still a unique
solution to (\full_system), regular away from $\Bx_P$, which can be constructed
with the help of the Green's function for $\Omega_F$.

Assume for the moment that there exists a unique $\Bu_0$, the solution to \eqref{eq:stokes}
with homogeneous boundary conditions:
\begin{equation}
    \label{eq:homoBC}
    \Bu(\Bx) = 0,\qquad \Bx \in \partial \Omega_F.
\end{equation}
Then, the existence and uniqueness of solution $\Bu$ to (\full_system) is equivalent
to the existence and uniqueness of $\Bu_1 = \Bu - \Bu_0$, the solution to (\full_system)
with $\BF = 0$ and the balance conditions modified to account the forces and torques due
to the point force flow $\Bu_0$.
Specifically, the balance conditions (\balance) are replaced with the following:
\begin{eqnarray}
    \BF^i_{_{H}}+\BF^i_{_{T}} + \BF^i_{0,_H} + \BF^i_{0,_T} + \BF^i_{_{P}} = 0,&\qquad &\text{\emph{balance of forces}},\label{eq:new_force_balance}\\
    \BT^i_{_{H}}+\BT^i_{_{T}} + \BT^i_{0,_H} + \BT^i_{0,_T} = 0, &\qquad &\text{\emph{balance of torques}}\label{eq:new_torque_balance},
\end{eqnarray}
\def\new_balance{\ref{eq:new_force_balance} - \ref{eq:new_torque_balance}}
where $\BF^i_{0,*}$ and $\BF^i_{0,*}$ are the hydrodynamics forces and torques on
the balls $B^i_*,\ * = H,T$ due to a constant flow $\Bu_0$ and computed using
$\Bu_0$ in place of $\Bu$ in \eqref{eq:def:forcesNtorques}.
Now the method of the previous subsection applies to the system satisfied by $\Bu_1$ with the only
modification: the linear functional $b(\cdot)$ defined on $V$ by \eqref{eq:def:L}
is replaced by $L_1(\cdot)$:
\begin{equation}
    b_1(\Bw):= - \left(\BF_{_P} + \BF_{0,_H} + \BF_{0,_T}\right)\cdot\Bv^\Bw_{_C} -
        \left(\BT^i_{0,_H} + \BT^i_{0,_T}\right) \cdot \omega^\Bw.
    \label{eq:def:L1}
\end{equation}
To complete the proof, it remains to show the existence of $\Bu_0$, which is done in \ref{app:Greens_function}.


\subsection{Green's function for $\Omega_F$}
\label{app:Greens_function}
Here we briefly indicate how to show the existence of the Green's function (the Green's tensor
for $\Omega_F$).  The result is well-known and we include it for the sake of completeness.
The Green's velocity tensor $\mathcal{G}$ and the corresponding pressure tensor $\mathcal{P}$
are analogous to the $(G,\,P)$-pair of (\ref{app:point_force}) in that they solve
\begin{equation}
\left\{
    \begin{array}{l}
      \mu \lap \mathcal{G}(\cdot - \Bx_0) = \D \mathcal{P}(\cdot - \Bx_0) - \delta(\cdot - \Bx_0)\BI \\
      \diiv(\Bu)(\cdot - \Bx_0) = 0
    \end{array}
    \right..
 \label{eq:stokes_green}
\end{equation}
with $\Bx$ and $\Bx_0$ in $\Omega_F$ and subject to the homogeneous boundary conditions on $\Omega_F$.
Once the existence of $(\mathcal{G},\,\mathcal{P})$ has been shown, the existence of $\Bu_0$
used in \ref{app:existence} is trivially established:
$$
        \Bu_0(\Bx) = f_P \sum_i \mathcal{G}(\Bx - \Bx^i_P) \cdot \Btau^i.
$$

The sought for Green's tensors are constructed by canceling the boundary values of $G$ and $P$ on
$\partial \Omega_F$ as follows:
$$
        \mathcal{G} = G - \tilde G,\quad
        \mathcal{P} = P - \tilde P,
$$
where $\tilde G$ and $\tilde P$ solve \eqref{eq:stokes_green} with the zero right-hand side
and the boundary conditions
$$
        \tilde G(\cdot - \Bx_0)|_{\partial \Omega_F} = G(\cdot - \Bx_0)|_{\partial \Omega_F},\quad
        \tilde P(\cdot - \Bx_0)|_{\partial \Omega_F} = P(\cdot - \Bx_0)|_{\partial \Omega_F}.
$$
The existence and uniqueness of $(\tilde G,\,\tilde P)$ is easily established by standard methods
(e.g., \cite{Galdi94}) both in the case of a bounded $\Omega_F$ and the exterior $\Omega_F = \R^3 \setminus \Omega_B$.
The only requirement in the bounded case is the compatibility condition
$$
    \int_{\partial \Omega_F} G(\Bx - \Bx_0) dS(\Bx) \cdot n(\Bx) = 0.
$$
This equality easily follows from the divergence theorem applied to $G$, whose divergence is
zero in $L_1(\Omega_F)$, as the following simple calculation shows (summation on $j$ implied
and $|\Bx|^2 = x_j x_j$):
\begin{gather*}
        G_{ij,j} = -\frac{x_j}{|\Bx|} \left( \delta_{ij} + \frac{x_i x_j}{|\Bx|^3}\right)
        - \frac{1}{|\Bx|} \left(\frac{\delta_{ij} x_j + \delta_{jj} x_i}{|\Bx|^3}
                                - 3 \frac{x_i x_j x_j}{|\Bx|^5}
                          \right) =
        -\frac{1}{|\Bx|^3} \left( x_i + x_i \frac{x_j\,x_j}{|\Bx|^3}
                  - 4 \frac{x_i}{|\Bx|} + 3 \frac{x_i}{|\Bx|}\right) = - \frac{x_i}{|\Bx|^3},\\
        |G_{ij,j}| < C\, \frac{1}{|\Bx|^2},\quad C = const.
\end{gather*}

\section{Asymptotic formulas}\label{sect:formulas}

The velocities of the head and tail balls in terms of $\{\a^j\}$ and $\Bx^j_*$
are
\begin{equation}\label{eq:vH}
\begin{split}
    \frac{1}{f_p}\Bv_{_H}^i =
    \sum_{j\neq i}&\Big[
    (1-\a^j)G(\Bx_{_H}^i-\Bx_{_H}^j)+
    \a^j    G(\Bx_{_H}^i-\Bx_{_T}^j)-
            G(\Bx_{_H}^i-\Bx_{_P}^j)
    \Big]\tau^j+\\
    +
    &\Big[
    \a^i    G(\Bx_{_H}^i-\Bx_{_T}^i)+
            G(\Bx_{_H}^i-\Bx_{_P}^i)+
    \frac{1}{\gamma_0}I
    \Big]\tau^i,
\end{split}
\end{equation}
\begin{equation}\label{eq:vT}
\begin{split}
    \frac{1}{f_p}\Bv_{_T}^i =
    \sum_{j\neq i}&\Big[
    (1-\a^j)G(\Bx_{_T}^i-\Bx_{_H}^j)+
    \a^j    G(\Bx_{_T}^i-\Bx_{_T}^j)-
            G(\Bx_{_T}^i-\Bx_{_P}^j)
    \Big]\tau^j+\\
    +
    &\Big[
    (1-\a^i)G(\Bx_{_T}^i-\Bx_{_H}^i)+
            G(\Bx_{_T}^i-\Bx_{_P}^i)+
    \frac{1}{\gamma_0}I
    \Big]\tau^i.
\end{split}
\end{equation}
Note that
\begin{equation}\label{eq:Gproperty}
    G(-\Bx)=G(\Bx)\qquad \forall\Bx\in\R^3
\end{equation}
and
\begin{equation}\label{eq:tildeGproperty}
    \tilde{G}(-\Bx)=\tilde{G}(\Bx)\qquad \forall\Bx\in\R^3.
\end{equation}
Hence,
\begin{equation}\label{eq:1}
    G(\Bx_{_T}^i-\Bx_{_H}^i) = G(\Bx_{_H}^i-\Bx_{_T}^i)=G(2L\tau).
\end{equation}

Using the (\ref{eq:vH},\ref{eq:vT}) and the relation
\begin{eqnarray*}
  \Bv^i_{_C}    &=& \frac{1}{2}(\Bv^i_{_H}+\Bv^i_{_T}), \\
  \w^i          &=& \frac{1}{2L}(\Bv^i_{_H}-\Bv^i_{_T})\times \tau^i,
\end{eqnarray*}
we obtain
\begin{eqnarray}\label{eq:vC}
    \frac{2}{f_p}\Bv_{_C}^i & = &
    \sum_{j\neq i}\Big[
    (1-\a^j) \left\{G(\Bx_{_H}^i-\Bx_{_H}^j)+G(\Bx_{_T}^i-\Bx_{_H}^j)\right\}+\nonumber\\
    & & +\a^j    \left\{G(\Bx_{_H}^i-\Bx_{_T}^j)+G(\Bx_{_T}^i-\Bx_{_T}^j)\right\}-\nonumber\\
    & & -\left\{G(\Bx_{_H}^i-\Bx_{_P}^j)+G(\Bx_{_T}^i-\Bx_{_P}^j)\right\}
    \Big]\tau^j+\nonumber\\
    & & +
    \Big[
            G(2L\tau^i)+
            G(\Bx_{_H}^i-\Bx_{_P}^i)+G(\Bx_{_T}^i-\Bx_{_P}^i)+
            \frac{2}{\gamma_0}I
    \Big]\tau^i,
\end{eqnarray}
\begin{eqnarray}\label{eq:vC:1}
    \frac{2L}{f_p}\w^i & = &
    \tau^i\times\sum_{j\neq i}\Big[
    (1-\a^j)\left\{G(\Bx_{_H}^i-\Bx_{_H}^j)-G(\Bx_{_T}^i-\Bx_{_H}^j)\right\}+\nonumber\\
    & & +\a^j   \left\{G(\Bx_{_H}^i-\Bx_{_T}^j)-G(\Bx_{_T}^i-\Bx_{_T}^j)\right\}-\nonumber\\
    & & -       \left\{G(\Bx_{_H}^i-\Bx_{_P}^j)-G(\Bx_{_T}^i-\Bx_{_P}^j)\right\}
    \Big]\tau^j+\nonumber\\
    & & +
    \tau^i\times\Big[
    (2\a^i-1)G(2L\tau^i)+
            G(\Bx_{_H}^i-\Bx_{_P}^i)-G(\Bx_{_T}^i-\Bx_{_P}^i)
    \Big]\tau^i.
\end{eqnarray}

\subsection{Expansion of $G(\cdot)$}

The expansion of $\Bv^i_*$ in orders of $\ve$ is due to the expansions of the Green's function $G$, e.g.
\begin{eqnarray}
    \label{eq:sample expansion}
    G(\Bx_{_H}^i-\Bx_{_P}^j) 
    &=& 
    G(\Bx_{_H}^i-\Bx_{_C}^i+\Bx_{_C}^i-\Bx_{_C}^j+\Bx_{_C}^j-\Bx_{_P}^j)=\\
    \nonumber
    &=&
    G\bigg((\Bx_{_H}^i-\Bx_{_C}^i)+(\Bx_{_C}^j-\Bx_{_P}^j)+(\Bx_{_C}^i-\Bx_{_C}^j)\bigg).
\end{eqnarray}
Here the quantities $(\Bx_{_H}^i-\Bx_{_C}^i)$ and $(\Bx_{_C}^j-\Bx_{_P}^j)$ measure the distances in the
same bacteria, hence they do not depend on $\ve=|\Bx_{_C}^i-\Bx_{_C}^j|^{-1}$.
The only quantity that depends on $\ve$ is $(\Bx_{_C}^i-\Bx_{_C}^j)$.

\subsection{Asymptotic expansion for $\a^i$ in powers of $\ve=|\Bx^2_{_C}-\Bx^1_{_C}|$}
\label{subsect:alpha expansion}

Consider the system \eqref{eq:rigidity system} for $\alpha^i$.

Consider the asymptotic expansion
\begin{equation}\label{eq:anzats a}
    \a^i = a^0+\ve \a^i_1 + \ve^2 \a^i_2+\dots,
\end{equation}
and substitute it into
\eqref{eq:rigidity system}, where all terms are expanded in asymptotic series
in $\ve = |\Bx^2_{_C}-\Bx^1_{_C}|^{-1}$.
Note that the terms like $G(\Bx^i_{_H}-\Bx^i_{_T})$ do not depend on $\ve$;
hence, their expansion will have only $\ve^0$ order term (itself).
On the other hand, the terms like $G(\Bx^i_{_H}-\Bx^j_{_T})$ are of order $\ve^1$ and will not have
order $\ve^0$ terms.

Thus, at the order $\ve^0$ the equation \eqref{eq:rigidity system} becomes
\begin{equation}\label{eq:rigidity system order 0}
\begin{split}
    \a^i_0
    (\tau^i)^T\bigg[& G(\Bx^i_{_H}-\Bx^i_{_T})+G(\Bx^i_{_T}-\Bx^i_{_H})-\frac{2}{\gamma_0}\bigg]\tau^i=\\
    =
    (\tau^i)^T\bigg[
    & G(\Bx^i_{_T}-\Bx^i_{_H})+
    G(\Bx^i_{_H}-\Bx^i_{_P})-G(\Bx^i_{_T}-\Bx^i_{_P})
    -\frac{1}{\gamma_0}
    \bigg]\tau^i.
\end{split}
\end{equation}
Express the arguments of $G(\cdot)$ in terms of $\tau^i, L,$ and $\zeta$:
\begin{equation}\label{eq:aux 1}
    \Bx^i_{_H}-\Bx^i_{_T} = 2L\tau^i,
    \qquad
    \Bx^i_{_H}-\Bx^i_{_P} = (1-\zeta)L\tau^i,
    \qquad
    \Bx^i_{_T}-\Bx^i_{_P} = (-1-\zeta)L\tau^i,
\end{equation}
and substitute back into \eqref{eq:rigidity system order 0} to get
\begin{equation}\label{eq:rs o0 1}
\begin{split}
    \a^i_0
    (\tau^i)^T\bigg[& G(2L\tau^i)+G(-2L\tau^i)-\frac{2}{\gamma_0}\bigg]\tau^i=\\
    =
    (\tau^i)^T\bigg[
    & G(-2L\tau^i)+
    G((1-\zeta)L\tau^i)
    -G((-1-\zeta)L\tau^i)
    -\frac{1}{\gamma_0}
    \bigg]\tau^i.
\end{split}
\end{equation}

Using the definition of
\begin{equation}\label{eq:gama0}
    \gamma_0 = \frac{1}{8\pi\mu R}
\end{equation}
and the properties (\ref{eq:lem:prop 1}-\ref{eq:lem:prop 2}), simplify \eqref{eq:rs o0 1}:
\begin{equation}\label{eq:rs o0 2}
\begin{split}
    \frac{1}{4\pi\mu} \a^i_0
    \bigg[
    \frac{1}{ 2L}
    +
    \frac{1}{2L}-\frac{1}{R}\bigg]
    =
    \frac{1}{4\pi\mu}
    \bigg[
    \frac{1}{ 2L}+
    \frac{1}{|1-\zeta|L}
    -
    \frac{1}{|1+\zeta|L}
    -\frac{1}{2R}
    \bigg].
\end{split}
\end{equation}
Multiply through by $4\pi\mu R L$:
\begin{equation}\label{eq:rs o0 3}
\begin{split}
    -\a^i_0
    \bigg[
    1-\frac{R}{L}\bigg]
    =
    -\bigg[
    \frac{1}{2}
    -\frac{R}{ 2L}-
    \frac{R}{|1-\zeta|L}
    +
    \frac{R}{|1+\zeta|L}
    \bigg].
\end{split}
\end{equation}
Pull out $\frac{1}{2}$ from the RHS, and solve for $\a^i_0$
\begin{equation}\label{eq:rs o0 4}
\begin{split}
    \a^i_0
    =
    \frac{1}{2}\
    \frac{
    1
    -
    \frac{R}{L}
    \left(
    1+
    \frac{2}{|1-\zeta|}
    -
    \frac{2}{|1+\zeta|}
    \right)}
    {
    1-\frac{R}{L}
    }.
\end{split}
\end{equation}
Denote
\begin{equation}\label{eq:xi}
    \xi:=\frac{R}{L}\ll 1,
\end{equation}
and perform the expansion of \eqref{eq:rs o0 4} in terms of $\xi$
\begin{equation}\label{eq:rs o0 5}
\begin{split}
    \a^i_0
    =
    \frac{1}{2}\
    \bigg(
    1+z\xi+z\xi^2+\dots
    \bigg),
\end{split}
\end{equation}
where
\begin{equation}\label{eq:z}
    z=z(\zeta):= 2\left(\frac{1}{|1+\zeta|}-\frac{1}{|1-\zeta|}\right)
    =
    \left\{
    \begin{split}
    \frac{4}{\zeta^2-1}         \qquad &\text{if }\zeta<-1,\\
    \frac{4\zeta}{\zeta^2-1}    \qquad &\text{if }1<\zeta<1,\\
    \frac{-4\zeta}{\zeta^2-1}   \qquad &\text{if }1<\zeta.\\
    \end{split}
    \right.
\end{equation}

\begin{lem}[Properties of $G(\cdot)$]
    Note the following two properties of $G(\cdot)$:
    \begin{enumerate}
      \item Let $q\in\R$ and $\tau\in\R^3$, $|\tau|=1$. Then
        \begin{equation}\label{eq:lem:prop 1}
            G(q\tau) = \frac{1}{|q|} G(\tau).
        \end{equation}
      \item Let $\tau\in\R^3$, $|\tau|=1$. Then
        \begin{equation}\label{eq:lem:prop 2}
            \tau^T G(\tau) \tau = \frac{1}{4\pi\mu}.
        \end{equation}
    \end{enumerate}
\end{lem}

\begin{proof}

Property 1 follows simply from the definition of $G(\cdot)$:
\begin{equation}\label{eq:proof:prop 1}
    G(q\tau) =\frac{1}{8\pi \mu |q|}\left(\BI+\tau\tau^T\right)=\frac{1}{|q|}G(\tau).
\end{equation}

Property 2 follows by simple substitution:
\begin{equation}\label{eq:proof:prop 1:1}
    \tau^T G(\tau)\tau
    =
    \frac{1}{8\pi \mu}\tau^T\left(\BI+\tau\tau^T\right)\tau
    =
    \frac{1}{8\pi \mu}(1+1)=
    \frac{1}{4\pi \mu}.
\end{equation}
\end{proof}

\subsection{Sign of $A^i(f_p,L,R,\mu,\a_0)$}
\label{subsect:A sign}

\begin{lem}[Properties of $(1-\zeta-2\a_0^i)$]
    Assume $f_p>0$ in the expression
    \begin{equation}\label{eq:A 1}
        A^i(f_p,L,R,\mu,\a_0)=\frac{f_p L}{32 \pi\mu} (1-\zeta-2\alpha_0),
    \end{equation}
    where $\a_0$ is given by (\ref{eq:rs o0 5}-\ref{eq:z}).
    For pushers ($\zeta<0$)
    \begin{equation}\label{eq:pusher}
        A^i(f_p,L,R,\mu,\a_0)>0.
    \end{equation}
    For pullers ($\zeta>0$)
    \begin{equation}\label{eq:puller}
        A^i(f_p,L,R,\mu,\a_0)<0.
    \end{equation}
\end{lem}

\begin{proof}
Use the formula \eqref{eq:rs o0 5} for $\a^i_0$ to rewrite
\begin{equation}\label{eq:prop:2}
\begin{split}
    1-\zeta-2\a_0^i &=
    1-\zeta-\bigg(
    1+z\xi+z\xi^2+\dots
    \bigg)=
    -\bigg(\zeta+
    z\xi+z\xi^2+\dots
    \bigg)=\\
    &=-\left(\zeta +\frac{z\xi}{1-\xi}\right).
\end{split}
\end{equation}

Consider three case: $\zeta<-1$, $-1<\zeta<1$, $1<\zeta$.
For each of these cases use the formula \eqref{eq:z} for $z(\zeta)$ to evaluate
\eqref{eq:prop:2}.

Case \underline{$\zeta<-1$}: Here $z(\zeta)=\frac{4}{\zeta^2-1}$, and
\begin{equation}\label{eq:prop:3}
    -\left(\zeta +\frac{z\xi}{1-\xi}\right)=
    -\left(\zeta +\frac{4}{\zeta^2-1}\ \frac{\xi}{1-\xi}\right).
\end{equation}
This expression is always positive when $(-\zeta-1)\sim 1$ as $\xi\to 0$.

Case \underline{$-1<\zeta<1$}: Here $z(\zeta)=\frac{4\zeta}{\zeta^2-1}$, and
\begin{equation}\label{eq:prop:4}
    -\left(\zeta +\frac{z\xi}{1-\xi}\right)=
    -\left(\zeta +\frac{4\zeta}{\zeta^2-1}\ \frac{\xi}{1-\xi}\right)=
    -\zeta\left(1 +\frac{4}{\zeta^2-1}\ \frac{\xi}{1-\xi}\right).
\end{equation}
This expression changes sign from positive to negative only as $\zeta$ passes through 0.
So it is positive when $\zeta<0$ and negative when $\zeta>0$.

Case \underline{$1<\zeta$}: Here $z(\zeta)=\frac{-4\zeta}{\zeta^2-1}$ and
\begin{equation}\label{eq:prop:5}
    -\left(\zeta +\frac{z\xi}{1-\xi}\right)=
    -\left(\zeta -\frac{4\zeta}{\zeta^2-1}\ \frac{\xi}{1-\xi}\right).
\end{equation}
This expression is always negative when $(\zeta-1)\sim 1$ as $\xi\to 0$.
\end{proof}

%

\section{Stability of the ``mirror image'' configuration}\label{subsect:stability of MI}

Before analyzing the stability of the ``mirror image'' configuration we determine the quantity,
which does not depend of the orientation of the $\overrightarrow{ox}$-axis 
(i.e., it does not depend on the choice of observer),
that characterizes how close a given configuration is to a ``mirror image'' configuration.
Then we perturb this parameter by a small amount and check whether this parameter 
is decreasing. 
If this parameter is decreasing, the configuration is stable; otherwise it is unstable.

Note that for two swimmer, in the ``mirror image'' configuration
\begin{equation}\label{eq:characerization of mirror image}
    \pi + 2\phi-(\theta^1+\theta^2) = 2\pi n,\qquad n\in \mathbb{N}
\end{equation}
for any choice of the $\overrightarrow{ox}$-axis.
Moreover, if equation \eqref{eq:characerization of mirror image} holds, then two swimmers are in the ``mirror image'' configuration.
Therefore, the quantity
\begin{equation}\label{eq:perturbation of the mirror image}
    \delta := \pi + 2\phi-(\theta^1+\theta^2)
\end{equation}
can be viewed as a measure of deviation from the ``mirror image'' configuration.

To perform the stability analysis, we perturb the ``mirror image'' configuration.
That is, we choose $|\delta(0)| > 0$ small, and we check whether $|\delta(t)|$ decreases with time.
If $\delta^\prime(0)$ has the opposite sign to $\delta(0)$, then $|\delta(t)|$
decreases with time locally, and the configuration is stable; otherwise it is not stable.

We have
\begin{equation}\label{eq:dt delta}
    \delta^\prime
    =
    2\phi^\prime-({\theta^1}^\prime+{\theta^2}^\prime)
    =
    2\phi^\prime-(\w^1+\w^2)
    .
\end{equation}
The expressions for $\w^1$ and $\w^2$ can be found
from (\ref{eq:order 3},\ref{eq:A=Algebraic},\ref{eq:C=Trigonometric}).
The expression for $\phi^\prime$ can be found simply
as a projection of the translational velocity difference $(\Bv^2_{_C}-\Bv^1_{_C})$
onto the unit circle
\begin{equation*}
\begin{split}
    \phi^\prime
    &=
    \ve (\Bv^2_{_C}-\Bv^1_{_C})\cdot
    \left[
      \begin{array}{c}
        -\sin(\phi) \\
        \cos(\phi) \\
      \end{array}
    \right]
    =\\
    &=
    \ve v_0
    \left(
    -\sin(\phi) (\cos(\theta^2)-\cos(\theta^1))
    +\cos(\phi) (\sin(\theta^2)-\sin(\theta^1))
    \right)+ O(\ve^3),
\end{split}
\end{equation*}
where we used
\begin{equation*}
    \Bv^i_{_C} = v_0\tau^i+O(\ve^2), \qquad i=1,2.
\end{equation*}

Without loss of generality, choose the ${x}$-axis so that $\pi + 2\phi=0$, that is,
$\phi = -\pi/2$, which means the second swimmer is directly \emph{below} the first swimmer.
Then
\begin{equation}\label{eq:dt delta 2}
\begin{split}
    \delta^\prime
    &=
    \ve v_0 \left(\cos(\theta^2)-\cos(\theta^1)\right) + O(\ve^3)
    =
    -2 \ve v_0
    \sin\left(\frac{\theta^2-\theta^1}{2}\right)
    \sin\left(\frac{\theta^2+\theta^1}{2}\right) + O(\ve^3)
    =\\
    &=
    -2 \ve v_0
    \sin\left(\theta^1-\frac{\delta}{2}\right)
    \sin\left(\frac{\delta}{2}\right) + O(\ve^3).
\end{split}
\end{equation}
Thus,
for swimmers rotated outward ($0<\theta^1<\pi$) the ``mirror image'' configuration is stable and
for swimmers rotated inward  ($0>\theta^1>-\pi$) the ``mirror image'' configuration is unstable under small perturbations.
The results of the stability analysis are not affected by the type of the swimmer;
they are the same for all values of $\zeta$.

\bibliographystyle{ieeetr}
\bibliography{dumbbell}

\end{document}